\RequirePackage{fix-cm}
\documentclass[smallextended]{svjour3}       
\smartqed  
\usepackage{graphicx}
\usepackage{ subfigure }
\usepackage{url}
\usepackage{xcolor}
\usepackage{rotating}
\usepackage{pdflscape}
\definecolor{newcolor}{rgb}{.8,.349,.1}

\usepackage{amsmath}	
\usepackage{subfig}
\usepackage{amsmath}
\usepackage{amsfonts}

\usepackage{url}
\usepackage{multirow}

\newcommand{\OLdue}[0]{$\!\!\!\!\!\!\!\!\!\!\!\mathcal{O}(L_2)\!\!\!\!\!\!$}
\newcommand{\CPU}[0]{$\!\!\!\!\!\!\!\!\!\!\!\!$T$\!\!\!\!\!\!\!\!\!\!\!$}

\newcommand{\mcolt}[2]{\multicolumn{3}{|c|}{$\mathbf{P_{#1}P_{#2}}$}}

\newcommand{\sm}{$\!\!$}
\newcommand{\smm}{$\!\!\!$}
\newcommand{\smmm}{$\!\!\!\!$}

\newcommand{\mbf}[1]{\mathbf{#1}}			%
\newcommand{\x}{\mbf{x}}

\newcommand{\Omegai}{\Omega_{\mathbf{i}}}
\newcommand{\I}{\mathbf{i}}

\newcommand{\omegaialpha}{\omega_{\mathbf{i},\alpha}}

\newcommand{\Q}{\mathbf{Q}}

\newcommand{\0}{\mathbf{0}}
\renewcommand{\u}{\mathbf{u}}
\newcommand{\w}{\mathbf{w}}
\newcommand{\q}{\mathbf{q}}
\newcommand{\F}{\mathbf{F}}

\newcommand{\f}{\mathbf{f}}
\newcommand{\g}{\mathbf{g}}
\newcommand{\h}{\mathbf{h}}
\newcommand{\halb}{\frac{1}{2}} 
 
\renewcommand{\v}{\mathbf{v}}
\newcommand{\B}{\mathbf{B}}

\newcommand{\be}{\begin{equation} \begin{aligned} }
\newcommand{\ee}{\end{aligned} \end{equation}}

\newcommand{\PN}[2]{P_#1P_#2}
\newcommand{\BoldXi}{\boldsymbol{\xi}}

\newcommand\Tstrut{\rule{0pt}{2.6ex}}         
\newcommand\Bstrut{\rule[-0.9ex]{0pt}{0pt}}   

\newcommand{\lmax}{\ell_{\max}}
\newcommand{\rRef}{\mathfrak{r}}

\renewcommand{\epsilon}{\varepsilon}
\renewcommand{\phi}{\varphi}

\begin{document}

\title{\textit{A posteriori} subcell finite volume limiter for general $P_NP_M$ schemes: applications from gasdynamics to relativistic magnetohydrodynamics}

\titlerunning{\textit{A posteriori} subcell FV limiter for $P_NP_M$ schemes}        

\author{Elena Gaburro \and Michael Dumbser}


\institute{
E. Gaburro \at
Department of Civil, Environmental and Mechanical Engineering, University of Trento, Via Mesiano 77, 38123 Trento, Italy \\
\email{elena.gaburro@unitn.it}         
\and
M. Dumbser \at
Department of Civil, Environmental and Mechanical Engineering, University of Trento, Via Mesiano 77, 38123 Trento, Italy 
}

\date{Received: date / Accepted: date}

\maketitle

\begin{abstract} 
In this work, we consider the general family of the so called ADER $\PN{N}{M}$ schemes for the numerical solution of hyperbolic partial differential equations with \textit{arbitrary} high order of accuracy in space and time.

The family of one-step $\PN{N}{M}$ schemes was introduced in~\cite{Dumbser2008} and represents a unified framework
for classical high order Finite Volume (FV) schemes ($N=0$), the usual Discontinuous Galerkin (DG) methods ($N=M$),
as well as a new class of intermediate hybrid schemes for which a reconstruction operator of degree $M$ is applied 
over piecewise polynomial data of degree $N$ with $M>N$.  
In all cases with $M \geq N > 0 $ the $\PN{N}{M}$ schemes are \textit{linear} in the sense of  Godunov~\cite{godunov}, thus when considering phenomena characterized by discontinuities, spurious oscillations may appear and even destroy the simulation.  
Therefore, in this paper we present a new simple, robust and accurate \textit{a posteriori} subcell finite volume limiting strategy that is valid for the entire class of $\PN{N}{M}$ schemes. The subcell FV limiter is activated only where it is needed, i.e. in the neighborhood of shocks or other discontinuities, and is able to maintain the resolution of the underlying high order $\PN{N}{M}$ schemes, due to the use of a rather fine subgrid of $2N+1$ subcells per space dimension.  

The paper contains a wide set of test cases for different hyperbolic PDE systems, solved on adaptive Cartesian meshes (AMR) that show the capabilities of the proposed method both on smooth and discontinuous problems, 
as well as the broad range of its applicability. The tests range from compressible gasdynamics over classical MHD to relativistic magnetohydrodynamics.

\end{abstract}

\section{Introduction}
\label{intro}

In this work we want to improve the family of high order accurate ADER $\PN{N}{M}$ schemes 
first introduced in~\cite{Dumbser2008,ADERNSE} for the solution of hyperbolic partial differential equations. 
In this family of schemes the discrete solution is represented in space through  
high order \textit{piecewise polynomials} of degree $N$ at each timestep; 
the data are then evolved in time through a space-time \textit{reconstruction} procedure of order $M$.
The reconstruction procedure is divided into two steps: 
concerning the \textit{spatial} reconstruction, we employ a classical WENO reconstruction in the case of pure finite volume schemes ($N=0$), 
a reconstruction procedure based on $L^2$ projection that is linear in the sense of Godunov for $N>0$ and $M>N$,  and in the case or pure DG schemes ($N=M$) the reconstruction reduces to the identity operator;
concerning the reconstruction in \textit{time}, we employ a novel variant of the ADER approach of Toro and Titarev, see \cite{toro3,toro4,titarevtoro,Toro:2006a,BTVC2016},  
based on an element-local space-time Galerkin predictor, see \cite{Dumbser2008}.
In practice, we can see the Finite Volume (FV) schemes of order $M$ as a particular case of $\PN{N}{M}$ methods
when $N=0$, and also the Discontinous Galerkin (DG) methods are included in this family when choosing $N=M$.

Furthermore, this family contains another important class of \textit{hybrid} or \textit{reconstructed} DG schemes when taking $N>0, M>N$, which are the main object of study of this paper. 
Indeed, they offer many advantages, in particular their good \textit{compromise} between \textit{cost and resolution}. 
In fact, data are represented with polynomials of order $N$, 
so more accurately with respect to FV methods, 
but without the expensive cost of a full DG representation of approximation degree $M$;  
also the CFL stability constraint that limits the timestep size of any explicit scheme, 
only depends on $N$ and not on $M$, allowing for larger timesteps once the desired order of accuracy has been fixed, see \cite{Dumbser2008}. Last but not least, for $N>0$ the  $\PN{N}{M}$ schemes require a \textit{much smaller reconstruction stencil} than comparable finite volume schemes of degree $M$. 
The nominal order of accuracy of the scheme is given by $M+1$ so it can be at least in principle arbitrary high.

The family of \textit{reconstructed DG schemes}, which is similar to the $\PN{N}{M}$ framework, 
was forwarded independently by Luo et al. in a series of papers, see e.g. \cite{luo1,luo2,luo3,luo4,luo5,luo6} and references therein.  At this point we also highlight that the use of reconstruction and filtering operators as a \textit{post-processor} for improving the accuracy of DG schemes goes back to work of Ryan et al., see \cite{ryan1,ryan2,ryan3,ryan4,ryan5}.
Other related work on reconstruction-based DG schemes can be found in  \cite{vanLeerDGdiffusion,chiravalle20193d,wang2020reconstructed,WAO-ALE}.   

Moreover, the ADER $\PN{N}{M}$ family provides a useful framework for code developers 
because it allows to include in a \textit{unique code} 
both types of standard discretization methods for hyperbolic PDE (FV and DG schemes), together with the new class of intermediate hybrid schemes for $M>N>0$. It is then possible to let it up to the user to decide whether for a particular application the use of
a robust finite volume approach ($N=0$), a very accurate DG scheme ($N=M$), or a less expensive but still very accurate intermediate $\PN{N}{M}$ method with $M>N>0$, is the most appropriate. 

The main drawback so far of the intermediate $\PN{N}{M}$ schemes with $M>N>0$, as presented in~\cite{Dumbser2008}, 
is that they are \textit{linear} in the sense of Godunov~\cite{godunov}, hence not well suited for dealing with discontinuous problems. 
For this reason, here we propose a new simple, robust and accurate \textit{limiting strategy} that is 
able to stabilize the \textit{entire class} of $\PN{N}{M}$ schemes in such a way that they can be 
 employed for the numerical solution of hyperbolic equations with discontinuous solutions, which may arise even when starting with smooth initial conditions. Moreover, the new limiter does not substantially deteriorate the benefits of $\PN{N}{M}$ schemes in terms of computational cost and accuracy of the original unlimited schemes.
To the best knowledge of the authors, this is the first time that an \textit{a posteriori} subcell finite volume limiter is proposed for general $\PN{N}{M}$ schemes with $M>N>0$. So far, only the cases $N=0$ and $M=N$ were covered in \cite{Dumbser2008} and \cite{DGLimiter1}, respectively. 
 
Our limiter is based on the MOOD approach~\cite{CDL1,CDL2,CDL3}, 
which has already been successfully applied in the framework of ADER finite volume schemes \cite{ADERMOOD,ALEMOOD1,ALEMOOD2}
and Discontinous Galerkin finite element schemes, see \cite{DGLimiter1,DGLimiter2,DGLimiter3,FrontierADERGPR,SolidBodies2020}.
Specifically, the numerical solution is checked \textit{a posteriori} for nonphysical values and spurious oscillations,
and if it does not satisfy all admissibility detection criteria, given by both physical and numerical requirements,
in a certain cell, that cell is marked as \textit{troubled}. 
Then, instead of applying a limiter to the already computed solution, 
the solution is {\em locally recomputed} with a more robust scheme in the {troubled cells}, 
relying either on a second order TVD scheme, as proposed for pure DG schemes in~\cite{DGLimiter3,ALEDG,SonntagDG},
or on a higher order ADER-WENO finite volume method as employed in  \cite{DGLimiter1,DGLimiter2,DGCWENO,rannabauer2018ader,DeLaRosaMunzDGMHD}. 
Moreover, this second computation is performed on a finer subgrid generated within each troubled cell;
the subcell approach is employed in order to maintain the high resolution of the initial $\PN{N}{M}$ scheme even when passing to a less accurate 
(but more robust) FV scheme.
For the given reasons our limiter is called \textit{a posteriori} subcell finite volume limiter.

Finally, for a complete review of ADER $\PN{N}{M}$ schemes we refer to the recent paper~\cite{FrontierADERGPR}, 
where a complete introduction traces the historical developments of these methods up to its latest evolutions.


The rest of the paper is organized as follows. 
After an introduction of the class of physical phenomena that can be discretized with the proposed numerical method and the structure of our data representation, we present the family of ADER $\PN{N}{M}$ schemes in Section~\ref{sec.method}. 
In particular, we describe the reconstruction procedure in space, see Section~\ref{ssec.reconstruction}, and in time see Section~\ref{ssec.predictor};
these procedures provide a high order reconstructed polynomial of degree $M$ in space and time 
that will be used in the final one-step update formula given in Section~\ref{ssec.PNPM}.
Then, Section~\ref{ssec.limiter} is dedicated to our \textit{a posteriori} subcell FV limiter, 
which in addition can be combined with mesh adaptation techniques as described in Section~\ref{ssec.AMR}.

Next, in Section~\ref{sec.results} we present a large set of numerical results that shows the order of convergence of our scheme for smooth solutions and their capability of dealing with discontinuities, i.e. their robustness and resolution. 
We also compare the hybrid reconstructed schemes with pure DG schemes in order to show the resulting gain in terms of computational cost.  
Finally, we close the paper with some remarks and an outlook to future works in Section~\ref{sec.conclusion}.

\section{Numerical method}
\label{sec.method}

In this Section we carefully describe the \textit{a posteriori} subcell finite volume limiter for general $P_NP_M$ schemes, showing its simplicity, accuracy, robustness and versatility thanks to the following key ingredients: 
\begin{itemize} 
\item the use of the \textit{unified} $\PN{N}{M}$ framework for finite volume (FV), discontinous Galerkin (DG) and hybrid reconstructed DG schemes allows the user to decide freely which combination of $N$ and $M$ is the better choice for a particular application;
\item the ADER space-time \textit{predictor-corrector} formalism allows the implementation of a truly arbitrary high order accurate fully discrete one-step scheme that needs only one MPI communication per time step within a parallel HPC implementation, see Section~\ref{ssec.predictor};
\item the \textit{a posteriori} subcell finite volume limiter avoids spurious oscillations of high order $P_NP_M$ schemes without affecting the resolution of the underlying method, see Section~\ref{ssec.limiter}; 
\item the adaptive mesh refinement (AMR) technique allows to use a fine grid only where necessary, 
resorting to cheaper coarse grids in smooth regions of the solution, see Section~\ref{ssec.AMR}.
\end{itemize}

\subsection{Governing PDE system}
\label{ssec.PDE}
We consider a very general formulation of the governing equations in order to model a wide class of physical phenomena, namely all those which are described by hyperbolic systems of conservation laws 
that can be cast into the following form, 
\be
\label{eq.generalform}
\partial_t\Q + \nabla \cdot \F(\Q) = 0, \qquad \x \in \Omega(t) \subset \mathbb{R}^d, 
\qquad \Q \in \Omega_{\Q} \subset \mathbb{R}^{m},   
\ee
where $\x = (x,y,z)$ is the spatial position vector, $d$ is the number of space dimensions, $t$ represents the time, 
$\Q = (q_1,q_2, \dots, q_{m})^T$ is the vector of conserved variables defined in the space of the admissible states 
$\Omega_{\Q} \subset \mathbb{R}^{m}$ and $ \F(\Q) = (\,\f(\Q), \g(\Q), \h(\Q)\,) = \f^i (\Q)$ ($i=1,2,3$) is the non-linear flux tensor. 
This kind of system~\eqref{eq.generalform} is said to be \textit{hyperbolic} if for all directions $\mathbf{n} \neq \mathbf{0}$ 
the matrix 
\begin{equation*}
\mathbf{A}_n =  \partial \mathbf{F} / \partial{\mathbf{Q}}  \cdot \mathbf{n}
\label{eqn.A.sys} 
\end{equation*}
has $m$ real eigenvalues and a full set of $m$ linearly independent eigenvectors. 
Examples of hyperbolic equations are the Euler equations of gasdynamics, the Shallow Water equations~\cite{Casulli1990,tavelli2014high} 
and many multiphase models~\cite{BaerNunziato1986,dumbser2013diffuse,gaburro2018diffuse} used in fluid mechanics, 
the magnetohydrodynamics system (MHD) for plasma flow \cite{BalsaraSpicer1999,Balsara2004}, the unified first order hyperbolic formulation of continuum mechanics 
by Godunov, Peshkov and Romenski (GPR)~\cite{GodRom1972,PeshRom2014,GodRom2003,GPRmodel,GPRmodelMHD,DFTBW2018} 
as well as the special and general relativistic formulations of MHD, see e.g. \cite{BalsaraRMHD,RMHD,Banyuls97,Aloy1999c,DelZanna2007,ADERGRMHD}, 
or for the Einstein field equations (CCZ4)~\cite{Alic:2009,Alic:2012,dumbser2020glm,ADERCCZ4}. 
We will test the method proposed in this paper on some of those systems in order to verify its applicability in different physical domains.

\subsection{Domain discretization and high order data representation (order $N$)}
\label{ssec.data}

On grid level $\ell=0$ the computational domain $\Omega$ is discretized with a uniform Cartesian grid, called \textit{main grid} or the \textit{level zero grid}, 
composed of $N_E =N_x \times N_y\times N_z$ conforming elements (quadrilaterals if $d=2$, or hexahedra if $d=3$) 
denoted by $\Omegai = \Omega_{ijk}, \mathbf{i} =(i,j,k)$ with $|\mathbf{i}| = 1, \dots, N_E$, $i = 1, \dots, N_x$, $j=1, \dots, N_y$, $k=1, \dots, N_z$, with volume $|\Omega_{ijk}| = \int_{\Omega_{ijk}} d\x$ and such that
\be
& \Omegai = \Omega_{ijk} \!=\! [x_{i-\frac{1}{2}},x_{i+\frac{1}{2}}]\times [y_{j-\frac{1}{2}},y_{j+\frac{1}{2}}]\times [z_{k-\frac{1}{2}},z_{k+\frac{1}{2}}], \\  
&\text{with } \ 
\Delta x_{i} \!=\! x_{i+\frac{1}{2}}-x_{i-\frac{1}{2}}, \quad 
\Delta y_{j} \!=\! y_{j+\frac{1}{2}}-y_{j-\frac{1}{2}}, \quad 
\Delta z_{k} \!=\! z_{k+\frac{1}{2}}-z_{k-\frac{1}{2}}.
\label{eq.MainGrid}
\ee 
For each element we define a reference frame of coordinates $\BoldXi=(\xi, \eta, \zeta)$ linked to the Cartesian coordinates $\x=(x,y,z)$ of $\Omega_{ijk}$ by
\be
\label{eq.mapping}
x \!=\! x_{i-\frac{1}{2}} + \xi \Delta x, \quad  y\!=\!y_{j-\frac{1}{2}} + \eta \Delta y, \quad
z\!=\!z_{k-\frac{1}{2}} + \zeta \Delta z, \quad \xi, \eta, \zeta \in [0,1].
\ee 

Then, we represent the conserved variables $\Q$ of~\eqref{eq.generalform} in each cell $\Omegai$ by a $d-$dimensional tensor product of piecewise polynomials of degree $N$ 
\be
&\mathbf{u}_h(\x,t^n) = \mathbf{u}_h(\BoldXi(\x)) = \sum \limits_{\ell=0}^{\mathcal{N}-1} \phi_\ell(\BoldXi) \, \hat{\mathbf{u}}_{\ell,\I} 
:= \phi_\ell(\BoldXi) \, \hat{\mathbf{u}}_{\ell,\I}, \\[2pt]
&\x \in \Omegai, \quad \mathcal{N} = (N+1)^d,
\label{eqn.uh}
\ee
where $\phi_\ell(\BoldXi)$ are \textit{nodal} spatial basis functions 
given by the tensor product of a set of Lagrange interpolation polynomials of maximum degree $N$ such that
\be
\label{eq.nodalGeneralBasis}
\phi_\ell(\BoldXi_{\text{GL}}^m) & = \phi_{\ell_1}\!\left(\xi_{\text{GL}}^m\right)\phi_{\ell_2}\!\left(\eta_{\text{GL}}^m\right)\phi_{\ell_3}\!\left(\zeta_{\text{GL}}^m\right) = \delta_{\ell m},
\ee 
where $\BoldXi_{\text{GL}}^m$ are the set of $(N+1)^d$ Gauss-Legendre (GL)
quadrature points obtained by the tensor product of the GL quadrature points $\xi_{\text{GL}}^m, \eta_{\text{GL}}^m, \zeta_{\text{GL}}^m$ in the unit interval $[0,1]$, see \cite{stroud}.

The discontinuous finite element data representation in~\eqref{eqn.uh} leads naturally to 
i) a Discontinuous Galerkin (DG) scheme if $N>0$ and $N=M$, where the desired order of accuracy $M$ already coincides 
	with the degree $N$ of the polynomial approximating the data ($M=N$), so that high order of accuracy in space can be obtained without the use of any spatial reconstruction operator,
and to ii) a Finite Volume (FV) scheme in the case $N=0$. 
	This indeed means that for $N=0$ we have $\phi_\ell(\BoldXi) = 1$ with $\ell=0$, 
	and~\eqref{eqn.uh} reduces to the classical piecewise constant data representation that is typical of finite volume schemes, 
	where the only degree of freedom per element is the usual cell average $\hat{\mathbf{u}}_{0}$. 
	In this case the order of accuracy $M$ in space will be obtained through the reconstruction procedure described in next Section~\ref{ssec.reconstruction}.
However, iii) also a family of hybrid reconstructed Discontinuous Galerkin methods is included in this representation, where 
	a Hermite-type reconstruction of degree $M > N$ is performed on cell data represented by polynomials of degree $N$, see the next Section~\ref{ssec.reconstruction}.

Thus, within the general $P_NP_M$ formalism one can simultaneously deal with arbitrary high order FV and DG schemes 
and reconstructed hybrid methods inside a unified framework, 
with only very few differences between the different schemes (substantially the reconstruction procedure and the type of limiter).

\subsection{High order spatial reconstruction (order $M$)}
\label{ssec.reconstruction}

\begin{figure}[h]
	\centering
	\subfigure[Data $\PN{1}{2}$]{\includegraphics[trim= 50 10 50 15,clip,width=0.33\linewidth]{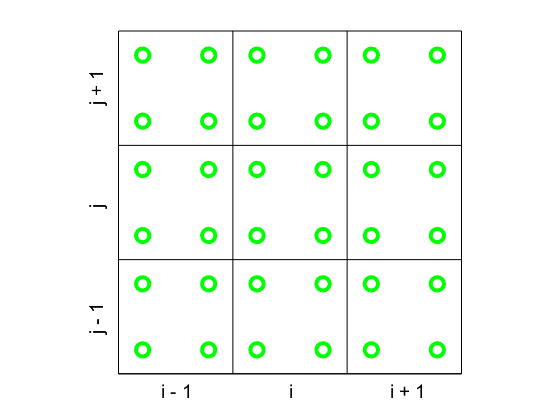}}%
	\subfigure[Data $x$-rec, part 1/2]{\includegraphics[trim= 50 10 50 15,clip,width=0.33\linewidth]{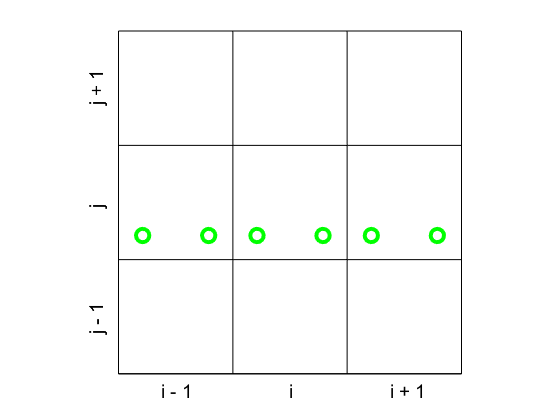}}%
	\subfigure[$x$-rec, part 1/2]{\includegraphics[trim= 50 10 50 10,clip,width=0.33\linewidth]{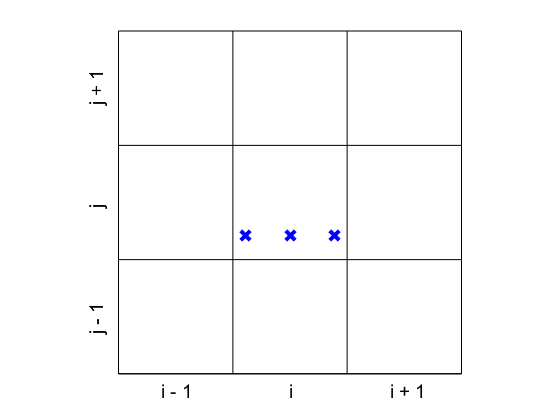}}\\
	\subfigure[Data $x$-rec, part 2/2]{\includegraphics[trim= 50 10 50 10,clip,width=0.33\linewidth]{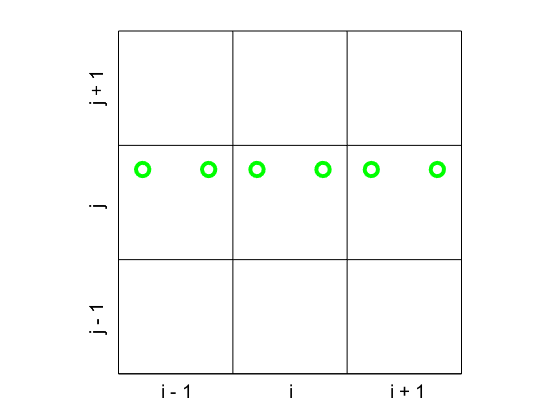}}%
	\subfigure[$x$-rec, part 2/2]{\includegraphics[trim= 50 10 50 10,clip,width=0.33\linewidth]{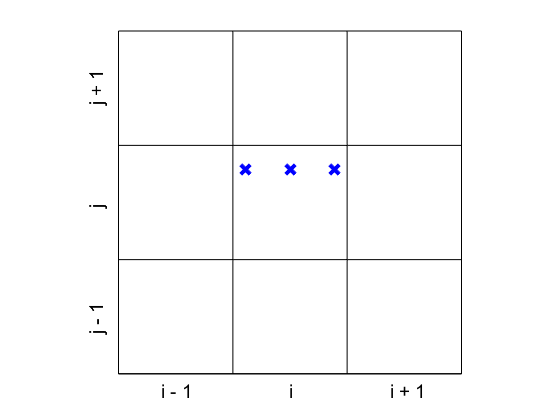}}%
	\subfigure[$x$-rec cell $\Omega_{ij}$]{\includegraphics[trim= 50 10 50 10,clip,width=0.33\linewidth]{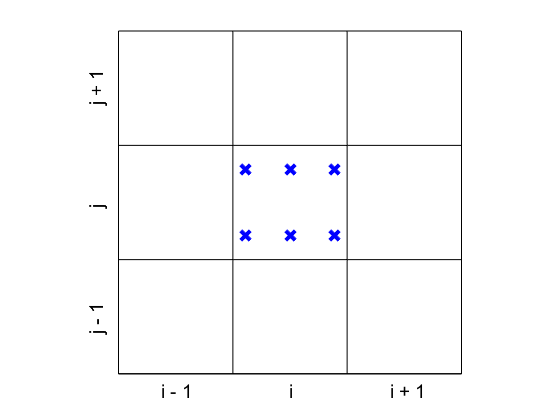}}\\
	\caption{Reconstruction $\PN{1}{2}$ in cell $\Omega_{ij}$ along the $x$-direction in $d=2$ dimensions.
		Since we are employing nodal basis functions, we can represent the available information at each stage of our $\PN{N}{M}$ scheme
		in each cell by a symbol located at a certain GL point inside the cell. 
		In a $\PN{1}{2}$ scheme $\u_h$ is represented by a $P_N=P_1$ polynomial, 
		so we have $(N+1)^2=4$ information in each cell (a, the green circles).
		By selecting the $N+1=2$ information along the same horizontal section in $\Omega_{ij}$ and 
		in its two immediate neighbors $\Omega_{i-1,j}, \Omega_{i+1,j}$ (b), 
		we have enough information ($3(N+1)=6>3=M+1$) in order to reconstruct a $P_M=P_2$ polynomial
		in $x$-direction (c); 
		then we have to repeat the same procedure for each $N+1=2$ horizontal section of cell $\Omega_{ij}$ (d-e). 
		In this way we obtain our reconstructing polynomial in the $x$-direction, represented by $(M+1)(N+1)=6$ information (f, the blue crosses).}
	\label{fig.recP1P2_x}
\end{figure}

\begin{figure}[p]
	\centering
	\subfigure[Data after $x$-rec]{\includegraphics[trim= 50 10 50 15,clip,width=0.33\linewidth]{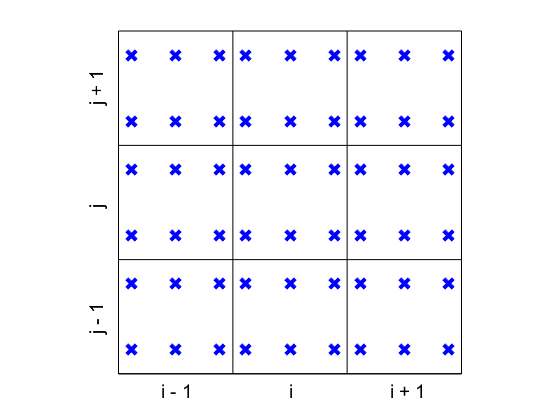}}
	\subfigure[Data after $x$-rec in $\mathcal{S}^y$]{\includegraphics[trim= 50 10 50 15,clip,width=0.33\linewidth]{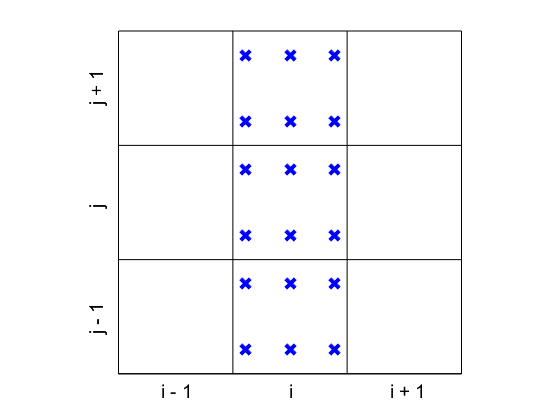}}%
	\subfigure[Data $y$-rec, part 1/3]{\includegraphics[trim= 50 10 50 15,clip,width=0.33\linewidth]{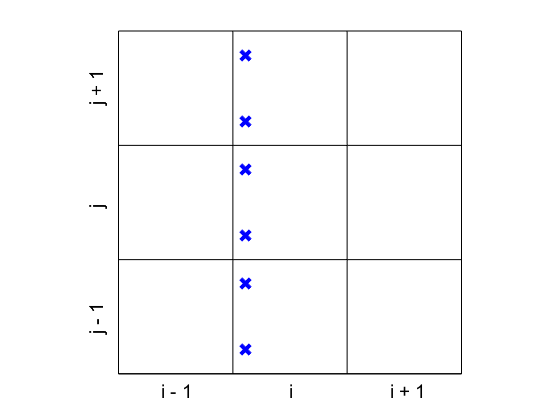}}\\
	\subfigure[$y$-rec, part 1/3]{\includegraphics[trim= 50 10 50 15,clip,width=0.33\linewidth]{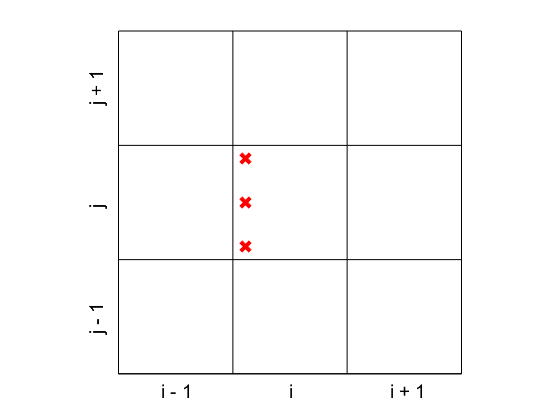}}%
	\subfigure[Data $y$-rec, part 2/3]{\includegraphics[trim= 50 10 50 15,clip,width=0.33\linewidth]{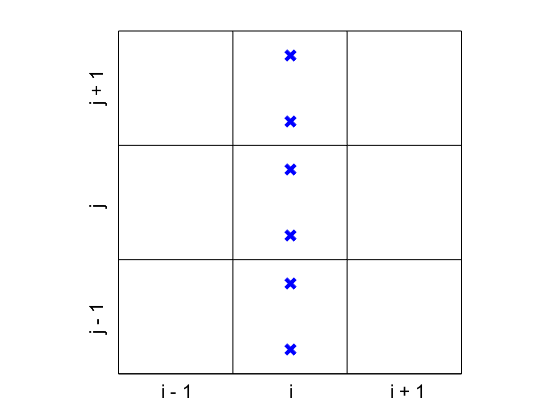}}%
	\subfigure[$y$-rec, part 1/3]{\includegraphics[trim= 50 10 50 15,clip,width=0.33\linewidth]{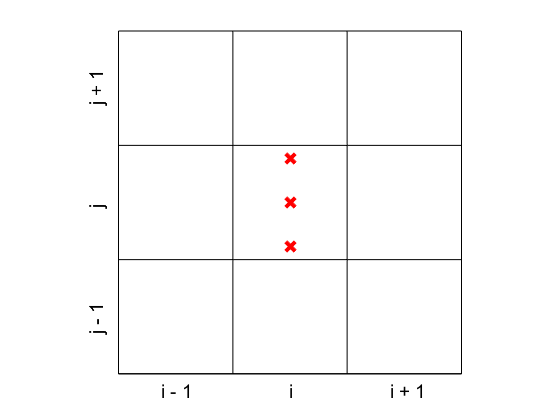}}\\
	\subfigure[Data $y$-rec, part 3/3]{\includegraphics[trim= 50 10 50 15,clip,width=0.33\linewidth]{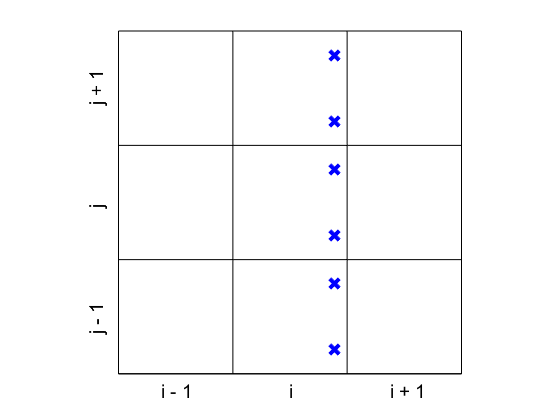}}%
	\subfigure[$y$-rec, part 1/3]{\includegraphics[trim= 50 10 50 15,clip,width=0.33\linewidth]{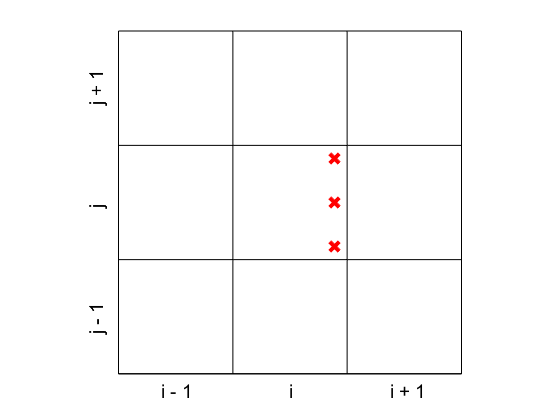}}%
	\subfigure[Reconstruction in $\Omega_{ij}$]{\includegraphics[trim= 50 10 50 15,clip,width=0.33\linewidth]{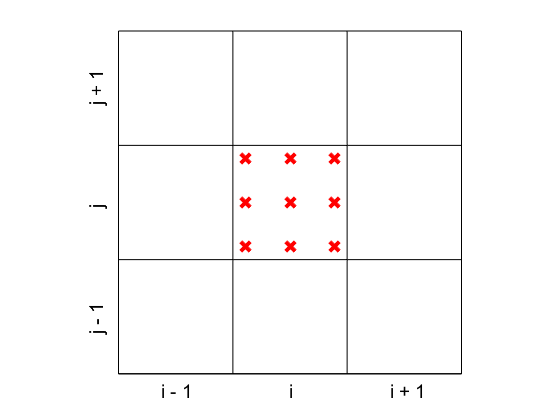}}%
	\caption{Reconstruction $\PN{1}{2}$ in cell $\Omega_{ij}$ along the $y$-direction in $d=2$ dimensions.
		After having performed the reconstruction along the $x$ direction, in each cell we have $(M+1)(N+1) = 6$ information (a, the blue crosses).
		Now, by selecting $(N+1)=2$ information along the same vertical section in $\Omega_{ij}$ and 
		in its two immediate neighbors $\Omega_{i,j-1}, \Omega_{i,j+1}$ (c), 
		we have enough information ($3(N+1)=6>3=M+1$) in order to reconstruct a $P_M=P_2$ polynomial
		in $y$-direction (d); then we have to repeat the same procedure for each $M+1=3$ vertical section of cell $\Omega_{ij}$ (e-h). 
		In this way we obtain our final $P_M=P_2$ reconstructing polynomial for the cell $\Omega_{ij}$ (i, the red crosses).}
	\label{fig.recP1P2_y}
\end{figure}

\begin{figure}[p]
	\centering
	\subfigure[Data $\PN{2}{4}$]{\includegraphics[trim= 50 10 50 15,clip,width=0.33\linewidth]{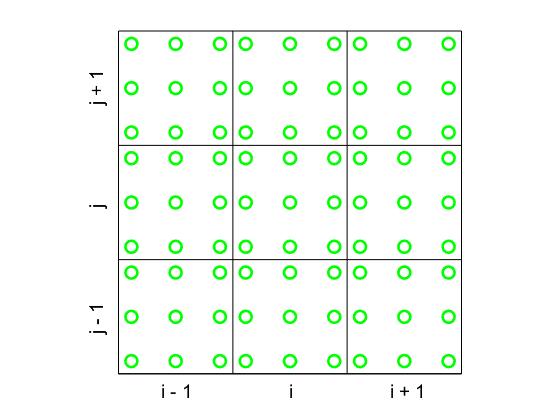}}%
	\subfigure[Data $x$-rec, part 1/3]{\includegraphics[trim= 50 10 50 15,clip,width=0.33\linewidth]{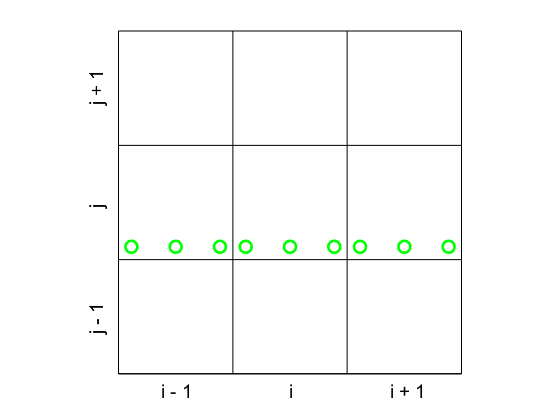}}%
	\subfigure[$x$-rec, part 1/3]{\includegraphics[trim= 50 10 50 15,clip,width=0.33\linewidth]{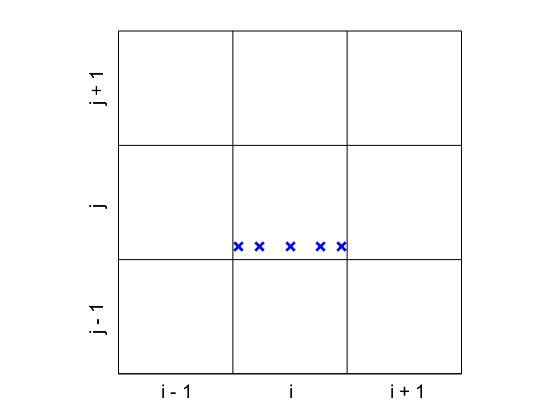}}\\
	\subfigure[$x$-rec cell $\Omega_{ij}$]{\includegraphics[trim= 50 10 50 15,clip,width=0.33\linewidth]{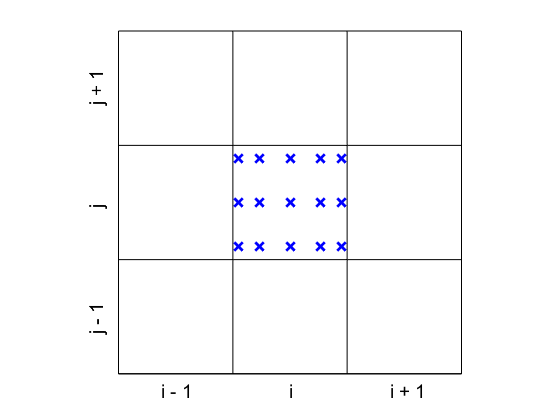}}%
	\subfigure[Data $y$-rec, part 2/5]{\includegraphics[trim= 50 10 50 15,clip,width=0.33\linewidth]{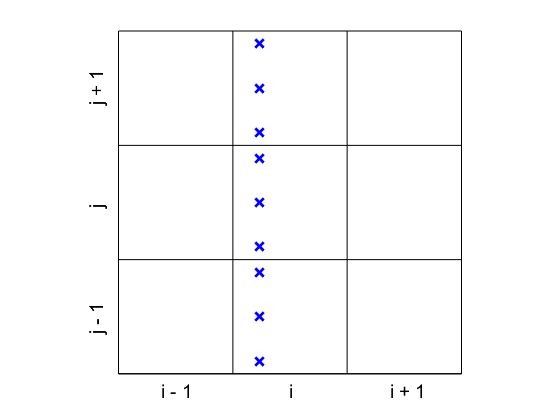}}%
	\subfigure[$y$-rec, part 2/5]{\includegraphics[trim= 50 10 50 15,clip,width=0.33\linewidth]{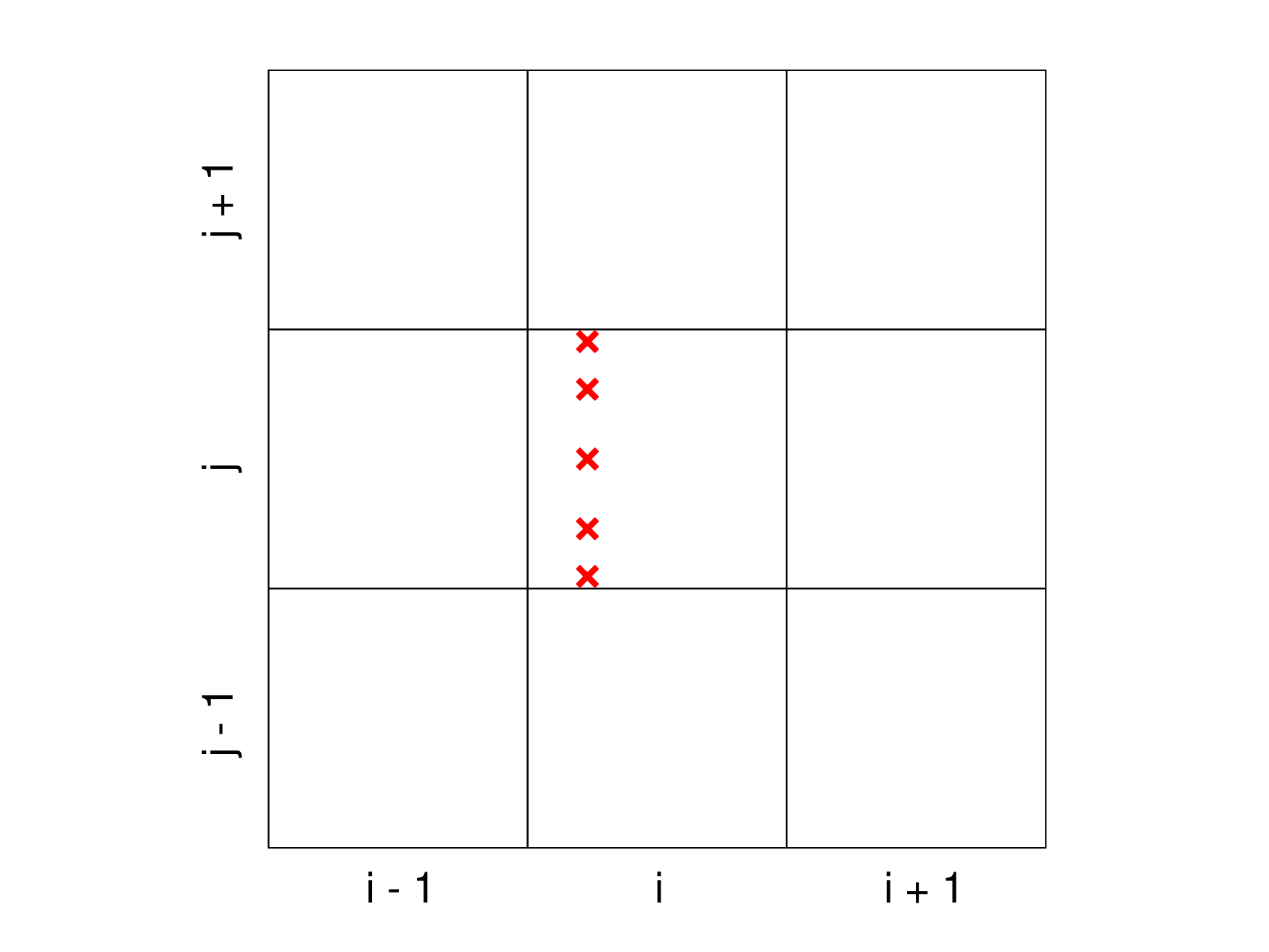}}\\
	\subfigure[Data $y$-rec, part 4/5]{\includegraphics[trim= 50 10 50 15,clip,width=0.33\linewidth]{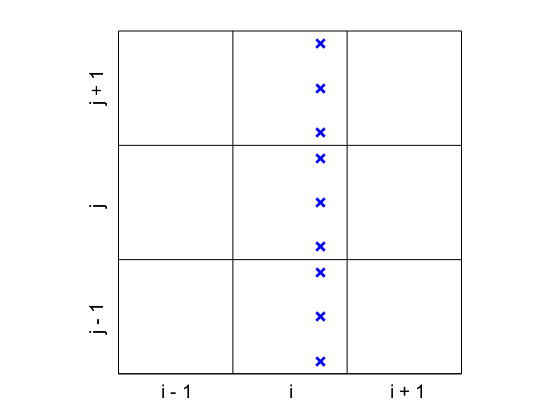}}
	\subfigure[$y$-rec, part 4/5]{\includegraphics[trim= 50 10 50 15,clip,width=0.33\linewidth]{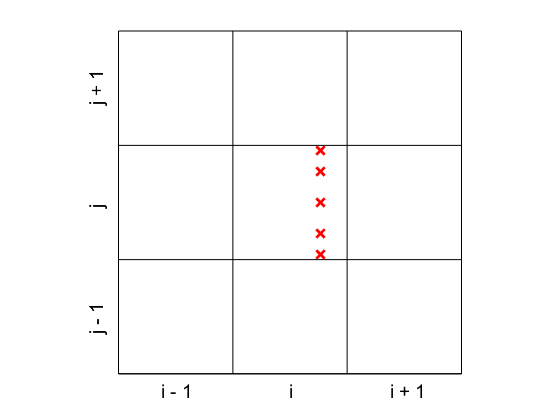}}%
	\subfigure[Reconstruction in $\Omega_{ij}$]{\includegraphics[trim= 50 10 50 15,clip,width=0.33\linewidth]{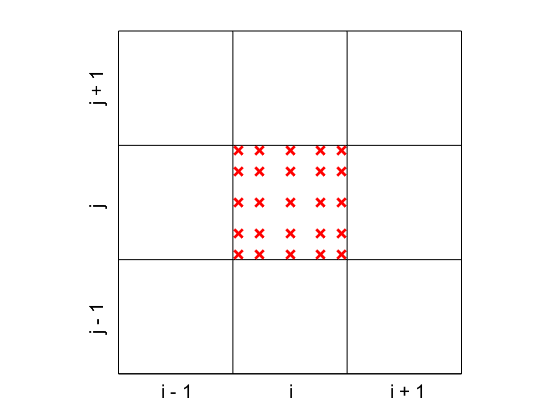}}%
	\caption{Reconstruction $\PN{2}{4}$ in cell $\Omega_{ij}$ in $d=2$ dimensions. 
		The available data (a, green circles) are provided by the $P_N=P_2$ polynomial $\u_h$. 
		By selecting $(N+1)=3$ information along the same horizontal section in $\Omega_{ij}$ and 
		in its two immediate neighbors $\Omega_{i-1,j}, \Omega_{i+1,j}$ (b), 
		we have enough information ($3(N+1)=9>5=M+1$) in order to reconstruct a $P_M=P_4$ polynomial
		in $x$-direction (c); 
		then the same procedure has to be repeated for each $N+1=3$ horizontal section of cell $\Omega_{ij}$ obtaining (d), and finally for each cell of the domain.		 
		At this point in each cell we have $(M+1)(N+1) = 15$ information (d, the blue crosses) and 		
	 	by selecting $(N+1)=3$ information along the same vertical section in $\Omega_{ij}$ and 
		in its two immediate neighbors $\Omega_{i,j-1}, \Omega_{i,j+1}$ (e), 
		we have enough information to reconstruct a $P_M=P_4$ polynomial
		in $y$-direction (f); then we have to repeat the same procedure for each $M+1=5$ vertical section of cell $\Omega_{ij}$ (g-h). 
		In this way we obtain our final $P_M=P_5$ reconstructing polynomial for the cell $\Omega_{ij}$ (i, the red crosses).}
	\label{fig.recP2P4}
\end{figure}

In the framework of $P_NP_M$ schemes, $M$ indicates the highest polynomial approximation degree used for the representation of the discrete solution within the method.  
Hence, in this Section we describe the reconstruction procedure that is needed to obtain approximation degree $M$ in \textit{space} from an underlying data representation $\mathbf{u}_h(\mathbf{x},t^n)$ of lower or equal degree $N \leq M$, i.e. 
the procedure that generates a spatially high order accurate reconstruction polynomial  $\mathbf{w}_h(\x,t^n)$ of degree $M$  
\begin{equation}
\mathbf{w}_h(\x,t^n) \!=\!\! \sum \limits_{\ell=0}^{\mathcal{M}-1} \!\psi_\ell(\x,t^n) \, \hat{\mathbf{w}}_{\ell,\I} 
:= \psi_\ell(\x,t^n) \, \hat{\mathbf{w}}_{\ell,\I}, \ \ \x \in \Omegai, \ \ \mathcal{M}\!=\!(M+1)^d\!,
\label{eqn.wh}
\end{equation}
where we formally employ the same nodal basis functions for the reconstruction and for the data representation, see~\eqref{eqn.uh}.
However, note that when $M \ne N$ of course $\psi_l(\x,t^n)$ does not coincide with $\phi_l(\x,t^n)$, since the polynomial degree and the positions of the GL points are not the same.

For the sake of a uniform notation, when $M=N$, we trivially impose that 
the reconstruction polynomial is given by the DG polynomial, i.e.  
$\mathbf{w}_h(\x,t^n)=\mathbf{u}_h(\x,t^n)$, which automatically implies
that in the case $N=M$ the reconstruction operator is simply the identity.

In the other cases, we employ a polynomial reconstruction procedure implemented in a \textit{dimension by dimension} fashion 
in order to compute the coefficients $\hat{\mathbf{w}}_{\ell,\I}$ in~\eqref{eqn.wh}.
To better follow the following reasoning we refer the reader also to the Figures~\ref{fig.recP1P2_x},~\ref{fig.recP1P2_y}, and~\ref{fig.recP2P4}.
Focusing on the reconstruction procedure along the $x$-direction, given an element $\Omegai=\Omega_{ijk}$, 
we write the reconstruction polynomial in $x$-direction $\w_h^x$ in terms of one dimensional basis functions as 
\begin{eqnarray} 
\mathbf{w}^{x}_{h}(\mathbf{x},t^n)  &=& 
\sum_{\ell_{1}=0}^{M}  \sum_{r_2=0}^{N}  \sum_{r_3=0}^{N} \psi_{\ell_{1}}\left(\xi \right) \phi_{r_2}\left(\eta \right) \phi_{r_3}\left(\zeta \right) \hat{\mathbf{w}}^x_{\ell_{1},r_2,r_3,\mathbf{i}} \nonumber \\ 
&:=&
					   \psi_{\ell_{1}}\left(\xi \right) \phi_{r_2}\left(\eta \right) \phi_{r_3}\left(\zeta \right) \hat{\mathbf{w}}^x_{\ell_{1},r_2,r_3,\mathbf{i}}.
\end{eqnarray}
Then, we integrate on a set $\mathcal{S}^x$ of neighbors of $\Omegai$ in $x$-direction, obtaining an algebraic system for the polynomial coefficients 
$ \hat{\mathbf{w}}_{\ell_{1},r_2,r_3,\mathbf{i}}$ (one for each horizontal section of $\Omegai$) 
\be
& \frac{1}{\Delta x_{m}}  \int_{x_{m-\frac{1}{2}}}^{x_{m+\frac{1}{2}}} \psi_{\ell_{1}} \phi_{r_2} \phi_{r_3} \hat{\mathbf{w}}^x_{\ell_{1},r_2,r_3,\mathbf{i}} dx = 
  \frac{1}{\Delta x_{m}}  \int_{x_{m-\frac{1}{2}}}^{x_{m+\frac{1}{2}}} \phi_{r_1} \phi_{r_2} \phi_{r_3} \hat{\mathbf{u}}_{r_1,r_2,r_3,\mathbf{m}}\,dx,  \\[3pt]
& \forall \Omega_{\mathbf{m}} \in \mathcal{S}^x, \quad \forall \ell_1 \in [1,M+1], \quad \forall r_1, r_{2}, r_{3} \in [1, N+1],
\label{eq.recX}
\ee
where $\mathcal{S}^x$ contains the neighbors along the $x$ axis, i.e $\Omega_{\mathbf{m}}=\Omega_{mjk}$ such that 
$m \in \left\{ i-1, i, i+1\right\}$. Note that we are using a stencil made always of only $3$ elements in each direction, thus a \textit{very compact} one. 
Indeed, a stencil composed of $3$ elements is enough for any $\PN{N}{M}$ scheme with $M\le3N+2$, 
because any $P_N$ cell contains $N+1$ degrees of freedom, thus $3$ cells provide $3N+3$ degrees of freedom, 
which are sufficient for a polynomial reconstruction of degree up to $3N+2$.
Moreover, when the provided information are more than the minimum required, 
the system~\eqref{eq.recX} results to be overdetermined; so, to solve it, we employ 
a constrained least-squares technique (CLSQ)~\cite{DumbserKaeser06b}, i.e. we impose that the reconstructed polynomial satisfies 
\be
& \frac{1}{\Delta x_{i}}  \int_{x_{i-\frac{1}{2}}}^{x_{i+\frac{1}{2}}} \psi_{\ell_{1}} \phi_{r_2} \phi_{r_3} \hat{\mathbf{w}}^x_{\ell_{1},r_2,r_3,\mathbf{i}} dx = 
\frac{1}{\Delta x_{i}}  \int_{x_{i-\frac{1}{2}}}^{x_{i+\frac{1}{2}}} \phi_{r_1} \phi_{r_2} \phi_{r_3} \hat{\mathbf{u}}_{r_1,r_2,r_3,\mathbf{i}}\,dx,  \\[3pt]
&  \quad \forall \ell_1 \in [1,N+1], \quad \forall r_1, r_{2}, r_{3} \in [1, N+1],
\label{eq.recX}
\ee
\textit{exactly}. In other words, all moments of the reconstructed solution $\mathbf{w}_h$ and the original solution $\mathbf{u}_h$ up to degree $N$ must \textit{coincide exactly} within cell $\Omegai$  
and match on the remaining stencil elements in the least-square sense.

To complete the reconstruction polynomial, we now repeat the above procedure in the $y$-direction, 
so we write the reconstruction polynomial in terms of one-dimensional basis functions as 
\be
\mathbf{w}^{y}_{h}(x,y,t^n)  =  \psi_{\ell_{1}}\left(\xi \right) \psi_{\ell_{2}}\left(\eta \right) \phi_{r_3}\left(\zeta \right)  \hat{\mathbf{w}}^y_{\ell_{1},\ell_{2},r_3,\mathbf{i}},
\ee
and we solve the algebraic system
\be
& \frac{1}{\Delta y_{n}}  \int_{y_{n-\frac{1}{2}}}^{y_{n+\frac{1}{2}}} \! \psi_{\ell_{1}}  \psi_{\ell_{2}}  \phi_{r_3}  \hat{\mathbf{w}}^y_{\ell_{1},\ell_{2},r_3,\mathbf{i}}\, dy \!=\! 
  \frac{1}{\Delta y_{n}}  \int_{y_{n-\frac{1}{2}}}^{y_{n+\frac{1}{2}}} \! \psi_{\ell_{1}}  \phi_{r_2}  \phi_{r_3}  \hat{\mathbf{w}}^x_{\ell_{1},r_2,r_3,\mathbf{n}}\, dy, \\[3pt]
& \forall \Omega_{\textbf{n}} \!\in\! \mathcal{S}^{y}\!, \quad \forall \ell_1, \ell_2 \in [1,M+1] , \quad \forall r_2, r_3 \in [1,N+1],
\ee
with $\mathcal{S}^y$ being the set of the neighbors along the $y$ axis, i.e $\Omega_{\mathbf{n}}=\Omega_{ink}$ such that $n \in \left\{ j-1, j, j+1 \right\}$. Again, for overdetermined systems we impose that the reconstruction exactly satisfies 
\be
& \frac{1}{\Delta y_{j}}  \int_{y_{j-\frac{1}{2}}}^{y_{j+\frac{1}{2}}} \! \psi_{\ell_{1}}  \psi_{\ell_{2}}  \phi_{r_3}  \hat{\mathbf{w}}^y_{\ell_{1},\ell_{2},r_3,\mathbf{i}}\, dy \!=\! 
\frac{1}{\Delta y_{j}}  \int_{y_{j-\frac{1}{2}}}^{y_{j+\frac{1}{2}}} \! \psi_{\ell_{1}}  \phi_{r_2}  \phi_{r_3}  \hat{\mathbf{w}}^x_{\ell_{1},r_2,r_3,\mathbf{n}}\, dy, \\[3pt]
& \quad \forall \ell_1 \in [1,M+1], \quad \forall \ell_2 \in [1,N+1] , \quad \forall r_2, r_3 \in [1,N+1].
\ee
And then the same procedure can be repeated along the $z$ axis by looking for the unknown coefficients
$\hat{\mathbf{w}}_{\ell_{1},\ell_{2},\ell_{3},\mathbf{i}}$ of
\be
\mathbf{w}_{h}(x,y,z,t^n)  =  \psi_{\ell_{1}}\left(\xi \right) \psi_{\ell_{2}}\left(\eta \right) \psi_{\ell_{3}}\left(\zeta \right)  \hat{\mathbf{w}}_{\ell_1,\ell_{2},\ell_3,\mathbf{i}},
\ee
and solving the algebraic system
\be
& \frac{1}{\Delta z_{p}}  \int_{z_{p-\frac{1}{2}}}^{z_{p+\frac{1}{2}}} \! \psi_{\ell_{1}} \psi_{\ell_{2}} \psi_{\ell_{3}}    \hat{\mathbf{w}}_{\ell_{1},\ell_{2},\ell_{3},\mathbf{i}}\, dy  = 
  \frac{1}{\Delta z_{p}}  \int_{z_{p-\frac{1}{2}}}^{z_{p+\frac{1}{2}}} \!\! \psi_{\ell_{1}} \psi_{\ell_{2}} \phi_{r_{3}}  \hat{\mathbf{w}}^y_{\ell_{1}\ell_{2},r_3,\mathbf{p}}\, dy,  \\[3pt]
& \forall \Omega_{\textbf{p}} \!\!\in\! \mathcal{S}^{z}\!, \quad \forall \ell_{1}, \ell_{2}, \ell_{3} \in [1, M+1], \quad \forall r_3 \in [1,N+1], 
\ee
 $\mathcal{S}^z$ being the set of the neighbors along the $z$ axis, i.e $\Omega_{\mathbf{p}}=\Omega_{ijp}$ such that 
$p \in \left[k-1, k, k+1\right]$. For overdetermined systems, the constraint reads 
\be
& \frac{1}{\Delta z_{k}}  \int_{z_{k-\frac{1}{2}}}^{z_{k+\frac{1}{2}}} \! \psi_{\ell_{1}} \psi_{\ell_{2}} \psi_{\ell_{3}}    \hat{\mathbf{w}}_{\ell_{1},\ell_{2},\ell_{3},\mathbf{i}}\, dy  = 
\frac{1}{\Delta z_{k}}  \int_{z_{k-\frac{1}{2}}}^{z_{k+\frac{1}{2}}} \!\! \psi_{\ell_{1}} \psi_{\ell_{2}} \phi_{r_{3}}  \hat{\mathbf{w}}^y_{\ell_{1}\ell_{2},r_3,\mathbf{i}}\, dy,  \\[3pt]
& \quad \forall \ell_{1}, \ell_{2}, \in [1, M+1], \quad \forall  \ell_{3} \in [1, N+1], \quad \forall r_3 \in [1,N+1]. 
\ee
Finally, the coefficients $\hat{\mathbf{w}}_{\ell_{1},\ell_{2},\ell_{3},\mathbf{i}}$
represent the $\hat{\mathbf{w}}_{\ell,\I}$ of~\eqref{eqn.wh} 
that give us the desired polynomial representation of order $M$ in space.

We would like to emphasize that the reconstructed $\PN{N}{M}$ schemes with $N>0$ are very compact because for the reconstruction 
they need a much smaller stencils than classical finite volume schemes and that for regular Cartesian meshes, 
the coefficients of the above constraint least squares systems depend only on the choice of the basis
functions, hence the integrals can be precomputed once and for all on the reference element before starting the simulation. 

Furthermore, in the specific case $N=0$, i.e. when the $\PN{N}{M}$ reduces to a FV scheme, 
the above polynomial reconstruction procedure must be made nonlinear; 
this can be easily done, for example, by adopting the WENO strategy specifically described in the context of $P_NP_M$ type schemes 
(thus with the same notation adopted here) on Cartesian meshes in~\cite{AMR3DCL,FrontierADERGPR}. 
We recall that the nonlinearity introduced through ENO/WENO type procedures \textit{essentially} avoids 
the spurious oscillations typical of high order linear schemes modeling discontinuous processes   see~\cite{godunov}, 
and was already introduced in the 80s and subsequently largely developed~\cite{HartenENO,eno,shu1,JiangShu1996,balsarashu,ZhangShu3D,shu2016high}.
Due to the already exhaustive  literature available on FV schemes, 
here, for what concerns the strategies that guarantee robustness on discontinuities, 
we focus on $\PN{N}{M}$ schemes only with $N>0$: 
indeed, it is for those schemes that we propose in this work a \textit{new} strategy, 
i.e. the new \textit{a posteriori} subcell FV limiter described in Section~\ref{ssec.limiter}.

\subsection{High order in time via a local space-time Galerkin predictor}
\label{ssec.predictor}

We recall that high order of accuracy in \textit{space} is provided by the piecewise polynomial data representation $\w_h$ of~\eqref{eqn.wh},
obtained in the previous Section~\ref{ssec.reconstruction}.

Now, in order to achieve also high order of accuracy in \textit{time}, 
relying on the ADER predictor-corrector approach, 
we need to compute the so-called space-time Galerkin \textit{predictor},
i.e. a space-time polynomial  $\q_h$ of degree $M$ in $(d+1)$-dimensions ($d$ for the space plus $1$ for the time) which takes the following form 
\be
& \q_h(\x, t) = \q_h(\BoldXi(\x), \tau(t)) = \sum_{\ell=0}^{\mathcal{Q}-1} \theta_\ell (\BoldXi, \tau) \hat{\q}_\ell
= \theta_\ell (\BoldXi, \tau) \hat{\q}_\ell, \\
& \x \in \Omegai,\quad t \in [t^n, t^{n+1}],\quad \mathcal{Q} = (M+1)^{d+1},
\label{eqn.qh}
\ee 
where again $\theta_\ell (\BoldXi, \tau)$ is given by the tensor product of Lagrange interpolation polynomials 
$\phi_{\ell}\left(\BoldXi (\x) \right) \phi_{\ell_\tau}\left(\tau\right)$, with 
$\BoldXi(\x)$ given by~\eqref{eq.mapping} and the mapping for the time coordinate given by $t = t^n + \tau \Delta t, \tau \in [0,1]$.
This high order polynomial in space and time will serve as a predictor solution, only valid inside $\Omegai\times[t^n,t^{n+1}]$,
to be used for evaluating the numerical fluxes and the sources when integrating the PDE 
in the final corrector step of the ADER scheme, see Section~\ref{ssec.PNPM}.

In order to determine the unknown coefficients $\hat{\q}_\ell$ of~\eqref{eqn.qh} we search $\q_h$ such that it satisfies 
a weak form of the governing PDE~\eqref{eq.generalform} integrated in space and time locally \textit{inside} each $\Omegai$ 
{(with $\Omegai^{\circ} = \Omegai \backslash \partial \Omegai$ being the interior of $\Omegai$)}
\be
& \int_{t^n}^{t^{n+1}} \!\! \int_{\Omegai^{\circ} } \theta_k \, \partial_t \q_h \,d\mbf{x} \, dt 
+ \int_{t^n}^{t^{n+1}} \!\! \int_{\Omegai^{\circ} } \theta_k \,\nabla \cdot \F(\q_h) \,d\x\,dt 
= \mbf{0},
\label{eq.predictorEq}
\ee
where the first term is integrated in time by parts exploiting the \emph{causality principle}
(upwinding in time) 
\be
& \int_{\Omegai^{\circ}}  \theta_k(\x,t^{n+1}) \q_h(\x,t^{n+1}) \, d\x -
\int_{\Omegai^{\circ}} \theta_k(\x,t^{n}) \mbf{w}_h(\x,t^{n}) \, d\x - \\
&\int_{t^n}^{t^{n+1}} \!\!\!  \int_{\Omegai^{\circ} } \!\!\! \partial_t \theta_k(\x,t) \q_h(\x,t) \,d\x\, dt + 
 \int_{t^n}^{t^{n+1}} \!\!\!  \int_{\Omegai^{\circ} } \!\!\! \theta_k(\x,t) \nabla \cdot \F(\q_h(\x,t)) \,d\x\,dt 
= \mbf{0},
\label{eq:DOFpredictor}
\ee
and $\mbf{w}_h(\x,t^{n})$ is the known initial condition at time $t^n$. 

Now, the system~\eqref{eq:DOFpredictor}, which contains only volume integrals to be calculated inside $\Omegai$ and no surface integrals,  
can be solved via a simple discrete Picard iteration for each element $\Omegai$, 
and there is no need of any communication with neighbor elements. 
Indeed, the so-called \textit{predictor} step consists in a {\it local} solution of the governing PDE~\eqref{eq.generalform} \textit{in the small}, 
see~\cite{eno}, inside each space-time element $\Omegai\times[t^n,t^{n+1}]$.
It is called \textit{local} because it is obtained by only considering cell $\Omegai$ with initial data $\w_h$, 
the governing equations~\eqref{eq.generalform} and the geometry, without taking into account any interaction between $\Omegai$ and its neighbors.  
We also want to emphasize that this procedure is exactly the same whatever $N$ and $M$ are.

We recall that this procedure has been introduced for the first time in~\cite{Dumbser2008} for unstructured meshes, 
it has been extended for example to moving meshes in~\cite{Lagrange2D} and
to degenerate space time elements in~\cite{GaburroAREPO}; 
finally, its convergence has been formally proved in~\cite{FrontierADERGPR}.

\subsection{High order fully-discrete one-step ADER $\PN{N}{M}$ scheme}
\label{ssec.PNPM}


Last, the update formula of our ADER $\PN{N}{M}$ scheme is recovered starting from 
the weak formulation of the governing equations~\eqref{eq.generalform}
(where the test functions $\phi_k$ coincide with the basis functions $\phi_\ell$ of~\eqref{eq.nodalGeneralBasis})
\be
\int_{t^n}^{t^{n+1}} \int_{\Omegai} \phi_k \left(\partial_t \Q  + \nabla \cdot \F(\Q) 
\right) 
= \mbf{0};
\label{eqn.pde.weak}
\ee
we then substitute $\Q$ with~\eqref{eqn.uh} at time $t=t^n$ (the known initial condition) 
and at $t=t^{n+1}$ (to represent the unknown evolved conserved variables),
and with the high order predictor $\q_h$ previously computed for $t\in[t^n, t^{n+1}]$, obtaining
\be
& \left( \, \int_{\Omegai} \phi_k \phi_l \, d\x\right)
\left( \hat{\u}_\ell^{n+1} - \hat{\u}_\ell^{n} \, \right) 
+ \int_{t^n}^{t^{n+1}} \!\! \int_{\partial \Omegai} \!\! \phi_k \mathcal{F}\left(\q_h^-, \q_h^+ \right) \cdot \mbf{n} \, dS \, dt \ - \\
& \int_{t^n}^{t^{n+1}} \!\! \int_{\Omegai } \!\!\! \nabla \phi_k \cdot \F(\q_h) \,d\x\,dt 
\, =\,  \mbf{0}.
\label{eq:ADER-DG}
\ee
The use of $\q_h$ allows to compute the integrals appearing in~\eqref{eq:ADER-DG} with high order of accuracy in both space and time.

The boundary fluxes  $\mathcal{F}\cdot\mbf{n}$ are obtained by a Riemann solver, 
thus providing the coupling between neighbors, which was neglected in the predictor step. 
In particular, in this work we will employ three types of standard fluxes, namely the Rusanov flux and  the HLL flux,
whose description can be found in~\cite{ToroBook}, and the HLLEM flux for which we refer to~\cite{HLLEM,Dumbser2015}.
For the sake of completness, we report here the expression of the Rusanov flux that reads as follows
\be 
\label{eqn.fluxPC}
\mathcal{F}(\q_h^{-},\q_h^{+}) \cdot \mathbf{n} 
=& \ \frac{1}{2} \left( {\mbf{F}}(\q_h^{+}) + {\mbf{F}}(\q_h^{-}) \right) \cdot {\mbf{n}} 
- \frac{1}{2} s_{\max} \left( \q_h^{+} - \q_h^{-} \right), 
\ee
where $s_{\max}$ is the maximum eigenvalue of the system matrices $\mathbf{A}(\q_h^{+})$ and $\mathbf{A}(\q_h^{-})$ 
being $\mathbf{A}(\Q)=\frac{\partial \mathbf{F}}{\partial \Q}$.
We remark also that due to the discontinuous character of $\q_h$ at the interfaces $\partial \Omegai$, 
$\mathcal{F}\cdot\mbf{n}$ is computed through a numerical flux function 
evaluated over the boundary-extrapolated data $\q_h^{-}$ and $\q_h^{+}$ (i.e the predictors $\q_h$ of two neighbors elements evaluated at the common interface).

Finally, we stress again that the update procedure in~\eqref{eq:ADER-DG} is the same whatever $N$ and $M$ are, 
and allows the contemporary evolution of all the $(N+1)^d$ degrees of freedom of $\u_h$.

\subsubsection{CFL stability constraint}
\label{ssec.CFL}

A very important feature of $P_NP_M$ schemes is linked to the CFL stability constraint.
Since this family of scheme is explicit, the time step $\Delta t$ 
has to be computed according to a (global) Courant-Friedrichs-Levy (CFL) stability condition  
given by 
\be 
\Delta t_{\PN{N}{M}} \, < \,  \text{CFL}_{\PN{N}{M}}\frac{h_{\text{min}} }{d} \frac{1}{|\lambda_{\text{max}}|} 
\,  < \,  \frac{\text{CFL}}{\left(2N+1\right)} \frac{h_{\text{min}} }{d} \frac{1}{|\lambda_{\text{max}}|} 
\label{eq.timestep}
\ee
where $h_{\text{min}}$ is the minimum characteristic mesh-size,  $|\lambda_{\max}|$ is the spectral radius of the system matrix $\mathbf{A}$
and the maximum admissible $\text{CFL}_{\PN{N}{M}}$ number is given in Table~\ref{tab.CFLPnPm}.
In the above formula we wanted also to recall the 
classical CFL condition of Runge-Kutta DG schemes (the one written on the right, with $\text{CFL}<1$) which is just a bit less restrictive 
than the one needed for ADER $\PN{N}{M}$ schemes,
but easier to remember and helpful 
in justifying the stability of our subcell limiter, see formula~\eqref{eq.CFLFV_limiter}.

\begin{table*}[h]
	\caption{Maximum admissible CFL number for $\PN{N}{M}$ schemes from second to fifth-order of accuracy} 
	\label{tab.CFLPnPm}
	\begin{center} 	
		\begin{tabular}{|c|ccccccc|} 
			\hline 
			$\text{CFL}_{\PN{N}{M}}$	&	N=0	&	N=1	&	N=2	&	N=3	&	N=4	 &	N=5	 &	N=6  \Tstrut\Bstrut\\
			\hline
			M=0	&  1.0  &   &   &   &   &  &  \Tstrut \\
			M=1	&  1.0  &  0.33 &   &   &  &  &  \\
			M=2	&  1.0  &  0.33 &  0.17 &    &  &  &  \\
			M=3	&  1.0  &  0.33 &  0.17 &   0.1 &  &  &  \\
			M=4	&  1.0  &  0.33 &  0.17 &   0.1 &  0.069 &  &  \\
			M=5	&  1.0  &  0.33 &  0.17 &   0.1 &  0.069 & 0.045 &  \\
			M=6	&  1.0  &  0.33 &  0.17 &   0.1 &  0.069 & 0.045 & 0.038 \Bstrut \\
			\hline 			
		\end{tabular}
	\end{center}
\end{table*}

Furthermore, we would like to emphasize that it is the degree $N$ of the data representation 
that governs the stability of the method and not the polynomial degree $M$ of the reconstruction operator. 
Hence, the reconstructed hybrid methods with $N>0$ and $M>N$ allow for larger
time steps than the pure DG methods ($N=M$) of the same order always maintaining a superior resolution with respect to FV schemes ($N=0$),
fact that further justifies the interest in their development.

\subsection{\textit{A posteriori} subcell finite volume limiter}
\label{ssec.limiter}

Up to now, the presented $P_NP_M$ scheme is high order accurate in space and time and, formally, 
the differences between the FV case ($N=0$) the pure DG case ($N=M$) and the hybrid reconstructed case ($N>0, M>N$) 
are basically only due to the procedure for achieving high order of accuracy in space, 
which is obtained through a WENO reconstruction in the FV case,
a linear reconstruction in the hybrid case 
and is automatic by construction for DG, see Section~\ref{ssec.reconstruction}. 
But this is actually a major difference, 
because the WENO operator provides a non-linear stabilization of the FV scheme, 
while the $\PN{N}{M}$ schemes with $N>0$ presented so far are unlimited and, as such, 
they are affected by the so-called Gibbs phenomenon, 
i.e. oscillations are likely to appear in presence of shock waves or other discontinuities.
These oscillations can be explained by the Godunov theorem~\cite{godunov}, 
because in this case the scheme is linear in the sense of Godunov. 

As a consequence, a \textit{limiting technique} is required. 
Our strategy is described in detail below and it will be applied whenever $N>0$.

First, we need to consider the numerical solution computed so far $\u_h^{n+1}$ only as a \textit{candidate} solution:
we denote it with $\mathbf{u}^{n+1,*}_{h}(\x,t^{n+1})$. 

Then, following~\cite{DGLimiter1,DGLimiter2,Zanotti2015d,ADERDGVisc,FrontierADERGPR,SolidBodies2020}, 
each element $\Omega_\mathbf{i}$ is divided into $N_\omega=(2N+1)^d$ equal non-overlapping subgrid cells $\omegaialpha$ whose volume is denoted by $|\omegaialpha|$; 
for any cell we define the corresponding \textit{subcell average} of the $\PN{N}{M}$ solution at time $t^n$  
\begin{equation}
\label{eq.projection}
\mathbf{v}_{\I,\alpha}^n(\x,t^n) = \frac{1}{|\omega_{\I,\alpha}|} \int_{\omega_{\I,\alpha}} 
\mathbf{u}_{h}^n(\x,t^n) \, d\x 
=\mathcal{P}(\mathbf{u}_h^n), \quad \forall \alpha \in [1,N_\omega],
\end{equation}
and the candidate subcell averages at time $t^{n+1}$
\begin{equation}
\label{eq.projectionn1}
\mathbf{v}_{\I,\alpha}^{n+1,*}(\x,t^{n+1}) = \frac{1}{|\omega_{\I,\alpha}|} \int_{\omega_{\I,\alpha}} 
\mathbf{u}_{h}^{n+1,*}(\x,t^{n+1}) \, d\x 
=\mathcal{P}(\mathbf{u}_h^{n+1,*}), \quad \forall \alpha \in [1,N_\omega],
\end{equation}
where $\mathcal{P}(\mathbf{u}_h)$ is the $L_2$ projection operator into the space of piecewise constant cell averages. 

Now, we have to mark the \textit{troubled cells}, 
i.e. we have to identify those cells where the solution found through the $\PN{N}{M}$ scheme cannot be accepted
because it may lead to spurious oscillations.
Thus, the candidate solution $\mathbf{v}^{n+1,*}_h$ is checked against a set of \textit{detection criteria}.
Here we follow the criteria described in~\cite{ALEDG}, however also other 
specific physical bounds or more elaborate choices as those of invariant domain preserving methods \cite{guermond2018second} could be considered.

First, we require that the computed solution is physically acceptable, 
i.e. that it belongs to the phase space of the conservation law being solved. 
For instance, if the compressible Euler equations for gas dynamics are considered, 
density and pressure should be positive and in practice we require that they are greater than a prescribed tolerance $\epsilon=10^{-12}$. 
Then, the solution should verify a relaxed discrete maximum principle (DMP) 
\be
\min_{\mbf{m} \in \mathcal{V}(\Omegai)} \!\!\left(\min_{\beta\in[1,N_\omega]}\!(\mathbf{v}_{\mbf{m},\beta\,}^n) \!\!\right) 
\!\!- \delta \leq  \mathbf{v}^{n+1,*}_{\I,\alpha} \!\! \leq \!\!\!
\max_{\mbf{m} \in \mathcal{V}(\Omegai)} \!\!\left( \max_{\beta\in[1,N_\omega]} \!(\mathbf{v}_{\mbf{m},\beta\,}^n) \!\!\right) 
\!\!+\! \delta, \
 \forall \alpha \! \in \! [1,N_\omega],
\label{eqn.RDMP}
\ee
where ${\cal{V}}(\Omegai)$ is the set containing all the neighbors of $\Omegai$
sharing a common node with $\Omegai$,
and $\delta$ is a parameter which, according to~\cite{ALEDG,DGLimiter1,DGLimiter2}, reads
\begin{equation}
\delta \!=\! \max \Biggl( \delta_0, \, \epsilon \cdot\! 
\biggl[
\max_{\mbf{m} \in \mathcal{V}(\Omegai)} \!\biggl( \max_{\beta\in[1,N_\omega]} (\mathbf{v}_{\mbf{m},\beta\,}^n) \biggr) \!-\!\!\!
\min_{\mbf{m} \in \mathcal{V}(\Omegai)} \!\biggl( \min_{\beta\in[1,N_\omega]} (\mathbf{v}_{\mbf{m},\beta\,}^n) \biggr) 
\! \biggr]\!
\Biggr),
\label{eqn.deltaRDMP}
\end{equation}
with $\delta_0=10^{-5}$ and $\epsilon=10^{-4}$.
If a cell does not fulfill the detection criteria in \textit{all} its subcells, 
then it is marked as \textit{troubled}. 
It is possible that some false positive activations of the limiter occur; 
however these local effects do not reduce the overall quality of the simulation thanks to the highly accurate limiter procedure adopted on troubled cells.

Then \textit{only} on these troubled cells we apply either a second-order accurate MUSCL-Hancock TVD finite volume scheme with \textit{minmod}
slope limiter~\cite{ToroBook} (in particular in presence of strong shock waves or low density atmospheres),
or a more accurate ADER-WENO FV scheme~\cite{DGLimiter1,AMR3DCL} that better captures local extrema.  
In this way we can re-compute the solution in order to evolve the cell averages $\mathbf{v}_{\I,\alpha}^n$ in time and obtain $\mathbf{v}_{\I,\alpha}^{n+1}$. 

Note that, due to the fact of applying a high order scheme and to do so on a subgrid instead that on the main grid, 
the subcell average representation given by $\mathbf{v}_{\I,\alpha}^{n+1}$ maintains the high resolution of the underling $\PN{N}{M}$ scheme.
Indeed now, we can recover from these cell averages a 
polynomial $\mathbf{u}_{h}^{n+1}$ of degree $N$; this is done 
by applying a reconstruction operator $\mathcal{R}$ such that 
\begin{equation}
\!\int_{\omega_{\I,\alpha}^n} \!\!\!\!\! \mathbf{u}_{h}^{n+1}(\x,t^{n+1}) \, d\x =\!\!\!\int_{\omega_{\I,\alpha}^n} \!\!\!\!\!\mathbf{v}_{\I,\alpha}^{n+1}(\x,t^n) \, d\x :=\mathcal{R}(\mathbf{v}_{\I,\alpha}^{n+1}(\x,t^n)),  \ \forall \alpha \!\in\! [1,N_\omega],
\label{eq.reconstruction}
\end{equation}
which is \textit{conservative} on the main cell $\Omegai$ thanks to the additional linear constraint
\begin{equation}
\int_{\Omegai} \mathbf{u}_{h}^{n+1}(\x,t^{n+1}) \, d\x = \int_{\Omegai} \mathbf{v}_{h}^{n+1}(\x,t^{n+1}) \, d\x.
\label{eqn.LSQ}
\end{equation}
Moreover, the projection operator $\mathcal{P}$ in~\eqref{eq.projection} and the reconstruction operator $\mathcal{R}$ in~\eqref{eq.reconstruction}
satisfy the property $\mathcal{P} \cdot \mathcal{R}=\mathcal{I}$, with $\mathcal{I}$ being the identity operator. 

However, we have to remark that the reconstruction operator~\eqref{eq.reconstruction}-\eqref{eqn.LSQ} might still lead to an oscillatory solution, 
since it is based on a linear unlimited least squares technique. 
If this is the case, the cell $\Omegai$ will be marked again automatically as troubled during the next timestep $t^{n+2}$, 
therefore the same finite volume subcell limiter will be used again in that cell and in particular
the subcell averages from which to start as initial data at time $t^{n+1}$ will be the $\mathbf{v}_{i,\alpha}^{n+1}$
kept in memory from the previous limited step. 

Furthermore, if a cell $\Omegai$ is acceptable but has at least one troubled neighbor in $\mathcal{V}(\Omegai)$, 
then we cannot accept its candidate solution $\mathbf{u}^{n+1,*}_{h}(\x,t^{n+1})$ as it is because the scheme would be nonconservative,
since the numerical flux $\mathcal{F}\cdot\mbf{n}$ at the common interface would have been  
computed in two different ways in $\Omegai$ and its neighbor 
and by using a $\q_h$ that may already be non acceptable.
Thus, the final $\PN{N}{M}$ solution in these cells, neighbors of troubled ones, is rearranged as follows:  
we keep the already computed values of the volume integral and of the surface integrals interacting with non-troubled cells, 
while the numerical flux across the troubled faces is substituted with the one computed through the limiter procedure.

Finally, note that for the subcell FV scheme we have a different CFL stability condition
\be
\label{eq.CFLFV_limiter}
\Delta t_{\text{FV}} < \text{CFL}_\text{FV}\frac{h_{\text{min}}}{d\, N_\omega}\frac{1}{|\lambda_{\text{max}}|},
\ee
with $\text{CFL}_\text{FV}<1$ and $h_{\min}$ the minimum cell size referred to $\Omegai$.
Condition~\eqref{eq.CFLFV_limiter} guides us in choosing the number of employed subcells $N_\omega$.
In particular, our choice  $N_\omega = (2N + 1)^d$ respects the stability condition~\eqref{eq.CFLFV_limiter}
maintaining the original timestep size fixed for the current timestep (see~\eqref{eq.timestep}) 
but also taking into account the maximum possible number of subgrid elements allowed by that timestep size. 

We would like to stress again that our interest for $\PN{N}{M}$ schemes is motivated 
by the high resolution that they are able to provide and the reduced cost offered by the possibility 
of representing data with lower order polynomials of degree $N$ and to achieve, however, an order of accuracy $M$, with $M>N$, using a very compact stencil and a simple reconstruction procedure.

The presented \textit{a posteriori} subcell FV limiter is applied only where it is needed by detecting 
spurious oscillations \textit{a posteriori} and it is based on strong stability preserving FV schemes developed precisely for dealing with discontinuous solutions. Since FV schemes are less accurate than $P_NP_M$ schemes with $N>0$, the limiter is applied on a finer subgrid than the original main grid in order to avoid a loss of useful information.

\subsection{Adaptive Cartesian mesh refinement}
\label{ssec.AMR}

The last ingredient that further increases the resolution of the proposed approach 
is the possibility of activating an Adaptive Mesh Refinement (AMR) technique based on a \textit{cell-by-cell} refinement approach;
indeed, the combined action of our subcell limiter and of AMR allows a sharp detection of all discontinuities.
For details we refer  to~\cite{Berger-Oliger1984,Berger-Colella1989,Khokhlov1998,AMR3DCL,AMR3DNC,DGLimiter2,ADERDGVisc,ADERGRMHD,Peano1,Peano2,Exahype} and we recall here just the main features. 
Our algorithm basically consists in, starting from the main grid~\eqref{eq.MainGrid},
introducing successive refinement levels, in regions of particular interest according to a prescribed refinement criterion.
In particular, we have to fix the following parameters:
\begin{itemize}
	\item
	the maximum level of refinement $\ell_{\rm max}$, typically chosen equal to $2$ or $3$ in our tests;
	\item
	the refinement factor $\rRef$, governing the number of subcells that are generated according to 
	\be
	\label{refine-factor}
	\Delta x_{\ell} = \mathfrak{r} \Delta x_{\ell+1}, \quad \Delta y_{\ell} = \mathfrak{r} \Delta y_{\ell+1}, \quad 
	\Delta z_{\ell} = \mathfrak{r} \Delta z_{\ell+1}, 
	\ee
	where $\Delta x_{\ell}$ is the size of the cell at refinement level number $\ell$ along the $x$-direction, and similarly for the other directions; 
	\item
	the refinement criterion that we base on oscillations of second derivatives, see~\cite{Lohner1987}. 
	In practice, we have to compute 
	\begin{equation}
	\chi_\I=\sqrt{\frac{\sum_{k,l} (\partial^2 \Phi/\partial x_k \partial x_l)^2 }{\sum_{k,l}[(|\partial \Phi/\partial x_k|_{i+1}+|\partial \Phi/\partial x_k|_i)/\Delta x_l+\varepsilon|(\partial^2 /\partial x_k \partial x_l )||\Phi|]^2} },
	\label{eqn.indicator}
	\end{equation}
	where the summation $\sum_{k,l}$ is taken over the number of space dimension of the problem in order to include the cross term derivatives, 
	the parameter $\varepsilon = 0.01$ acts as a filter preventing refinement in regions of small ripples, 
	and the function $\Phi=\Phi(\Q)$, that could be any suitable indicator function of the conserved variables $\Q$, in our test is chosen to be 
	simply $\Phi(\Q)=\rho$.
	Next, a cell $\Omegai$ is marked for refinement if $\chi_\I>\chi_{\rm ref}$, while it is marked for re-coarsening if $\chi_\I<\chi_{\rm rec}$. 	
	In our tests we have chosen $\chi_{\rm ref}$ in the range $[0.2,0.25]$ and $\chi_{\rm rec}$ in $[0.05,0.15]$. 
\end{itemize}

Finally, the numerical solution at the subcell level during a refinement step is obtained by a standard $L_2$ projection, 
while a reconstruction operator is employed to recover the solution on the main grid starting from the subcell level.
Moreover, in order to simplify the reconstruction procedure, the grid is treated as locally uniform for each cell independent of its grid level $\ell$, because the neighbors cells at a coarser level $\ell-1$ 
can be virtually refined in order to allow for the reconstruction procedure on locally uniform meshes detailed in section \ref{ssec.reconstruction}. We also note that our AMR algorithm is endowed with a time-accurate local time stepping (LTS) feature, see \cite{AMR3DCL} for details.

\section{Numerical results} 
\label{sec.results}

In this Section we present a large set of numerical test cases in order to show the accuracy, robustness and efficiency of the presented $P_NP_M$ family of schemes equipped with the \textit{a posteriori} subcell finite volume limiter.

In order to cover a wide variety of physical phenomena we have applied our schemes to three sets of equations of relevance in fluid-dynamical applications, 
namely the Euler equations of compressible hydrodynamics (HD), the magnetohydrodynamics equations (MHD), and the special relativistic magnetohydrodynamics equations (RMHD).  

In particular, for any set of equations we have selected both a smooth test case, to show the order of convergence of our schemes (up to order six),
and some problems containing strong discontinuities going from logically one-dimensional Riemann problems to classical challenging two-dimensional benchmarks, such as the Sedov explosion problem, the Double Mach Reflection problem, the MHD rotor problem, the RMHD blast wave, as well as the MHD \& RMHD Orszag-Tang vortex problems.
The presence of discontinuities allows to prove the robustness and resolution of our \textit{a posteriori} subcell limiting strategy.

Moreover, the results obtained with the intermediate $P_NP_M$ schemes (i.e. $N\ne 0$ and $M> N$) are compared with the pure DG approach (i.e. $N=M$) in order to show their gain in terms of computational efficiency, while maintaining a similar resolution.  
We also compare numerical results on AMR meshes against results obtained on fine uniform Cartesian meshes, 
demonstrating both the robustness of our schemes on adaptive meshes and the obtained savings in computational time.

\subsection{Euler equations of gasdynamics}

The first set of hyperbolic equations that we consider is given by the homogeneous Euler equations of compressible gasdynamics
that can be cast in form~\eqref{eq.generalform} by choosing   
\begin{equation}
\label{eulerTerms}
\Q = \left( \begin{array}{c} \rho   \\ \rho u  \\ \rho v \\ \rho w \\ \rho E \end{array} \right)\!, \quad
\mathbf{F} = \left( \begin{array}{ccc}  \rho u       & \rho v        \\ 
\rho u^2 + p & \rho u v          \\
\rho u v     & \rho v^2 + p      \\ 
\rho u w     & \rho v w          \\ 
u(\rho E + p) & v(\rho E + p)   
\end{array} \right). 
\end{equation}
The vector of conserved variables $\Q$ involves the fluid density $\rho$, the momentum density vector $\rho \v=(\rho u, \rho v)$ and the total energy density $\rho E$. 
The fluid pressure $p$ is related to the conserved quantities $\Q$ using the equation of state for an ideal 
gas    
\begin{equation}
\label{eqn.eos} 
p = (\gamma-1)\left(\rho E - \frac{1}{2} \rho \mathbf{v}^2 \right)\!, 
\end{equation}
where $\gamma$ is the ratio of specific heats so that the speed of sound takes the form $c=\sqrt{\frac{\gamma p}{\rho}}$.

\subsubsection{Isentropic vortex}

First of all, in order to verify the order of convergence of the proposed $\PN{N}{M}$ schemes, we consider 
a smooth isentropic vortex flow according to~\cite{HuShuVortex1999}. 
The computational domain is given by the square $\Omega=[0,10]\times[0,10]$ with periodic boundary conditions set everywhere. 
For the initial conditions we consider a homogeneous background field $\Q_0=(\rho,u,v,p)=(1,1,1,1)$ traveling 
with a constant velocity $\v_c = (1,1)$ and we superimpose on this field some perturbations for density and pressure of the following form 
\begin{equation}
\label{rhopressDelta}
\delta \rho = (1+\delta T)^{\frac{1}{\gamma-1}}-1, \quad \delta p = (1+\delta T)^{\frac{\gamma}{\gamma-1}}-1, 
\end{equation}
with the temperature fluctuation \[\delta T = -\frac{(\gamma-1)\epsilon^2}{8\gamma\pi^2}e^{1-r^2}\] and the vortex strength $\epsilon=5$.
The velocity field is also affected by the following perturbations
\begin{equation}
\label{ShuVortDelta}
\left(\begin{array}{c} \delta u \\ \delta v  \end{array}\right) = \frac{\epsilon}{2\pi}e^{\frac{1-r^2}{2}} \left(\begin{array}{c} -(y-5) \\ \phantom{-}(x-5)  \end{array}\right).
\end{equation}
The initial condition is thus given by $\Q = \Q_0 + \mathbf{\delta} \mathbf{Q}$.
The exact solution $\Q_e$ at the final time $t_f$ can be simply computed as the time-shifted initial condition, i.e. $\Q_e(x, t_f) = \Q(x - \v_c t_f, 0)$.

In Table~\ref{tab.orderOfconvergenceFV_shu}, we report the convergence rates from second up to sixth order of accuracy for the smooth vortex test problem 
run on a sequence of successively refined meshes up to the final time $t_f=1.0$. 
The optimal order of accuracy is achieved for the hybrid schemes $\PN{N}{M}$ with $M>N$ and for the pure DG schemes with $N=M$. 

\begin{landscape} 
	\begin{table*}
		\caption{Numerical convergence table for general $P_NP_M$ schemes for the isentropic vortex problem. The error norms refer to the variable $\rho$ at time $t_f=1.0$ in $L_2$ norm.} 
		\label{tab.orderOfconvergenceFV_shu}
		\begin{center} 	
			\begin{tabular}{|c|ccc|ccc|ccc|ccc|ccc|ccc|} 
				\hline 
				h                        & \CPU      & $L_2$             & \OLdue        &   \CPU      & $L_2$             & \OLdue       &     \CPU      & $L_2$             & \OLdue       &      \CPU      & $L_2$             & \OLdue       &    \CPU      & $L_2$             & \OLdue       &        \CPU     & $L_2$            & \OLdue          \\[1pt]
				\hline                                                                                                                                                                                                                                                                                                                                   
				$\mathbf{\mathcal{O}2}$      &               \mcolt{0}{1}                    &                \mcolt{1}{1}                    &               &                  &               &                &                  &               &              &                   &              &                 &                  &                   \\
				\smm 5.0E-02       \smmm &\sm 60 \sm & \sm 1.0E-02   \sm & \smm     \smm & \sm 146\sm & \sm 1.2E-04   \sm & \smm     \smm &   \sm     \sm & \sm          \sm & \smm     \smm &    \sm     \sm & \sm          \sm & \smm     \smm &  \sm     \sm & \sm          \sm & \smm     \smm &      \sm    \sm & \sm          \sm & \smm     \smm    \\
				\smm 4.0E-02       \smmm &\sm 91 \sm & \sm 7.4E-03   \sm & \smm 1.4 \smm & \sm 291\sm & \sm 7.9E-05   \sm & \smm 1.6 \smm &   \sm     \sm & \sm          \sm & \smm     \smm &    \sm     \sm & \sm          \sm & \smm     \smm &  \sm     \sm & \sm          \sm & \smm     \smm &      \sm    \sm & \sm          \sm & \smm     \smm    \\
				\smm 3.3E-02       \smmm &\sm 159\sm & \sm 5.7E-03   \sm & \smm 1.4 \smm & \sm 450\sm & \sm 5.4E-05   \sm & \smm 2.0 \smm &   \sm     \sm & \sm          \sm & \smm     \smm &    \sm     \sm & \sm          \sm & \smm     \smm &  \sm     \sm & \sm          \sm & \smm     \smm &      \sm    \sm & \sm          \sm & \smm     \smm    \\
				\smm 2.5E-02       \smmm &\sm 367\sm & \sm 3.8E-03   \sm & \smm 1.4 \smm & \sm1132\sm & \sm 3.1E-05   \sm & \smm 2.0 \smm &   \sm     \sm & \sm          \sm & \smm     \smm &    \sm     \sm & \sm          \sm & \smm     \smm &  \sm     \sm & \sm          \sm & \smm     \smm &      \sm    \sm & \sm          \sm & \smm     \smm    \\[1pt]
				\hline                                                                                                                                                                                                                                                                                                                               
				$\mathbf{\mathcal{O}3}$      &               \mcolt{0}{2}                    &                \mcolt{1}{2}                    &                  \mcolt{2}{2}                    &               &                   &               &             &                   &               &                 &                   &                    \\
				\smm 6.7e-02       \smmm &\sm 53 \sm & \sm 3.7e-04   \sm & \smm     \smm & \sm 180\sm & \sm 4.2e-06   \sm & \smm     \smm &   \sm 349\sm & \sm 1.7e-05   \sm & \smm     \smm &    \sm    \sm & \sm           \sm & \smm     \smm &  \sm    \sm & \sm           \sm & \smm     \smm &      \sm    \sm & \sm           \sm & \smm     \smm    \\
				\smm 5.0e-02       \smmm &\sm 120\sm & \sm 1.6e-04   \sm & \smm 2.9 \smm & \sm 415\sm & \sm 1.8e-06   \sm & \smm 3.0 \smm &   \sm 819\sm & \sm 7.9e-06   \sm & \smm 2.6 \smm &    \sm    \sm & \sm           \sm & \smm     \smm &  \sm    \sm & \sm           \sm & \smm     \smm &      \sm    \sm & \sm           \sm & \smm     \smm    \\
				\smm 4.0e-02       \smmm &\sm 245\sm & \sm 8.2e-05   \sm & \smm 2.9 \smm & \sm 791\sm & \sm 9.0e-07   \sm & \smm 3.0 \smm &   \sm 1615\sm &\sm 4.3e-06   \sm & \smm 2.6 \smm &    \sm    \sm & \sm           \sm & \smm     \smm &  \sm    \sm & \sm           \sm & \smm     \smm &      \sm    \sm & \sm           \sm & \smm     \smm    \\
				\smm 3.3e-02       \smmm &\sm 417\sm & \sm 4.7e-05   \sm & \smm 2.9 \smm & \sm1377\sm & \sm 5.2e-07   \sm & \smm 3.0 \smm &   \sm 2734\sm &\sm 2.6e-06   \sm & \smm 2.7 \smm &    \sm    \sm & \sm           \sm & \smm     \smm &  \sm    \sm & \sm           \sm & \smm     \smm &      \sm    \sm & \sm           \sm & \smm     \smm    \\[1pt]
				\hline                                                                                                                                                                                                                                                                                                                               
				$\mathbf{\mathcal{O}4}$      &               \mcolt{0}{3}                    &                \mcolt{1}{3}                    &                  \mcolt{2}{3}                    &                   \mcolt{3}{3}                    &             &                   &               &                 &                   &                    \\
				\smm 1.3e-01       \smmm &\sm 2.9\sm & \sm 4.0e-03   \sm & \smm     \smm & \sm 6.5\sm & \sm 2.2e-04   \sm & \smm     \smm &   \sm 15 \sm & \sm 2.2e-05   \sm & \smm     \smm &    \sm 23 \sm & \sm 6.0e-06   \sm & \smm     \smm &  \sm    \sm & \sm           \sm & \smm     \smm &      \sm    \sm & \sm           \sm & \smm     \smm    \\
				\smm 8.3e-02       \smmm &\sm 8.9\sm & \sm 6.0e-04   \sm & \smm 4.6 \smm & \sm 22 \sm & \sm 3.8e-05   \sm & \smm 4.3 \smm &   \sm 50 \sm & \sm 3.9e-06   \sm & \smm 4.2 \smm &    \sm 79 \sm & \sm 9.1e-07   \sm & \smm 4.6 \smm &  \sm    \sm & \sm           \sm & \smm     \smm &      \sm    \sm & \sm           \sm & \smm     \smm    \\
				\smm 6.3e-02       \smmm &\sm 20 \sm & \sm 1.5e-04   \sm & \smm 4.8 \smm & \sm 51 \sm & \sm 1.1e-05   \sm & \smm 4.1 \smm &   \sm 115\sm & \sm 1.1e-06   \sm & \smm 4.4 \smm &    \sm 183\sm & \sm 2.7e-07   \sm & \smm 4.2 \smm &  \sm    \sm & \sm           \sm & \smm     \smm &      \sm    \sm & \sm           \sm & \smm     \smm    \\
				\smm 5.0e-02       \smmm &\sm 38 \sm & \sm 4.9e-05   \sm & \smm 5.0 \smm & \sm 101\sm & \sm 4.6e-06   \sm & \smm 4.0 \smm &   \sm 235\sm & \sm 4.1e-07   \sm & \smm 4.4 \smm &    \sm 380\sm & \sm 1.1e-07   \sm & \smm 4.1 \smm &  \sm    \sm & \sm           \sm & \smm     \smm &      \sm    \sm & \sm           \sm & \smm     \smm    \\[1pt]
				\hline                                                                                                                                                                                                                                                                                                                               
				$\mathbf{\mathcal{O}5}$      &               \mcolt{0}{4}                    &                \mcolt{1}{4}                    &                  \mcolt{2}{4}                    &                   \mcolt{3}{4}                    &                 \mcolt{4}{4}                    &                 &                  &                     \\
				\smm 8.3e-02       \smmm &\sm 22 \sm & \sm 4.2e-04   \sm & \smm     \smm & \sm 47 \sm & \sm 1.7e-06   \sm & \smm     \smm &   \sm 96 \sm & \sm 2.9e-06   \sm & \smm     \smm &    \sm 181 \sm & \sm 4.8e-08   \sm & \smm    \smm &  \sm 242  \sm & \sm 1.1e-07   \sm & \smm    \smm &      \sm    \sm & \sm          \sm & \smm    \smm      \\
				\smm 7.7e-02       \smmm &\sm 28 \sm & \sm 2.7e-04   \sm & \smm 4.6 \smm & \sm 57 \sm & \sm 1.1e-06   \sm & \smm 5.3 \smm &   \sm 131\sm & \sm 2.0e-06   \sm & \smm 4.8 \smm &    \sm 229 \sm & \sm 3.3e-08   \sm & \smm 4.8\smm &  \sm 314  \sm & \sm 8.0e-08   \sm & \smm 4.2\smm &      \sm    \sm & \sm          \sm & \smm    \smm      \\
				\smm 7.1e-02       \smmm &\sm 34 \sm & \sm 2.0e-04   \sm & \smm 4.7 \smm & \sm 70 \sm & \sm 7.4e-07   \sm & \smm 5.3 \smm &   \sm 146\sm & \sm 1.5e-06   \sm & \smm 4.6 \smm &    \sm 294 \sm & \sm 2.3e-08   \sm & \smm 4.6\smm &  \sm 372 \sm & \sm 5.8e-08   \sm & \smm 4.3\smm &      \sm    \sm & \sm          \sm & \smm    \smm      \\
				\smm 6.7e-02       \smmm &\sm 41 \sm & \sm 1.4e-04   \sm & \smm 4.7 \smm & \sm 91 \sm & \sm 5.1e-07   \sm & \smm 5.3 \smm &   \sm 181\sm & \sm 1.1e-06   \sm & \smm 4.1 \smm &    \sm 346 \sm & \sm 1.7e-08   \sm & \smm 4.1\smm &  \sm 461 \sm & \sm 4.3e-08   \sm & \smm 4.2\smm &      \sm    \sm & \sm          \sm & \smm    \smm      \\[1pt]
				\hline                                                                                                                                                                                                                                                                                                                               
				$\mathbf{\mathcal{O}6}$      &               \mcolt{0}{5}                    &                \mcolt{1}{5}                    &                  \mcolt{2}{5}                    &                   \mcolt{3}{5}                    &                 \mcolt{4}{5}                    &                     \mcolt{5}{5}                       \\
				\smm 2.0e-01       \smmm &\sm 7  \sm & \sm 1.1e-02   \sm & \smm     \smm & \sm 8  \sm & \sm 2.4e-04   \sm & \smm     \smm &   \sm 15 \sm & \sm 1.3e-05   \sm & \smm     \smm &    \sm 26 \sm & \sm 1.3e-05  \sm & \smm     \smm &  \sm 43  \sm & \sm 7.2e-07   \sm & \smm     \smm &      \sm 57 \sm & \sm 6.7e-07   \sm & \smm     \smm    \\			
				\smm 1.7e-01       \smmm &\sm 11 \sm & \sm 4.4e-03   \sm & \smm 5.1 \smm & \sm 12 \sm & \sm 7.1e-05   \sm & \smm 6.6 \smm &   \sm 24 \sm & \sm 5.0e-06   \sm & \smm 5.6 \smm &    \sm 44 \sm & \sm 4.1e-06  \sm & \smm 6.5 \smm &  \sm 69  \sm & \sm 2.6e-07   \sm & \smm 5.6 \smm &      \sm 98  \sm & \sm 1.8e-07   \sm & \smm 7.3 \smm    \\
				\smm 1.4e-01       \smmm &\sm 16 \sm & \sm 1.8e-03   \sm & \smm 5.7 \smm & \sm 19 \sm & \sm 2.4e-05   \sm & \smm 7.0 \smm &   \sm 38 \sm & \sm 2.1e-06   \sm & \smm 5.5 \smm &    \sm 68 \sm & \sm 1.4e-06  \sm & \smm 7.0 \smm &  \sm 110 \sm & \sm 1.1e-07   \sm & \smm 5.8 \smm &      \sm 153 \sm & \sm 5.2e-08   \sm & \smm 8.0 \smm    \\
				\smm 1.3e-01       \smmm &\sm 23 \sm & \sm 8.3e-04   \sm & \smm 5.6 \smm & \sm 29 \sm & \sm 9.2e-06   \sm & \smm 7.2 \smm &   \sm 55 \sm & \sm 1.0e-06   \sm & \smm 5.5 \smm &    \sm 97 \sm & \sm 5.4e-07  \sm & \smm 7.2 \smm &  \sm 157 \sm & \sm 4.5e-08   \sm & \smm 5.9 \smm &      \sm 219 \sm & \sm 1.9e-08   \sm & \smm 7.5 \smm    \\
				\hline 
			\end{tabular}		
		\end{center}
	\end{table*}
\end{landscape}


\subsubsection{The Sod shock tube problem}
\begin{figure}[p]
	\centering
	\includegraphics[width=0.48\linewidth]{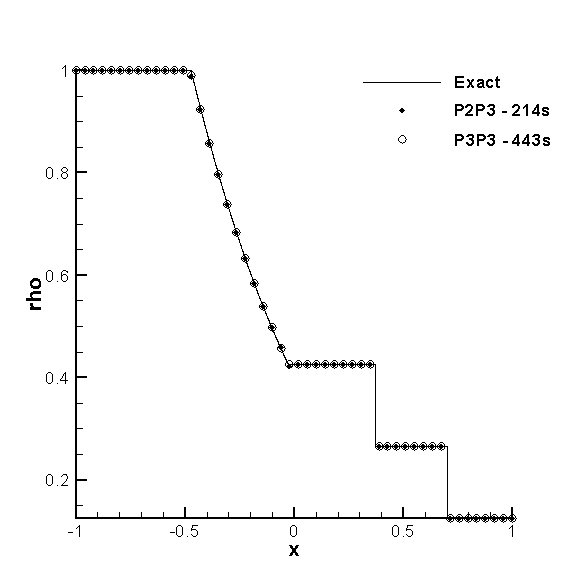}
	\includegraphics[width=0.48\linewidth]{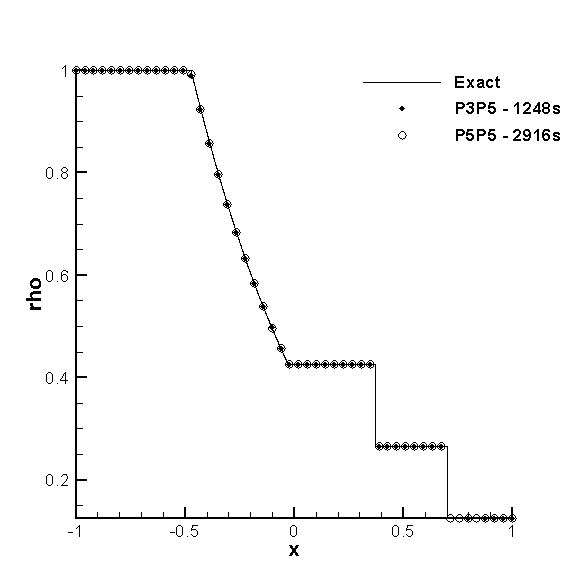}
	\caption{Sod shock tube problem at $t_f = 0.4$. Left: we compare our numerical results obtained with two fourth order methods, 
		namely the $P_2P_3$ and the $P_3P_3$ with the exact solution. Note that the two schemes show a similar resolution but the $P_2P_3$ 
		is twice faster than the $P_3P_3$,
	Right: we compare our numerical results obtained with two sixth order methods, namely the $P_3P_5$ and the $P_5P_5$ with the exact solution. 
	Note again that the two schemes show a similar resolution but the $P_3P_5$ is 2.33 times faster than the $P_5P_5$. 
	}
	\label{fig:Sod_P3_vsExact}
\end{figure}
\begin{figure}[p]
	\centering
	\includegraphics[width=0.48\linewidth]{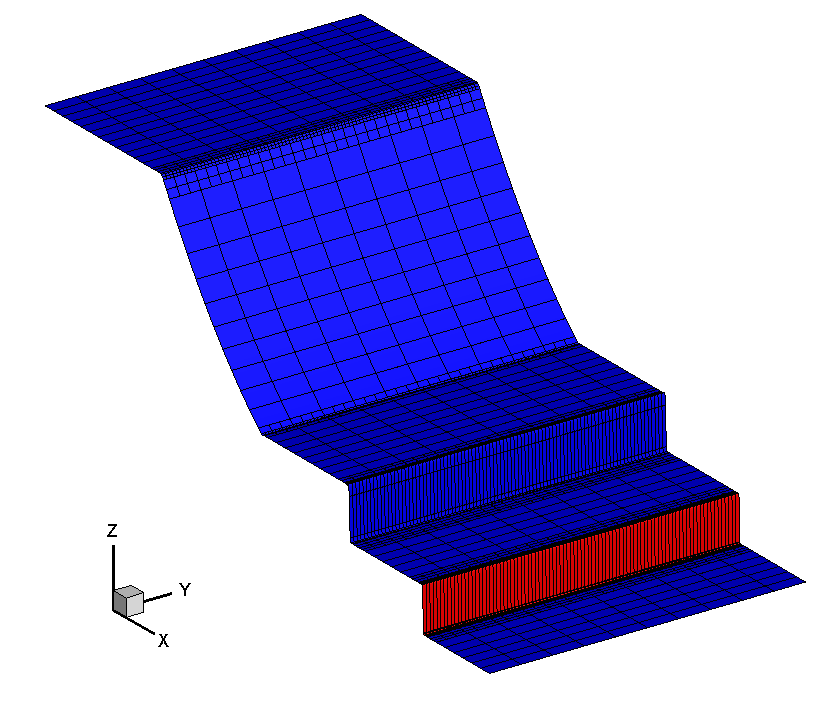}
	\includegraphics[width=0.48\linewidth]{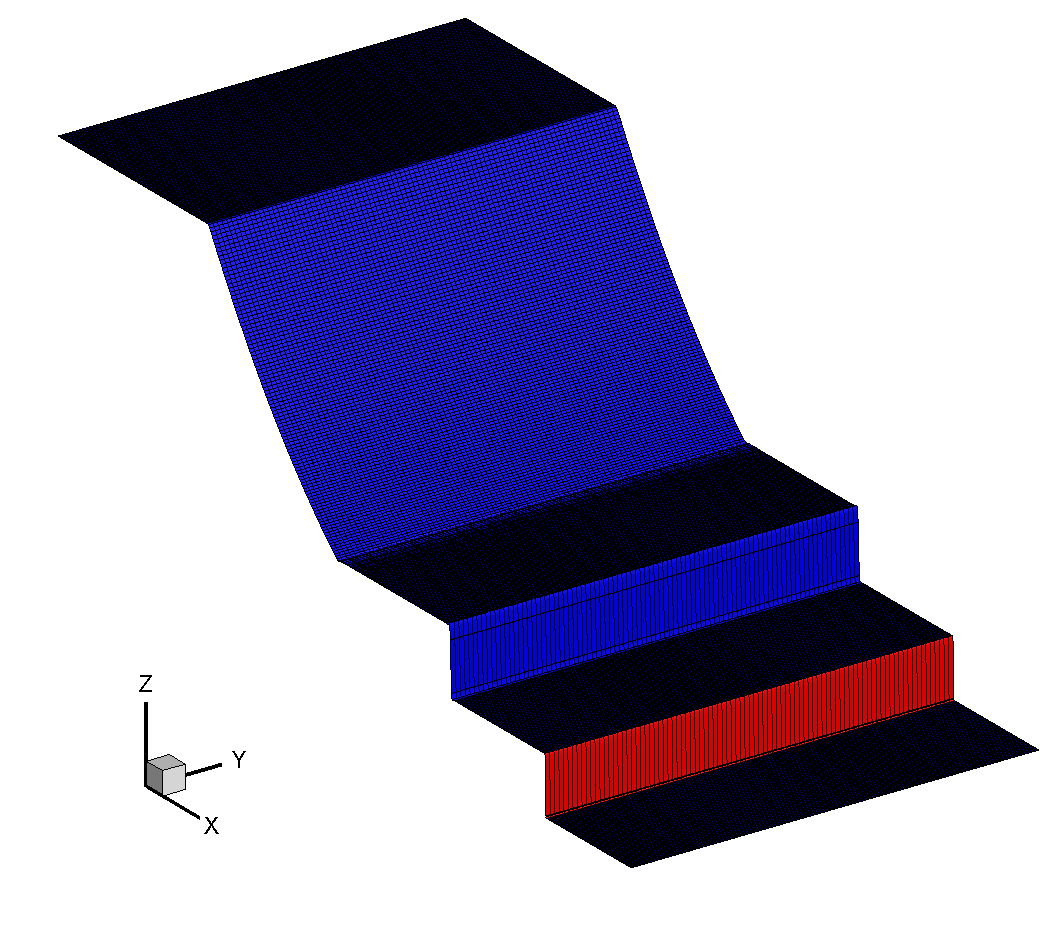}
	\caption{Sod shock tube problem at $t_f = 0.4$ solved with our $P_3P_5$ sixth order scheme. 
		We show in red the cells on which the limiter is activated and in blue the unlimited cells. 
		The left panel is obtained with an initial grid of $50\times 10$ elements and $\ell_{\max}=2$ levels of refinement with $\rRef=3$ in a total computational time of $1248$s.  
		The right image is obtained by using a fine uniform grid of $450\times90$ elements corresponding to the finest AMR grid level in a total computational time of $3670$s.  
		In both the cases the limiter perfectly activates only where the shock wave is located.}
	\label{fig:Sod_P3P5_50x10_Ref2_limiter3d}
\end{figure}

The Sod shock tube problem in 2D can be seen as a multidimensional extension of the classical Sod test case in 1D~\cite{ToroBook},
which allows to verify the robustness and the resolution capacity of the employed numerical method
on a rarefaction wave, a contact discontinuity and a shock at the same time, indeed the three waves are originated by the discontinuous initial condition.

Here, we consider as computational domain a square of dimension $[-1,1]\times[-1,1]$ covered with a uniform mesh of $50\times10$ control volumes, 
and the initial condition is composed of two  different states, separated by a discontinuity at $x_d=0$ 
\be
\begin{cases} 
	\rho_L = 1, \quad \mathbf{v}_L = 0, \quad p_L = 1, \quad  & x \leq x_d \\
	\rho_R = 0.125, \quad \mathbf{v}_R = 0, \quad  p_R = 0.1, \quad & x  > x_d.
\end{cases}
\ee
The final time is chosen to be $t_f=0.4$, so that the shock wave does not cross the external boundary of the domain, where wall boundary conditions are set.  
The algorithm for the calculation of the exact solution of this Riemann problem is given in \cite{toro-book}. 

We have run this problem with two fourth order methods, namely the $\PN{2}{3}$ and the $\PN{3}{3}$ schemes, and two sixth order methods, 
namely the $\PN{3}{5}$ and the $\PN{5}{5}$ schemes, 
equipped with the \textit{a posteriori} subcell TVD FV limiter and employing the Rusanov flux both in the main $\PN{N}{M}$ scheme and at the limiter stage.
The agreement of our numerical results with the exact solution is perfect and the hybrid schemes (i.e. $M<N$) are faster than the pure DG schemes (N=M), see Figure~\ref{fig:Sod_P3_vsExact}. 

Moreover, in Figure~\ref{fig:Sod_P3P5_50x10_Ref2_limiter3d}, one can see that the limiter activates exactly where the shock discontinuity is located also when 
the adaptive mesh refinement technique is employed.

\subsubsection{The Lax shock tube problem}
\begin{figure}[b]
	\centering
	\includegraphics[width=0.3\linewidth]{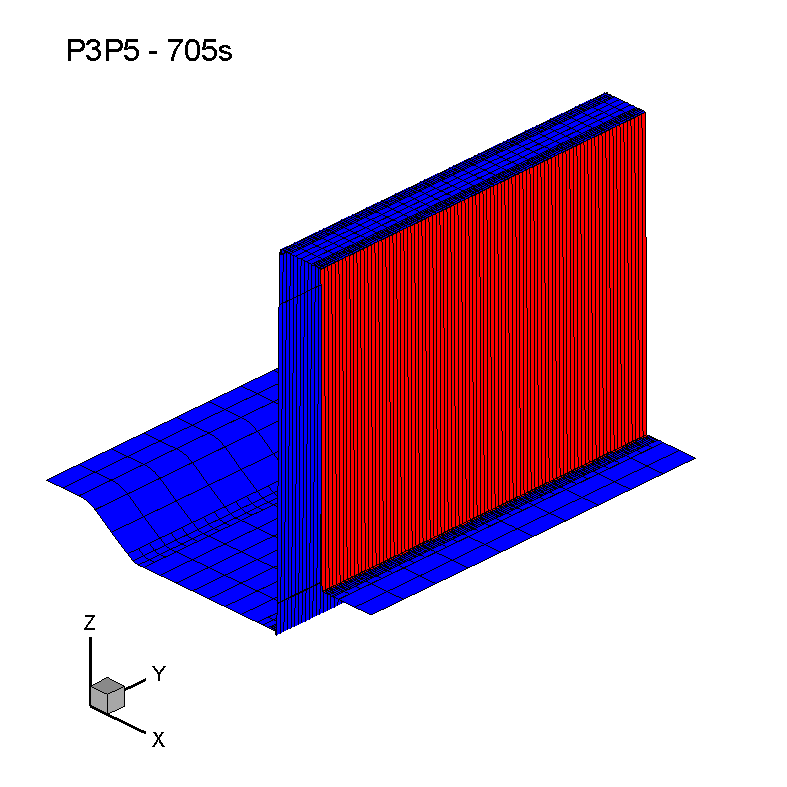}
	\includegraphics[width=0.3\linewidth]{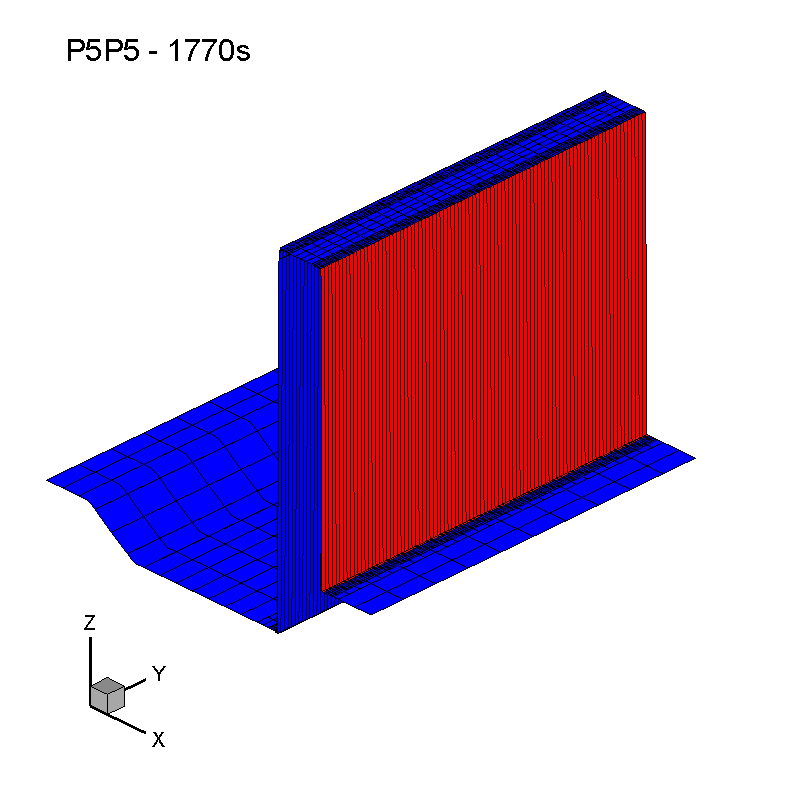}
	\includegraphics[width=0.35\linewidth]{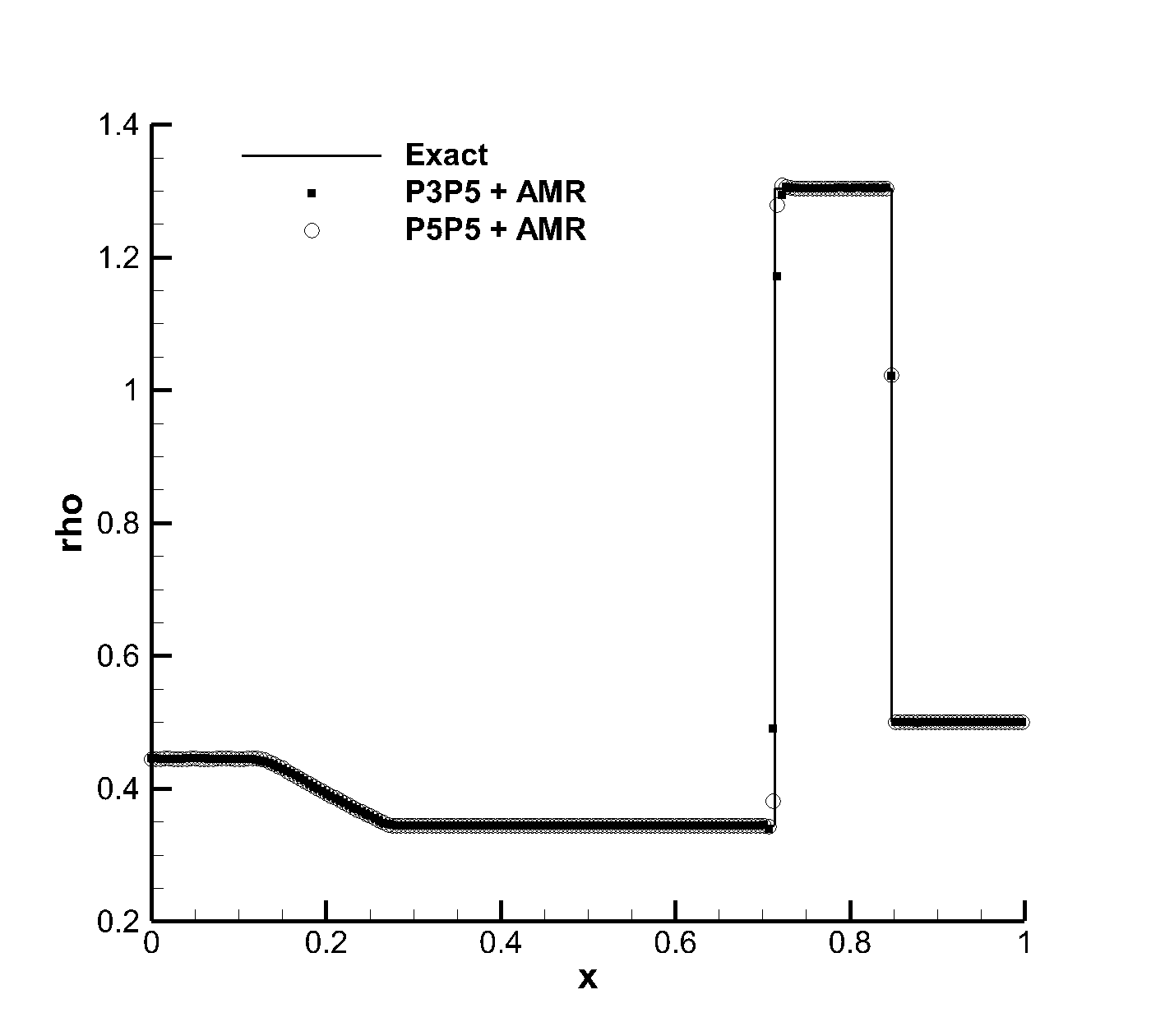}
	\caption{Lax shock tube problem at $t_f=0.14$, obtained with our sixth order schemes, namely the $P_3P_5$ and $P_5P_5$ schemes. 
		Left: we draw the density profile on the $z$ axis and we colour in red the cells on which the limiter is activated and in blue the unlimited cells. 
		Right: we compare our numerical results with the exact solution.}
	\label{fig:laxvsref}
\end{figure}

The Lax shock tube problem, introduced for the first time in~\cite{lax2}, is another classical benchmark for high order methods for the solution of the Euler equations.
The computational domain is the square $[-1,1]\times[-1,1]$ and the initial condition is composed of two  different states, separated by a discontinuity located at $x_d=0$ 
\be
\begin{cases} 
	\rho_L = 0.445, \quad \mathbf{v}_L = 0.698, \quad p_L = 3.528, \quad  & x \leq x_d \\
	\rho_R = 0.5, \quad \mathbf{v}_R = 0, \quad  p_R = 0.571, \quad & x  > x_d .
\end{cases}
\ee

In this case, we have covered our computational domain with a very coarse mesh of $20\times10$ elements activating the AMR procedure with  $\ell_{\max}=2$ levels of refinement and $\rRef=3$. 
In Figure~\ref{fig:laxvsref}, we present the results obtained with two sixth order methods, namely the $\PN{3}{5}$ 
and the $\PN{5}{5}$ schemes, used together with the HLLEM~\cite{HLLEM,Dumbser2015} numerical flux.
Both the numerical results perfectly agree with the reference solution and the hybrid scheme is $2.5$ times faster than the pure DG scheme. 
Also the limiter activates exactly only at the shock location and, 
due to the subcell resolution, it does not affect the quality of the profile which is sharply captured even on a very coarse mesh.

\subsubsection{The Shu-Osher shock tube problem}
\begin{figure}[!h]
	\centering
	\includegraphics[width=0.38\linewidth]{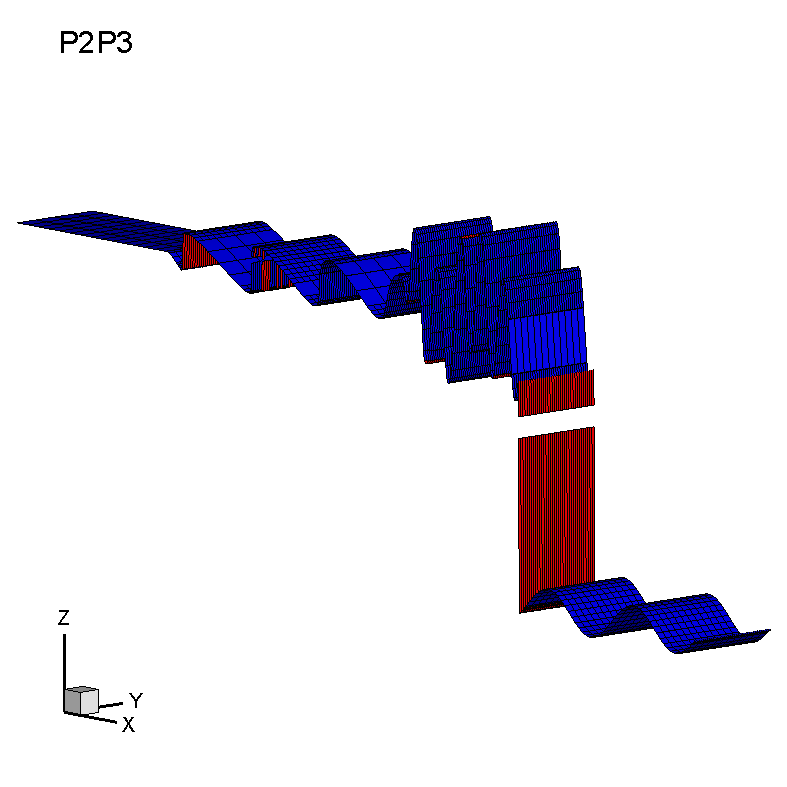}\qquad\qquad
	\includegraphics[width=0.38\linewidth]{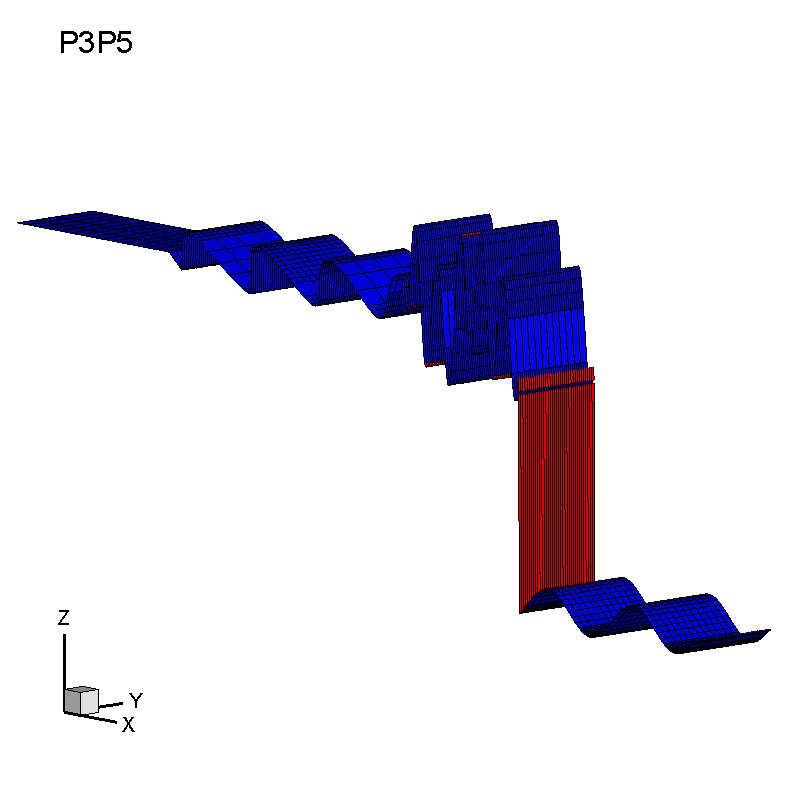}\\
	\includegraphics[width=0.38\linewidth]{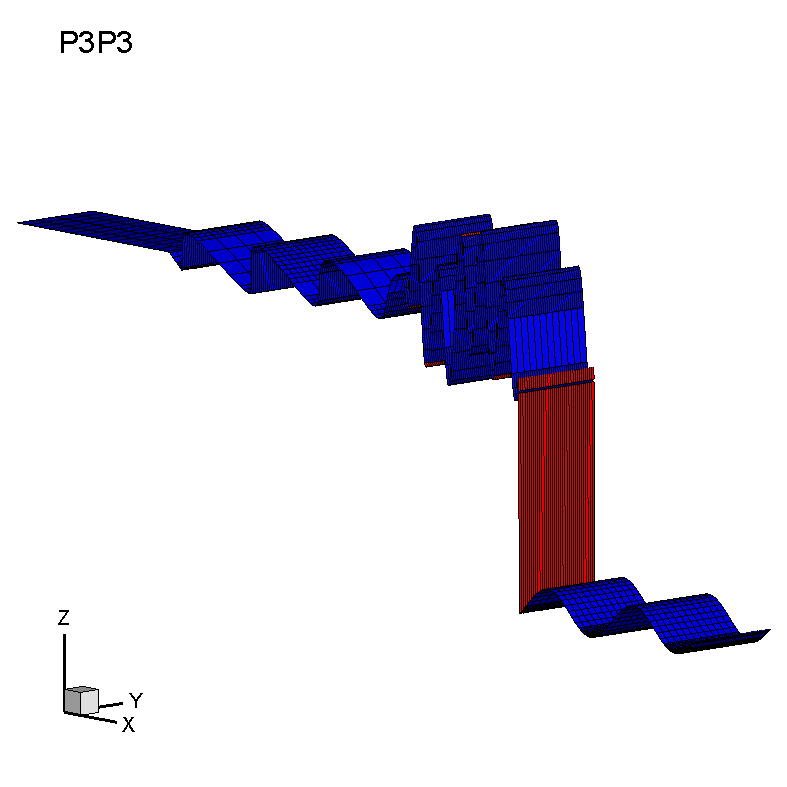}\qquad\qquad 
	\includegraphics[width=0.38\linewidth]{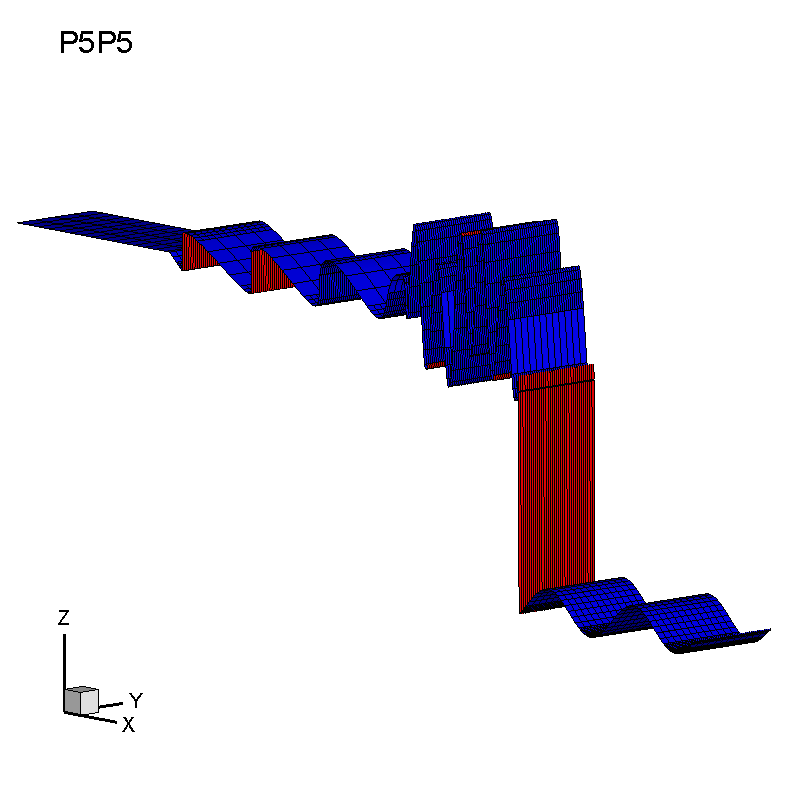}\\
	\includegraphics[width=0.42\linewidth]{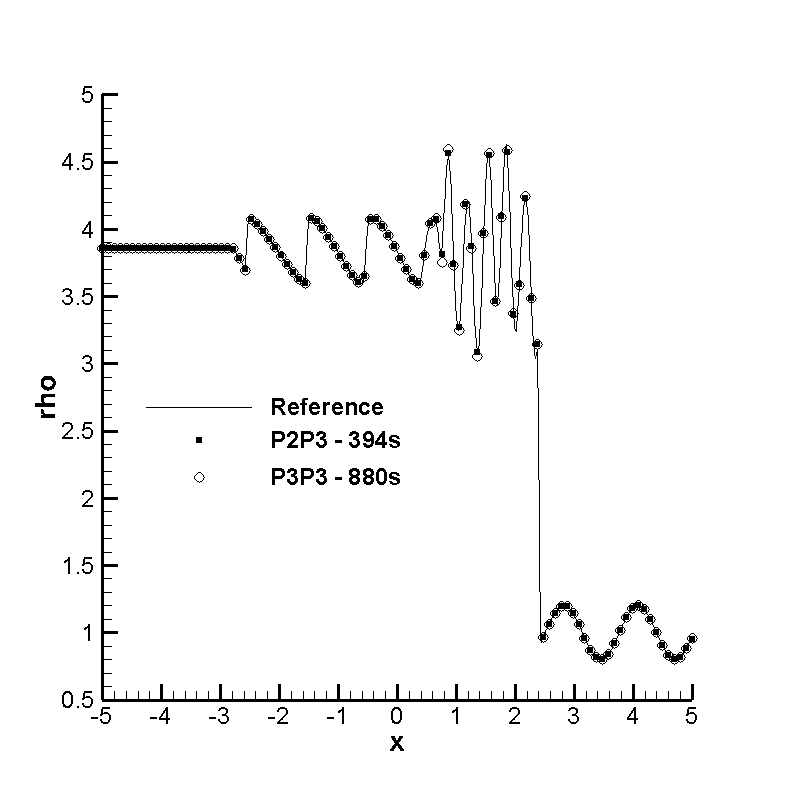}\qquad\qquad
	\includegraphics[width=0.42\linewidth]{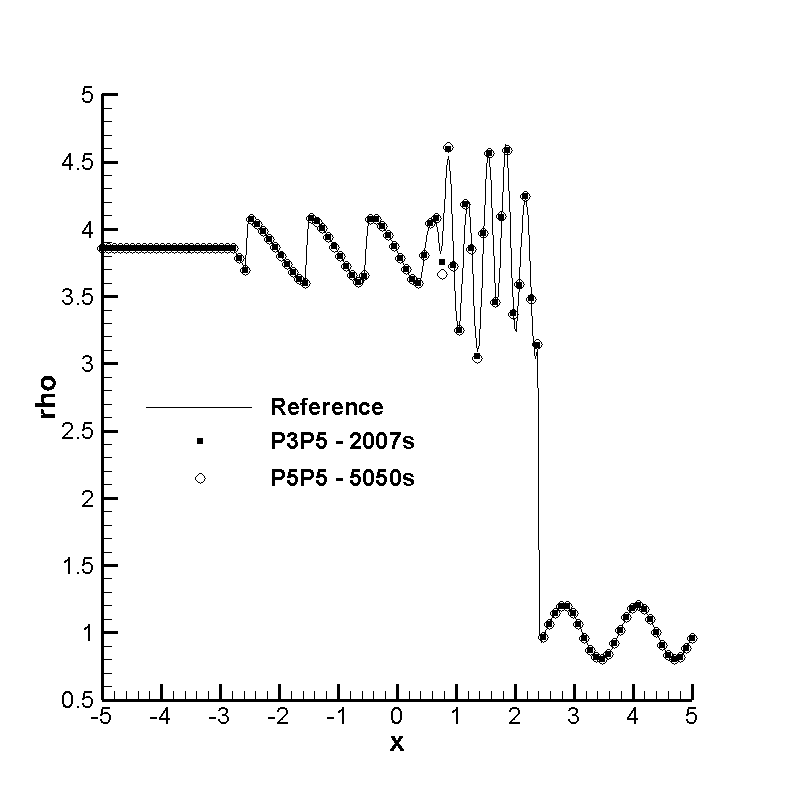}\\
	\caption{Shu-Osher shock tube problem at $t_f=1.8$ solved with our fourth order scheme, namely the $P_2P_3$ and $P_3P_3$ schemes (first column), 
		and our sixth order schemes, namely the $P_3P_5$ and $P_5P_5$ schemes (second column).	
		In the first two rows we plot the value of the density on the $z$ axis and we depict in red the cell where the limiter is activated. 
		In the third row our results are compared with a reference solution obtained with a WENO FV scheme on a very fine mesh. 
		Moreover the computational time required by the $P_NP_M$ schemes with $N<M$ (respectively $394$s and $2007$s) are significantly shorter than those required by the pure DG schemes ($N=M$) 
		(respectively $880$s and $5050$s); nevertheless all the results show an excellent agreement with the reference solution.
	}
	\label{fig.shuosher}
\end{figure}

The Shu-Osher problem was first introduced in~\cite{shuosher1} 
and it allows to check the capability of our new scheme to deal simultaneously with \textit{physical oscillations} and \textit{shock waves} appearing at the same time during the simulation.   
It consists of a one-dimensional Mach $3$ shock front interacting with a sinusoidal density disturbance 
that generates a combination of discontinuities and smooth structures, 
whose entropy fluctuations are amplified when passing through the shock. 

We discretize our computational domain $\Omega = [-5,5]\times [0,1]$ with an AMR grid with $64\times4$ elements 
on the coarsest level and a maximum refinement level $\lmax=2$ with $\rRef=3$.
The initial conditions are given by
\be
\begin{cases}
	\rho=3.8571, \quad u= 2.6294, \quad v= 0, \quad p= 10.333, \quad & x < 4, \\
	\rho= 1.0+0.2\sin(5x), \quad u= 0,\quad  v=0, \quad p=1, \quad & x \ge 4,
\end{cases}
\ee 
and the simulations run up to the final time $t_f=1.8$.
For this test case, we employ the HLL Riemann solver, 
and the TVD \textit{a posteriori} subcell finite volume limiter.  

We then compare the results obtained with two pure DG schemes $\PN{3}{3}$ and $\PN{5}{5}$, and the hybrid schemes $\PN{2}{3}$ and $\PN{3}{5}$ which have, 
respectively, the same order of accuracy. As shown in Figure~\ref{fig.shuosher}, the  
methods are accurate, and robust thanks to the employed limiter strategy, 
and the intermediate schemes $\PN{N}{M}$ with $M>N$ are computational more efficient than pure DG schemes.

%
%
%
%


\subsubsection{Sedov problem}
This test problem is widespread in literature~\cite{SedovExact,LoubereSedov3D,GaburroAREPO} and it describes the evolution of a blast wave 
that is generated at the origin $\mathbf{O}=(x,y)=(0,0)$ of the computational domain $\Omega(0)=[0,1.2]\times[0,1.2]$. 
An exact solution based on self-similarity arguments is available from~\cite{Sedov} and the fluid is assumed to be an ideal gas with $\gamma=1.4$, 
which is initially at rest and assigned with a uniform density $\rho_0=1$. 
The initial pressure is $p_0=10^{-6}$ everywhere except in the cell $V_{or}$ containing the origin $\mathbf{O}$ where it is given by
\[
p_{or} = (\gamma-1) \rho_0 \frac{ E_\text{tot}}{ |V_\text{or}|}
\]
\begin{figure}
	\centering
	\includegraphics[width=0.33\linewidth]{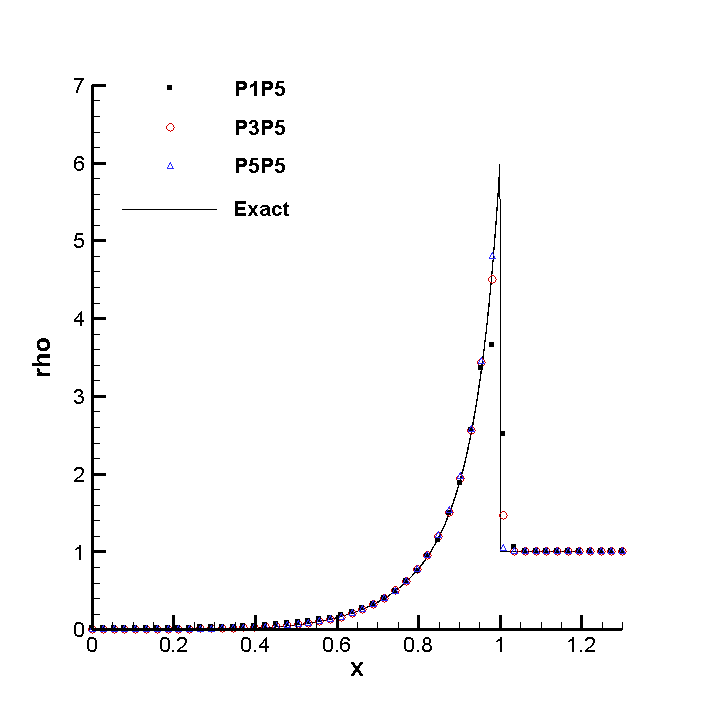}%
	\includegraphics[width=0.33\linewidth]{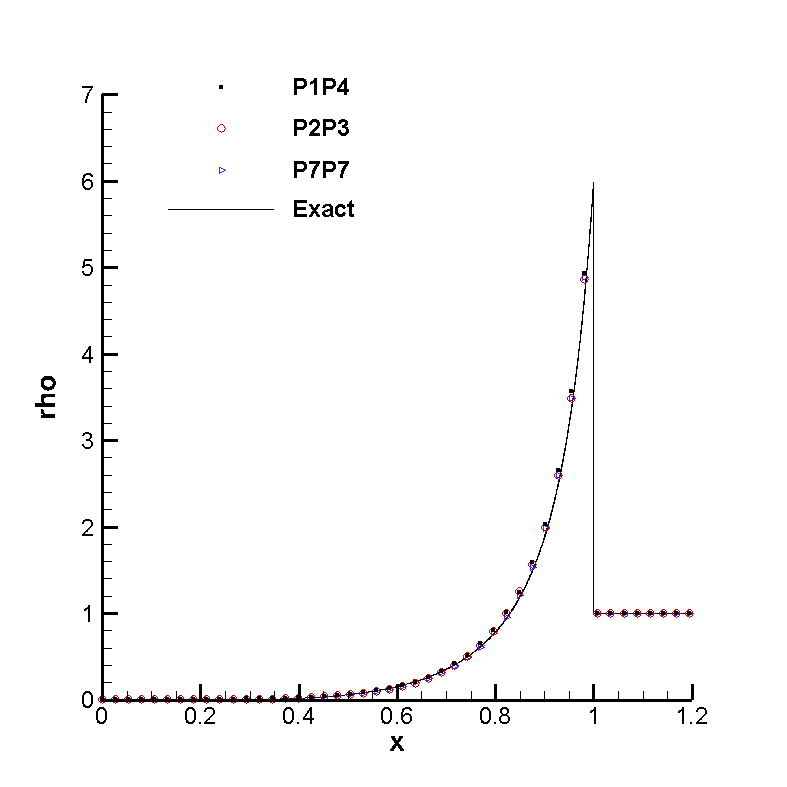}%
	\includegraphics[width=0.33\linewidth]{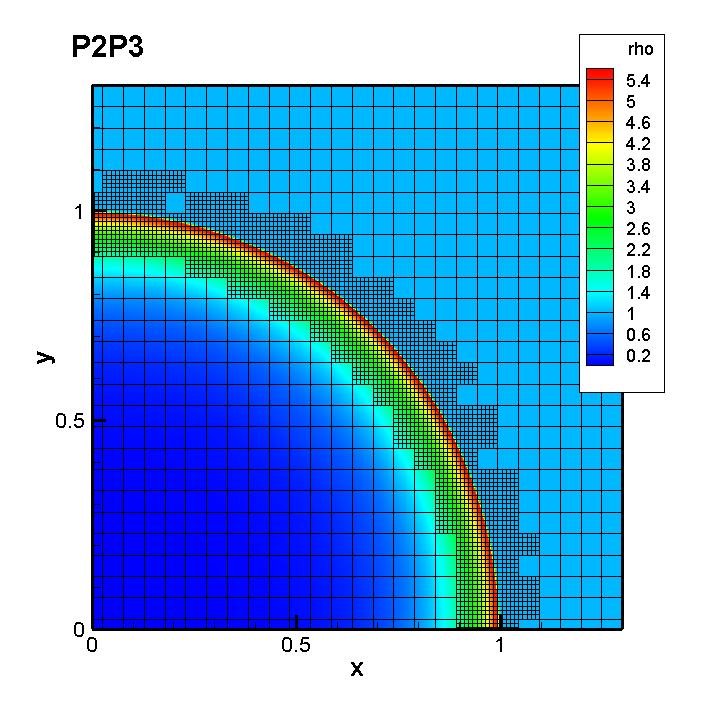}
	\caption{Sedov test problem. Comparison of numerical density obtained with a selection of our high order methods (left, middle)
		 and density contours obtained with the hybrid $\PN{2}{3}$ third order scheme on an AMR grid.}
	\label{fig:sedovRho}
\end{figure}
\begin{figure}
	\centering
	\includegraphics[width=0.33\linewidth]{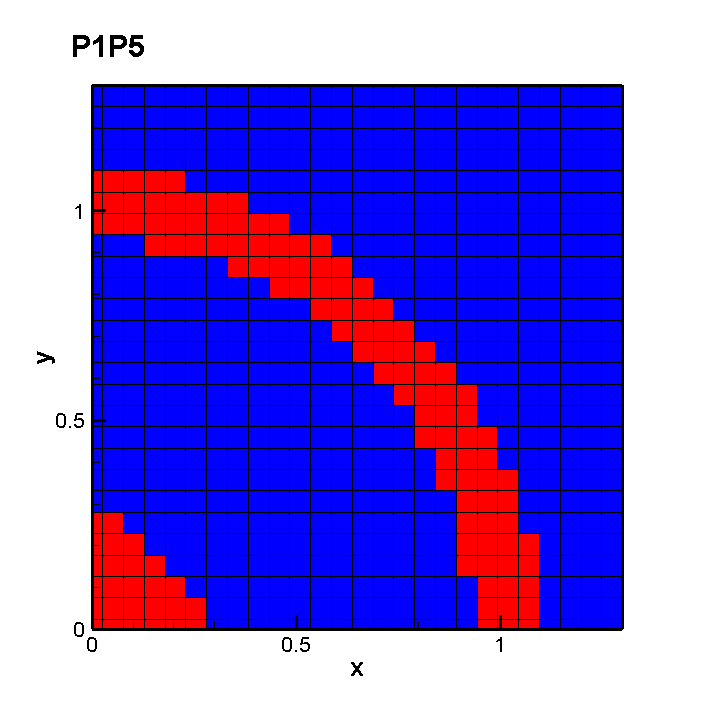}%
	\includegraphics[width=0.33\linewidth]{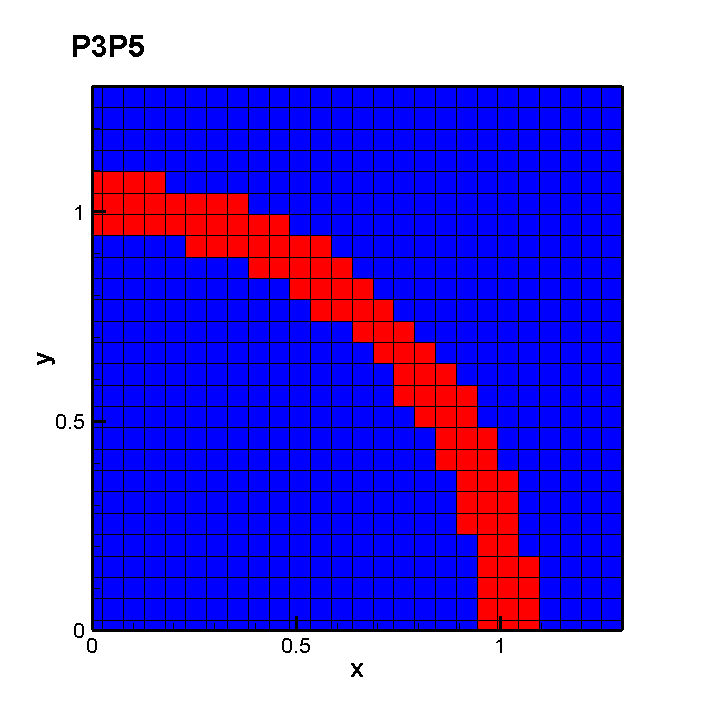}%
	\includegraphics[width=0.33\linewidth]{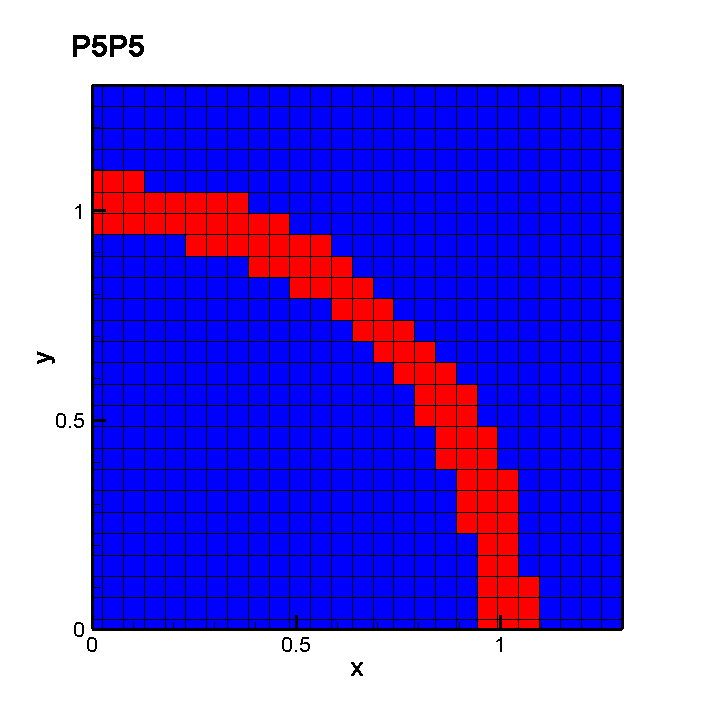}\\
	\includegraphics[width=0.33\linewidth]{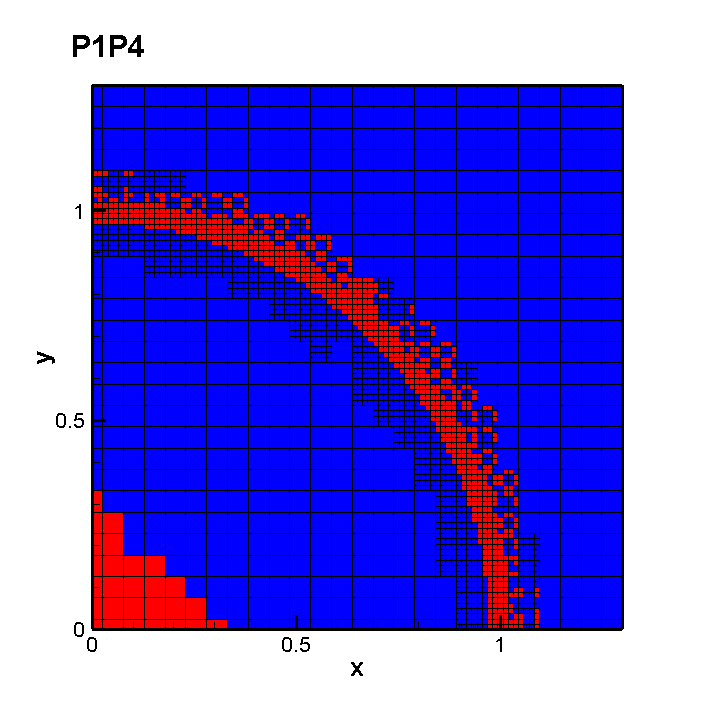}%
	\includegraphics[width=0.33\linewidth]{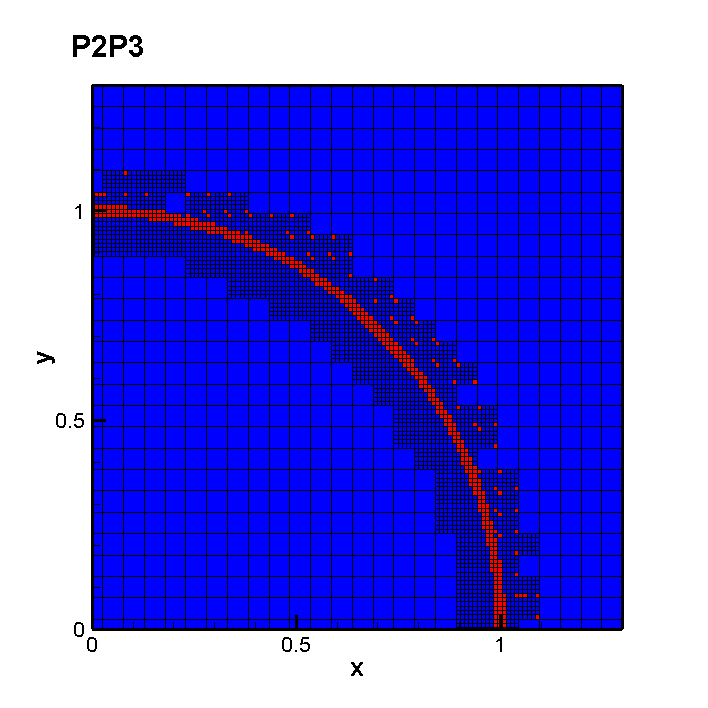}%
	\includegraphics[width=0.33\linewidth]{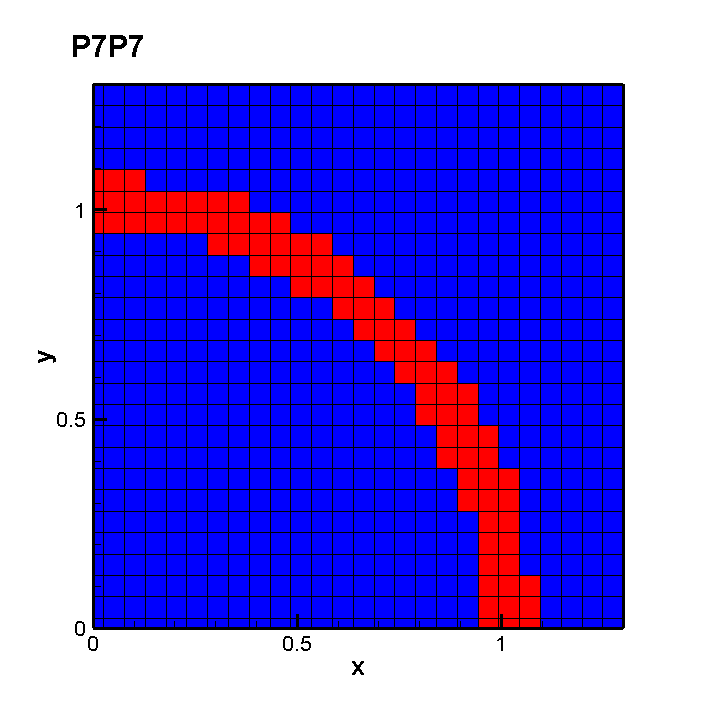}%
	\caption{Limited cells (red) and unlimited ones (blue) for the Sedov problem solved with a selection of high order methods, 
		i.e. $\PN{1}{5}, \PN{3}{5}, \PN{5}{5}, \PN{1}{4}, \PN{2}{3}, \PN{7}{7}$. A part some spurious oscillations of the $\PN{N}{M}$ schemes with $N=1$, 
		the limiter activates exactly at the shock location.}
	\label{fig:sedovLimiter}
\end{figure}
\begin{table}
	\caption{Sedov problem. We report the total CPU time in seconds and the values of the density peak (that should be equal to $6$) obtained with a selection of methods going from order $4$ to $8$. 
		The results are ordered with respect to the density peak value. One can notice that the hybrid schemes (as the $\PN{2}{3}$ and the $\PN{3}{5}$) have accurate results at a lower computational cost.} 
	\label{tab.Sedov}
	\begin{center} 	
		\begin{tabular}{lccc} 
			\hline 
			Method & Order & Total CPU time & Density peak  \\ 
			\hline
			P1P4 & 5 & 24 & 4.01 \\
			P1P5 & 6 & 37 & 4.02 \\
			P2P3 & 4 & 76 & 4.48 \\ 	
			P3P5 & 6& 361 & 4.94 \\			
			P5P5 & 6& 1458 & 5.08 \\
			P1P4 + AMR & 5 & 847 & 5.20\\
			P7P7 &8&  4428 & 5.25 \\	 
			P2P3 + AMR &4&  1168& 5.53\\
			P3P3 + AMR & 4 & 2454  & 5.61 \\
			\hline  
		\end{tabular}		
	\end{center}
\end{table}
being $E_\text{tot}=0.244816$  the the total energy concentrated at $\mathbf{x}=\mathbf{0}$ and $|V_\text{or}|$ the total volume of $V_{or}$.

We solve this numerical test with a selection of high order methods going from fourth to eighth order of accuracy on a mesh of $50\times50$ elements with and without AMR. 
When the adaptive mesh refinement is activated, we take $\lmax=2$ and $\rRef=3$.
For all these test cases we employ the Rusanov flux, $CFL=0.9$, and the WENO \textit{a posteriori} subcell finite volume limiter.  

We show the results on the activation of the limiter in Figure~\ref{fig:sedovLimiter}, and the obtained density profiles in Figure~\ref{fig:sedovRho}. 
Finally, we compare the performances of a selection of methods in Table~\ref{tab.Sedov}, 
in particular we highlight the needed computational times versus their capability of capturing the high density peak. 
We can remark that the hybrid schemes (as the $\PN{2}{3}$ and the $\PN{3}{5}$) have accurate and robust results at a lower computational cost; 
also the combination with adaptive mesh refinement helps in increasing the accuracy keeping the computational cost low.

\subsubsection{Double Mach Reflection}
\begin{figure}[h]
	\centering
	\includegraphics[width=1.0\linewidth]{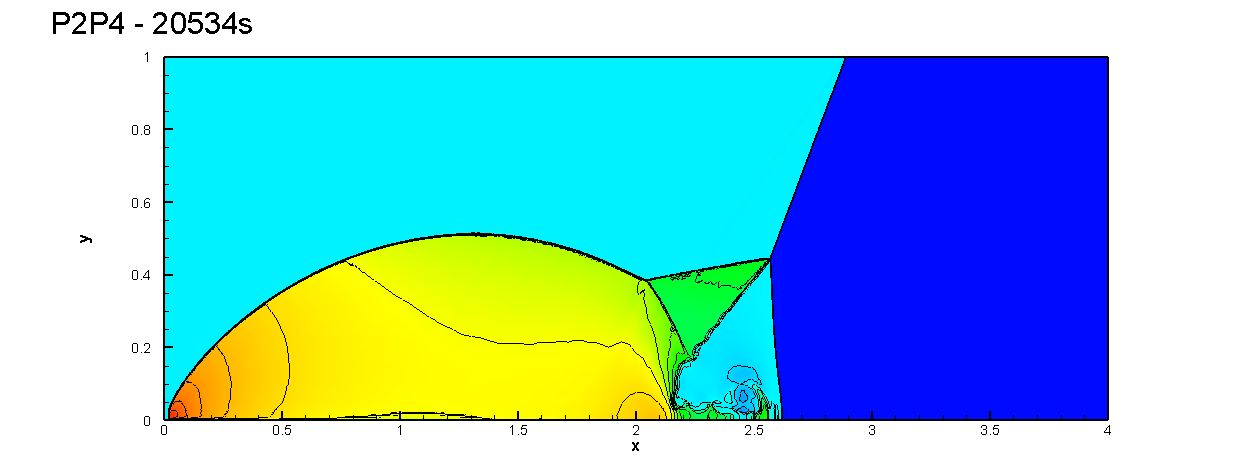}\\
	\includegraphics[width=1.0\linewidth]{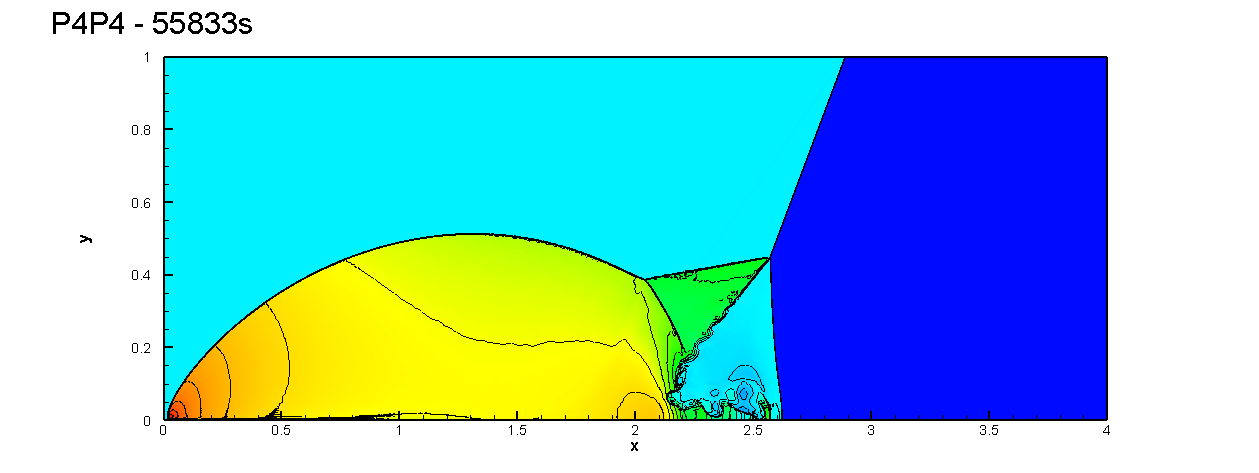} 
	\caption{Double Mach Reflection density contours obtained with two fifth order schemes, namely the $\PN{2}{4}$ and the $\PN{4}{4}$ schemes. 
		We plot 30 equally spaced contour lines from $1.5$ to $22.9705$ as suggested in~\cite{shi2003resolution}.}
	\label{fig.DMRorder5}
\end{figure}
\begin{figure}[h]
	\centering	
	\includegraphics[width=0.45\linewidth]{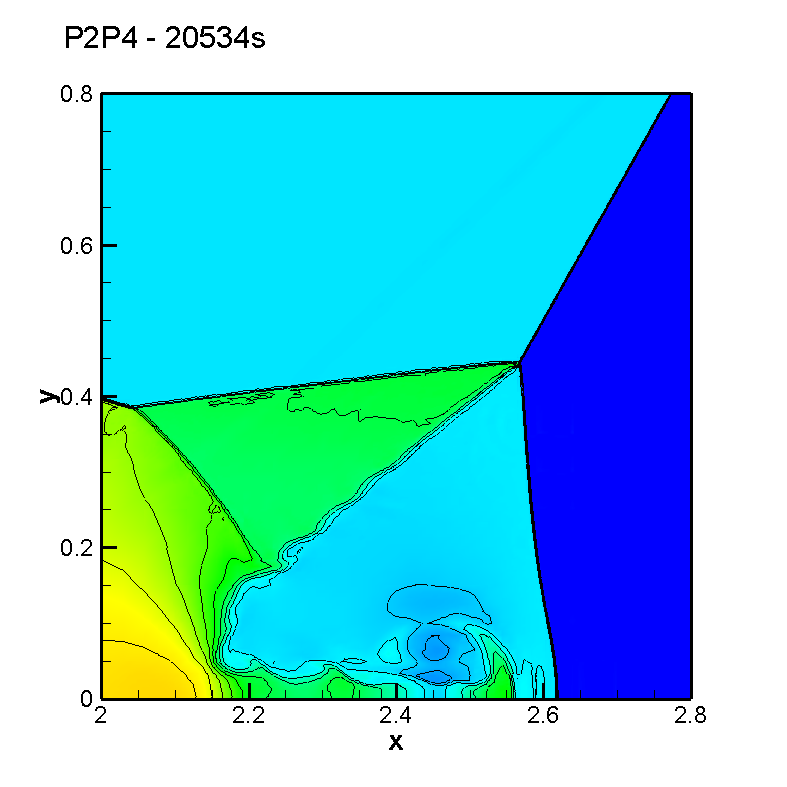}
	\includegraphics[width=0.45\linewidth]{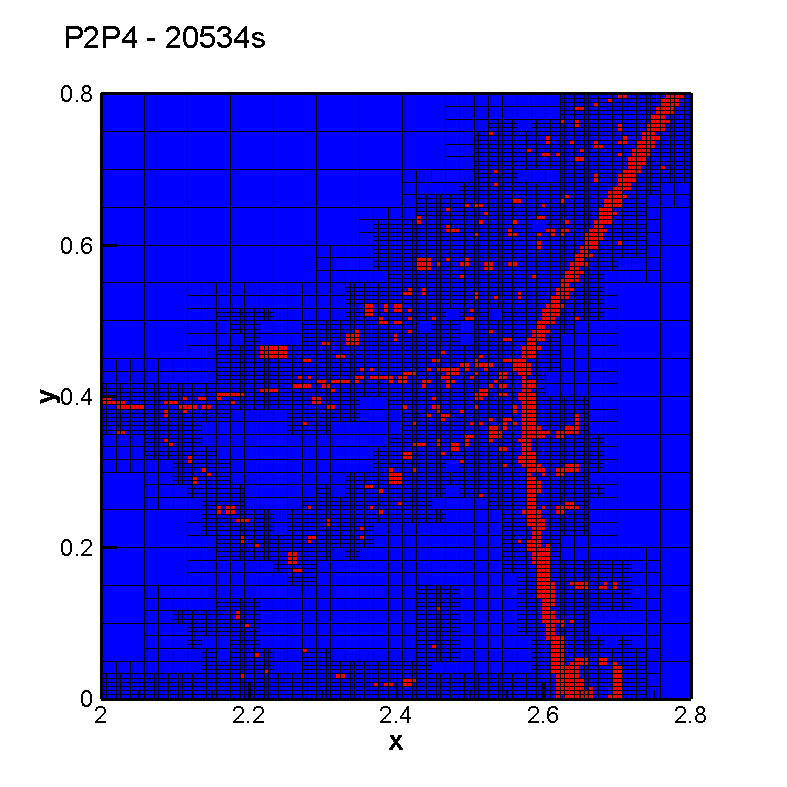}\\
	\includegraphics[width=0.45\linewidth]{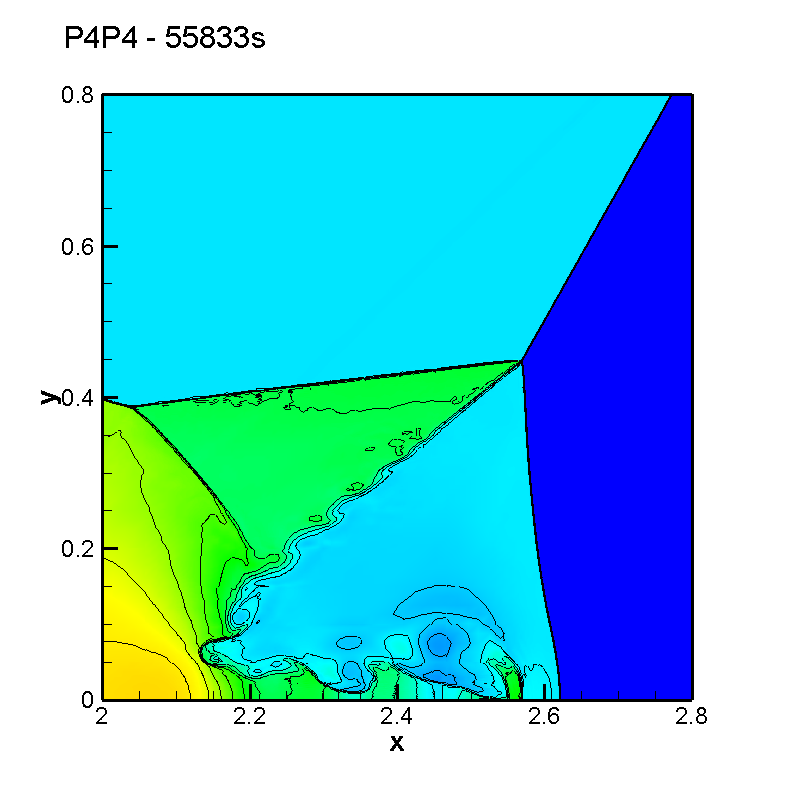}
	\includegraphics[width=0.45\linewidth]{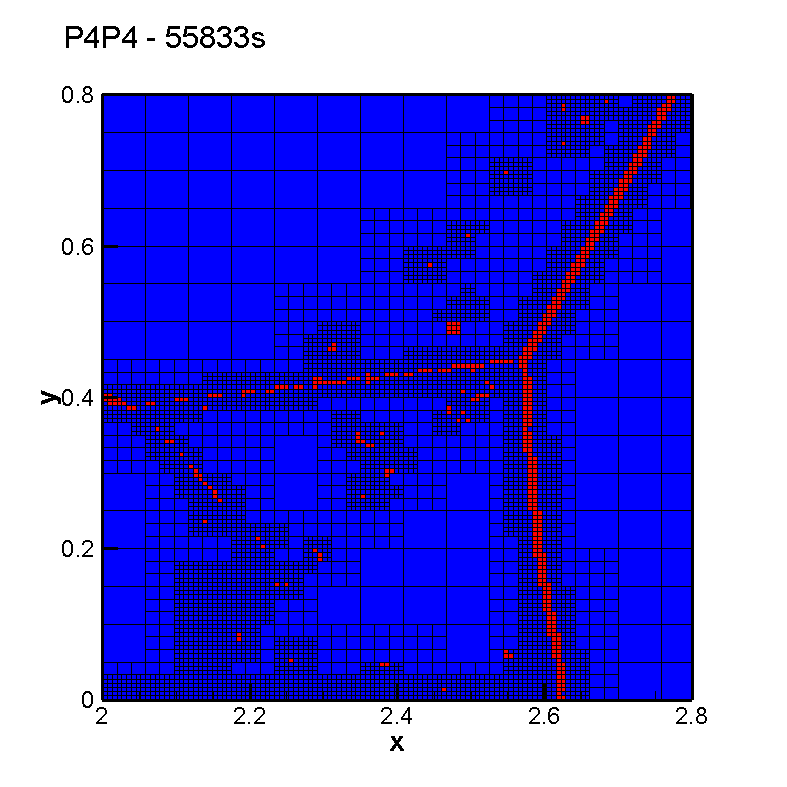}
	\caption{Double Mach Reflection. Left: density contours with 30 equally spaced contour lines from $1.5$ to $22.9705$. Right: limited cells (red) and unlimited cells (blue).}
	\label{fig.DMRorder5_zoom}
\end{figure}

The double Mach reflection problem was first studied by Woodward and Colella in~\cite{woodwardcol84} from which we take the setup also for our test.   

We consider a computational domain $\Omega = [0,4] \times [0,1]$ covered with a coarse mesh of  $72\times24$ elements where we activate the adaptive mesh refinement with $\lmax=2$ and $\rRef=3$, 
and we compare the behavior of two fifth order schemes, namely the hybrid $\PN{2}{4}$ scheme and the $\PN{4}{4}$ pure DG scheme.
For all the simulations, we employ the Rusanov flux 
and the second order TVD \textit{a posteriori} subcell finite volume limiter. 

At the beginning of the computation a shock wave moving at Mach number $10$ is positioned at $(x,y) = (1/6,0)$ with an angle of $60^\circ$ with respect to the $x$-axis 
and the initial pre-shock conditions (on the left of the shock) are given by a constant density equal to $1.4$ and a constant pressure $p=1$.
At the bottom boundary we employ reflective boundary conditions for $x>1/6$ where we suppose the presence of a wall, and the exact post-shock conditions for $0 \le x \le 1/6$ to mimic an angled wedge.
At the top boundary, the flow variables are set to describe the exact motion of the Mach $10$ shock. 
Finally at the left and right boundaries we set inflow and outflow boundaries.

The obtained numerical results are shown in Figure~\ref{fig.DMRorder5} for the entire domain; 
we also plot a zoom in Figure~\ref{fig.DMRorder5_zoom} where, one can notice the roll up of the Mach stem due to Kelvin-Helmholtz instabilities.

\subsection{Ideal MHD equations}
\label{sec.MHD}

Next, we consider the equations of ideal classical magnetohydrodynamics (MHD) which, with respect to the previous set of equations, 
take also into account the evolution of the magnetic field $\B$.
The vector of the conserved variables $\Q$ and the flux tensor $\F$ of the general form~\eqref{eq.generalform} are given by 
\begin{equation}
\label{MHDTerms}
\Q=\left( \begin{array}{c} \rho \\ \rho \mathbf{v} \\ \rho E \\ \B \\ \psi \end{array} \right)\!, \quad \F(\Q) = 
\left( \begin{array}{c}  \rho \v  \\ 
\rho \v \otimes \v + p_t \mathbf{I} - \frac{1}{4 \pi} \B \otimes \B \\ 
\v (\rho E + p_t ) - \frac{1}{4 \pi} \B ( \v \cdot \B ) \\ 
\v \otimes \B - \B \otimes \v + \psi \mathbf{I} \\
c_h^2 \B  \end{array} \right)\!.
\end{equation}
Here, $\B=(B_x,B_y,B_z)$ represents the magnetic field and $p_t=p+\frac{1}{8\pi}\mathbf{B}^2$ is the total pressure. 
The hydrodynamic pressure is given by the equation of state used to close the system, thus
\begin{equation}
p = \left(\gamma - 1 \right) \left(\rho E - \frac{1}{2}\mathbf{v}^2 - \frac{\mathbf{B}^2}{8\pi}\right).
\label{MHDeos}
\end{equation}
System~\eqref{MHDTerms} requires an additional constraint on the divergence of the magnetic field to be satisfied, that is
\begin{equation}
\nabla \cdot \mathbf{B} = 0.
\label{eq:divB}
\end{equation}  
Here,~\eqref{MHDTerms} includes one additional scalar PDE for the evolution of the variable $\psi$, 
which is needed to transport divergence errors outside the computational domain with an artificial \textit{divergence cleaning} speed $c_h$, see~\cite{MunzCleaning,Dedneretal}. 
A similar approach is adopted in~\cite{ADERDGVisc,boscheri2014high,boscheri2017efficient}. 
A more recent and more sophisticated methodology to fulfill this condition exactly at the discrete level also in the context of high order ADER WENO finite volume schemes on unstructured simplex meshes can be found in~\cite{MHDdivFree2015}.

\subsubsection{MHD vortex}

\begin{landscape} 
\begin{table*}
	\caption{Numerical convergence table for general $P_NP_M$ schemes for the MHD vortex. The error norms refer to the variable $\rho$ at time $t_f=1.0$ in the $L_2$ norm.} 
	\label{tab.orderOfconvergenceMHD}
	\begin{center} 	
		\begin{tabular}{|c|ccc|ccc|ccc|ccc|ccc|ccc|} 
			\hline 
			h                        & \CPU      & $L_2$             & \OLdue        &   \CPU      & $L_2$             & \OLdue       &     \CPU      & $L_2$             & \OLdue       &      \CPU      & $L_2$             & \OLdue       &    \CPU      & $L_2$             & \OLdue       &        \CPU     & $L_2$            & \OLdue          \\[1pt]
			\hline                                                                                                                                                                                                                                                                                                                                   
			$\mathbf{\mathcal{O}2}$      &               \mcolt{0}{1}                    &                \mcolt{1}{1}                    &               &                  &               &                &                  &               &              &                   &              &                 &                  &                   \\
			\smm 2.0e-02       \smmm &\sm 37 \sm & \sm 4.7e-05   \sm & \smm     \smm & \sm 69  \sm & \sm 2.6e-06   \sm & \smm     \smm &   \sm     \sm & \sm          \sm & \smm     \smm &    \sm     \sm & \sm          \sm & \smm     \smm &  \sm     \sm & \sm          \sm & \smm     \smm &      \sm    \sm & \sm          \sm & \smm     \smm    \\
			\smm 1.6e-02       \smmm &\sm 54 \sm & \sm 3.7e-05   \sm & \smm 1.4 \smm & \sm 132 \sm & \sm 1.8e-06   \sm & \smm 1.9 \smm &   \sm     \sm & \sm          \sm & \smm     \smm &    \sm     \sm & \sm          \sm & \smm     \smm &  \sm     \sm & \sm          \sm & \smm     \smm &      \sm    \sm & \sm          \sm & \smm     \smm    \\
			\smm 1.2e-02       \smmm &\sm 153\sm & \sm 2.4e-05   \sm & \smm 1.4 \smm & \sm 392 \sm & \sm 1.0e-06   \sm & \smm 1.9 \smm &   \sm     \sm & \sm          \sm & \smm     \smm &    \sm     \sm & \sm          \sm & \smm     \smm &  \sm     \sm & \sm          \sm & \smm     \smm &      \sm    \sm & \sm          \sm & \smm     \smm    \\
			\smm 1.0e-02       \smmm &\sm 249\sm & \sm 1.7e-05   \sm & \smm 1.6 \smm & \sm 574 \sm & \sm 6.7e-07   \sm & \smm 2.1 \smm &   \sm     \sm & \sm          \sm & \smm     \smm &    \sm     \sm & \sm          \sm & \smm     \smm &  \sm     \sm & \sm          \sm & \smm     \smm &      \sm    \sm & \sm          \sm & \smm     \smm    \\[1pt]
			\hline                                                                                                                                                                                                                                                                                                                               
			$\mathbf{\mathcal{O}3}$      &               \mcolt{0}{2}                    &                \mcolt{1}{2}                    &                  \mcolt{2}{2}                    &               &                   &               &             &                   &               &                 &                   &                    \\
			\smm 1.2e-01       \smmm &\sm 0.6\sm & \sm 1.5e-06   \sm & \smm     \smm & \sm 1.4\sm & \sm 9.4e-07   \sm & \smm     \smm &   \sm 2.7\sm & \sm 4.2e-06   \sm & \smm     \smm &    \sm    \sm & \sm           \sm & \smm     \smm &  \sm    \sm & \sm           \sm & \smm     \smm &      \sm    \sm & \sm           \sm & \smm     \smm    \\
			\smm 8.3e-02       \smmm &\sm 1.8\sm & \sm 4.1e-07   \sm & \smm 3.3 \smm & \sm 4.3\sm & \sm 2.5e-07   \sm & \smm 3.2 \smm &   \sm 8.7\sm & \sm 1.5e-06   \sm & \smm 2.4 \smm &    \sm    \sm & \sm           \sm & \smm     \smm &  \sm    \sm & \sm           \sm & \smm     \smm &      \sm    \sm & \sm           \sm & \smm     \smm    \\
			\smm 6.2e-02       \smmm &\sm 3.7\sm & \sm 1.6e-07   \sm & \smm 3.1 \smm & \sm 10 \sm & \sm 1.1e-07   \sm & \smm 2.8 \smm &   \sm 20 \sm & \sm 7.9e-07   \sm & \smm 2.4 \smm &    \sm    \sm & \sm           \sm & \smm     \smm &  \sm    \sm & \sm           \sm & \smm     \smm &      \sm    \sm & \sm           \sm & \smm     \smm    \\
			\smm 5.0e-02       \smmm &\sm 6.9\sm & \sm 8.6e-08   \sm & \smm 3.0 \smm & \sm 19 \sm & \sm 5.4e-08   \sm & \smm 3.3 \smm &   \sm 31 \sm & \sm 4.6e-07   \sm & \smm 2.4 \smm &    \sm    \sm & \sm           \sm & \smm     \smm &  \sm    \sm & \sm           \sm & \smm     \smm &      \sm    \sm & \sm           \sm & \smm     \smm    \\[1pt]
			\hline                                                                                                                                                                                                                                                                                                                               
			$\mathbf{\mathcal{O}4}$      &               \mcolt{0}{3}                    &                \mcolt{1}{3}                    &                  \mcolt{2}{3}                    &                   \mcolt{3}{3}                    &             &                   &               &                 &                   &                    \\
			\smm 2.5e-01       \smmm &\sm 0.3\sm & \sm 1.1e-05   \sm & \smm     \smm & \sm 0.5\sm & \sm 1.0e-05   \sm & \smm     \smm &   \sm 0.9\sm & \sm 6.3e-07   \sm & \smm     \smm &    \sm 2.2\sm & \sm 4.0e-07   \sm & \smm     \smm &  \sm    \sm & \sm           \sm & \smm     \smm &      \sm    \sm & \sm           \sm & \smm     \smm    \\
			\smm 1.6e-01       \smmm &\sm 1.5\sm & \sm 1.2e-06   \sm & \smm 5.5 \smm & \sm 1.3\sm & \sm 1.9e-06   \sm & \smm 4.1 \smm &   \sm 2.6\sm & \sm 1.2e-07   \sm & \smm 4.0 \smm &    \sm 3.7\sm & \sm 8.9e-08   \sm & \smm 3.7 \smm &  \sm    \sm & \sm           \sm & \smm     \smm &      \sm    \sm & \sm           \sm & \smm     \smm    \\
			\smm 1.2e-01       \smmm &\sm 4.4\sm & \sm 2.9e-07   \sm & \smm 5.0 \smm & \sm 2.8\sm & \sm 6.1e-07   \sm & \smm 3.9 \smm &   \sm 5.8\sm & \sm 3.4e-08   \sm & \smm 4.5 \smm &    \sm 11 \sm & \sm 2.6e-08   \sm & \smm 4.2 \smm &  \sm    \sm & \sm           \sm & \smm     \smm &      \sm    \sm & \sm           \sm & \smm     \smm    \\
			\smm 1.0e-01       \smmm &\sm 4.7\sm & \sm 1.5e-07   \sm & \smm 2.9 \smm & \sm 5.2\sm & \sm 2.5e-07   \sm & \smm 3.9 \smm &   \sm 11 \sm & \sm 1.3e-08   \sm & \smm 4.2 \smm &    \sm 15 \sm & \sm 1.0e-08   \sm & \smm 4.0 \smm &  \sm    \sm & \sm           \sm & \smm     \smm &      \sm    \sm & \sm           \sm & \smm     \smm    \\[1pt]
			\hline                                                                                                                                                                                                                                                                                                                               
			$\mathbf{\mathcal{O}5}$      &               \mcolt{0}{4}                    &                \mcolt{1}{4}                    &                  \mcolt{2}{4}                    &                   \mcolt{3}{4}                    &                 \mcolt{4}{4}                    &                 &                  &                     \\
			\smm 4.0e-01       \smmm &\sm 0.2\sm & \sm 2.4e-05   \sm & \smm     \smm & \sm 0.3\sm & \sm 1.1e-06   \sm & \smm     \smm &   \sm 0.5\sm & \sm 4.6e-06   \sm & \smm     \smm &    \sm 0.8 \sm & \sm 1.6e-07   \sm & \smm    \smm &  \sm 0.9\sm & \sm 1.9e-07   \sm & \smm    \smm &      \sm    \sm & \sm          \sm & \smm    \smm      \\
			\smm 3.3e-01       \smmm &\sm 0.3\sm & \sm 8.5e-06   \sm & \smm 5.8 \smm & \sm 0.4\sm & \sm 5.5e-07   \sm & \smm 3.8 \smm &   \sm 0.7\sm & \sm 2.0e-06   \sm & \smm 4.6 \smm &    \sm 1.2 \sm & \sm 5.1e-08   \sm & \smm 6.4\smm &  \sm 1.5\sm & \sm 8.7e-08   \sm & \smm 4.4\smm &      \sm    \sm & \sm          \sm & \smm    \smm      \\
			\smm 2.8e-01       \smmm &\sm 0.4\sm & \sm 3.7e-06   \sm & \smm 5.4 \smm & \sm 0.7\sm & \sm 2.8e-07   \sm & \smm 4.4 \smm &   \sm 1.1\sm & \sm 9.5e-07   \sm & \smm 4.9 \smm &    \sm 1.9 \sm & \sm 1.8e-08   \sm & \smm 6.5\smm &  \sm 2.3\sm & \sm 4.3e-08   \sm & \smm 4.4\smm &      \sm    \sm & \sm          \sm & \smm    \smm      \\
			\smm 2.5e-01       \smmm &\sm 0.6\sm & \sm 1.9e-06   \sm & \smm 4.7 \smm & \sm 1.0\sm & \sm 1.4e-07   \sm & \smm 5.2 \smm &   \sm 1.6\sm & \sm 4.9e-07   \sm & \smm 4.9 \smm &    \sm 2.9 \sm & \sm 8.6e-09   \sm & \smm 5.8\smm &  \sm 3.4\sm & \sm 2.4e-08   \sm & \smm 4.4\smm &      \sm    \sm & \sm          \sm & \smm    \smm      \\[1pt]
			\hline                                                                                                                                                                                                                                                                                                                               
			$\mathbf{\mathcal{O}6}$      &               \mcolt{0}{5}                    &                \mcolt{1}{5}                    &                  \mcolt{2}{5}                    &                   \mcolt{3}{5}                    &                 \mcolt{4}{5}                    &                     \mcolt{5}{5}                       \\
			\smm 7.1e-01       \smmm &\sm 0.1\sm & \sm 5.3e-04   \sm & \smm     \smm & \sm 0.1\sm & \sm 4.2e-05   \sm & \smm     \smm &   \sm 0.2\sm & \sm 3.6e-06   \sm & \smm     \smm &    \sm 0.3\sm & \sm 4.9e-06   \sm & \smm     \smm &  \sm 0.4\sm & \sm 3.6e-07   \sm & \smm     \smm &      \sm 0.6\sm & \sm 2.1e-07   \sm & \smm     \smm    \\			
			\smm 5.5e-01       \smmm &\sm 0.3\sm & \sm 1.3e-04   \sm & \smm 5.5 \smm & \sm 0.2\sm & \sm 1.2e-05   \sm & \smm 4.9 \smm &   \sm 0.3\sm & \sm 1.1e-06   \sm & \smm 4.7 \smm &    \sm 0.5\sm & \sm 1.2e-06   \sm & \smm 5.5 \smm &  \sm 0.8\sm & \sm 7.6e-08   \sm & \smm 6.2 \smm &      \sm 1.2\sm & \sm 5.8e-08   \sm & \smm 5.1 \smm    \\
			\smm 4.5e-01       \smmm &\sm 0.3\sm & \sm 4.1e-05   \sm & \smm 5.7 \smm & \sm 0.4\sm & \sm 3.9e-06   \sm & \smm 5.6 \smm &   \sm 0.5\sm & \sm 4.0e-07   \sm & \smm 5.0 \smm &    \sm 0.9\sm & \sm 3.8e-07   \sm & \smm 5.9 \smm &  \sm 1.4\sm & \sm 2.6e-08   \sm & \smm 5.3 \smm &      \sm 2.0\sm & \sm 1.7e-08   \sm & \smm 6.0 \smm    \\
			\smm 3.8e-01       \smmm &\sm 0.5\sm & \sm 1.5e-05   \sm & \smm 6.0 \smm & \sm 0.5\sm & \sm 1.4e-06   \sm & \smm 6.0 \smm &   \sm 0.9\sm & \sm 1.6e-07   \sm & \smm 5.5 \smm &    \sm 2  \sm & \sm 1.3e-07   \sm & \smm 6.0 \smm &  \sm 2.1\sm & \sm 1.0e-08   \sm & \smm 5.2 \smm &      \sm 2.7\sm & \sm 6.7e-09   \sm & \smm 5.7 \smm    \\
			\hline 
		\end{tabular}		
	\end{center}
\end{table*}
\end{landscape} 

First, for the numerical convergence studies, we solve the vortex test problem proposed by Balsara in~\cite{Balsara2004}.  
The computational domain is given by the box $\Omega=[0,10]\times[0,10]$ with wall boundary conditions imposed everywhere. 
The initial condition can be written in terms of the vector of primitive variables $\mathbf{V} = ( \rho, u, v, w, p, B_x, B_y,  B_z, \Psi )^T$ as 
\begin{equation}
\mathbf{V}(\mathbf{x},0) =
( 1, \delta u, \delta v, 0, 1+\delta p, \delta B_x, \delta B_y,  0, 0 )^T, 
\end{equation}
with $\delta \mathbf{v} = (\delta u, \delta v, 0)^T$, $\ \delta \mathbf{B} = ( \delta B_x, \delta B_y, 0 )^T$ and 
\begin{equation}
\begin{aligned} 
\label{eqn.mhd3d.ic1}
&\delta \mathbf{v} = \frac{\kappa}{2\pi} e^{ q(1-r^2)} \mathbf{e}_z \times \mathbf{r}  \\ 
&\delta \mathbf{B} = \frac{\mu}{2\pi}    e^{ q(1-r^2)} \mathbf{e}_z \times \mathbf{r},  \\ 
&\delta p = \frac{1}{64 q \pi^3} \left( \mu^2 (1 - 2 q r^2) - 4 \kappa^2 \pi \right) e^{2q(1-r^2)},
\end{aligned} 
\end{equation}
where $\mathbf{e}_z = (0,0,1)$, $\mathbf{r} = (x-5,y-5,0)$ and $r = \left\| \mathbf{r} \right\| = \sqrt{ (x-5)^2 + (y-5)^2 }$. The divergence cleaning speed 
is chosen as $c_h=2$. The other parameters are $q=\halb$, $\kappa=1$ and $\mu=\sqrt{4 \pi}$, according to~\cite{Balsara2004}. 

In Table~\ref{tab.orderOfconvergenceMHD}, 
we report the convergence rates from second up to sixth order of accuracy for the MHD vortex test problem run on a sequence of successively refined meshes up to the final time $t_f=1.0$. 
The optimal order of accuracy is achieved both in space and time both for the hybrid schemes $\PN{N}{M}$ with $M>N$ and for the pure DG schemes with $N=M$.

\subsubsection{MHD rotor problem}
\begin{figure}[b]
	\centering
	\includegraphics[width=0.33\linewidth]{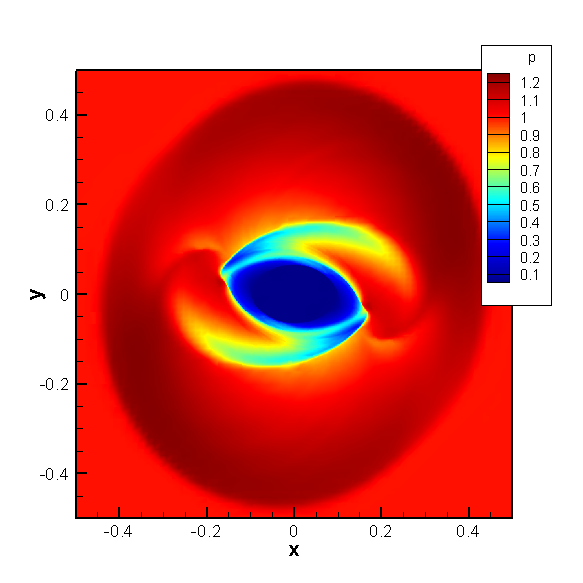}%
	\includegraphics[width=0.33\linewidth]{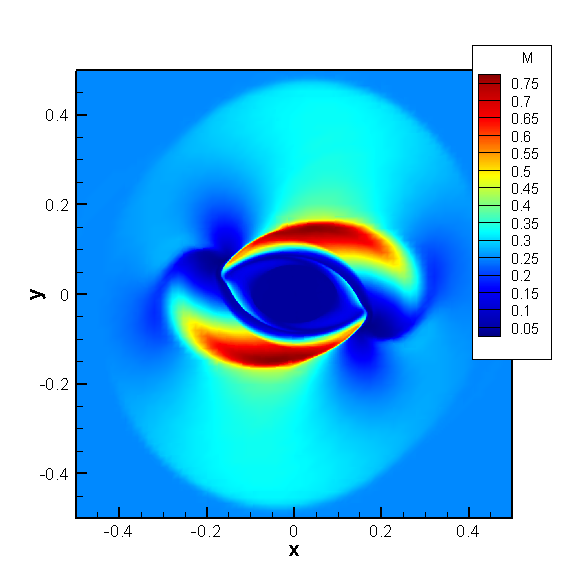}%
	\includegraphics[width=0.33\linewidth]{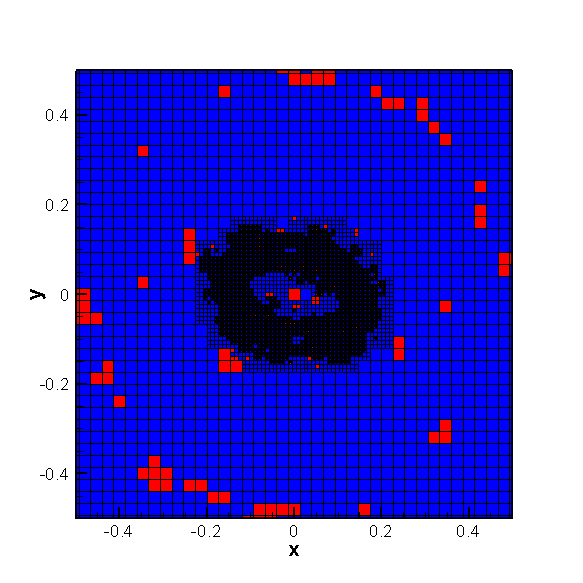}\\
	\includegraphics[width=0.33\linewidth]{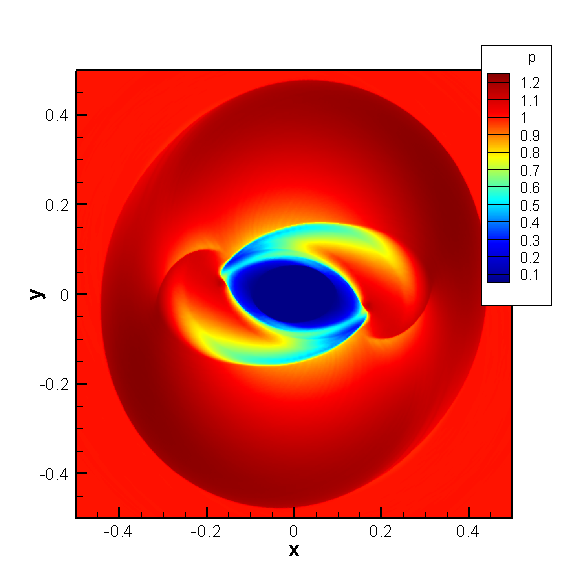}%
	\includegraphics[width=0.33\linewidth]{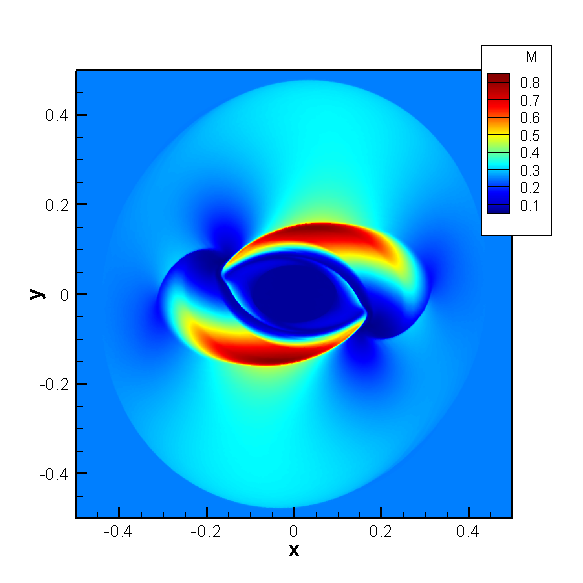}%
	\includegraphics[width=0.33\linewidth]{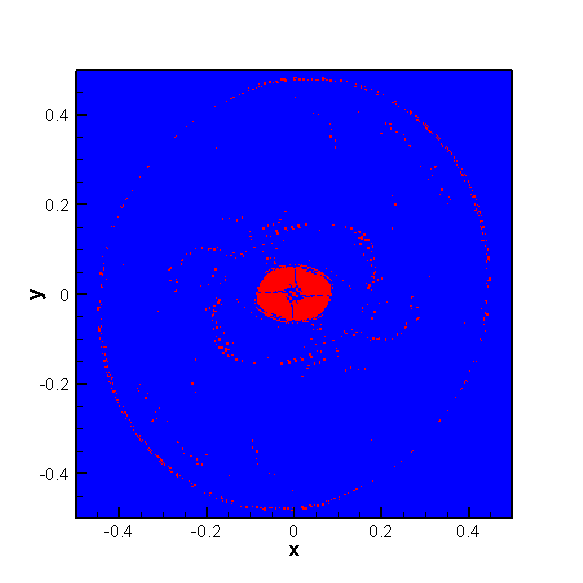}%
	\caption{MHD Rotor problem at the final time $t_f =0.25$ solved with our $P2P4$ fifth order scheme. 
		In the left column we plot the pressure contours, in the central column the magnetic density profile $M = \frac{(B_x^2 + B_y^2+B_z^2)}{(8\pi)}$ 
		and in the right column we depict in red the troubled cells and in blue unlimited cells. 
		The results on the first row are obtained with a coarse mesh of $45\times 45$ cells and $\ell_{\max}=2$ levels of refinement with $\rRef=3$. 
		The results on the second row are obtained by using a fine uniform grid of $405\times405$ elements corresponding to the finest AMR grid level. 
		The computation on the finer grids takes twice the time of the computation on a coarse mesh with AMR (i.e $7890$s instead of $3779$s).	}
	\label{fig:MHDRotor_P2P4_45x45_ref2_p}
\end{figure}

\begin{figure}[b]
	\centering
	\includegraphics[width=0.33\linewidth]{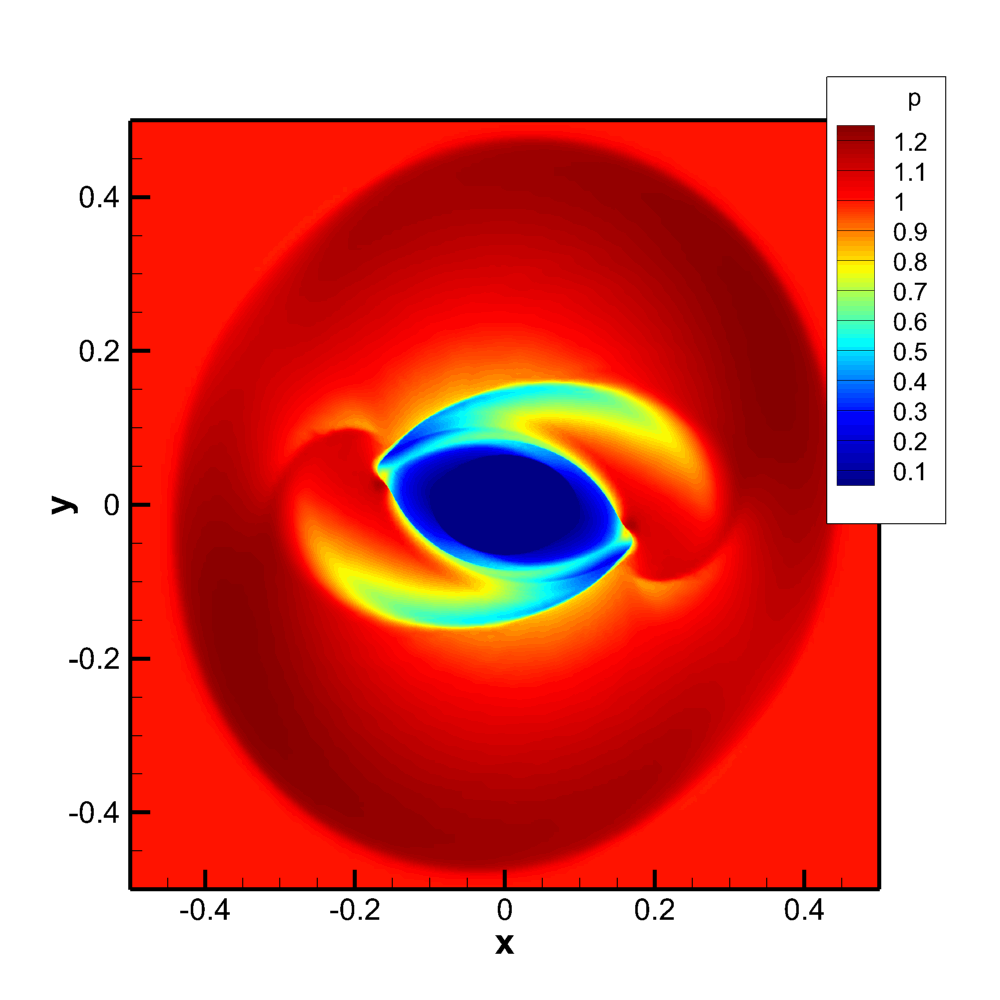}%
	\includegraphics[width=0.33\linewidth]{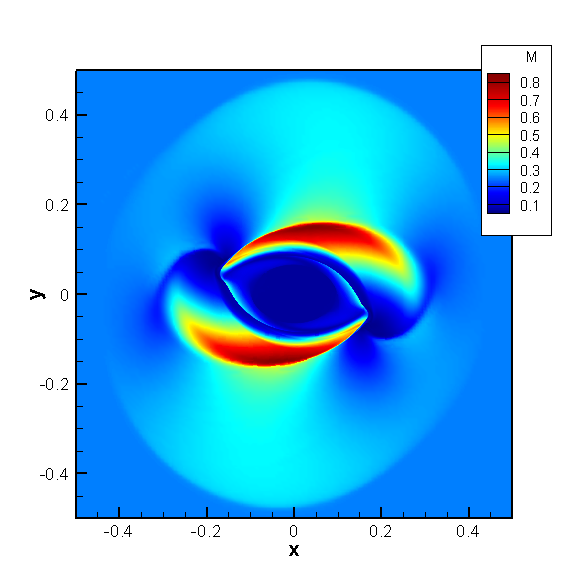}%
	\includegraphics[width=0.33\linewidth]{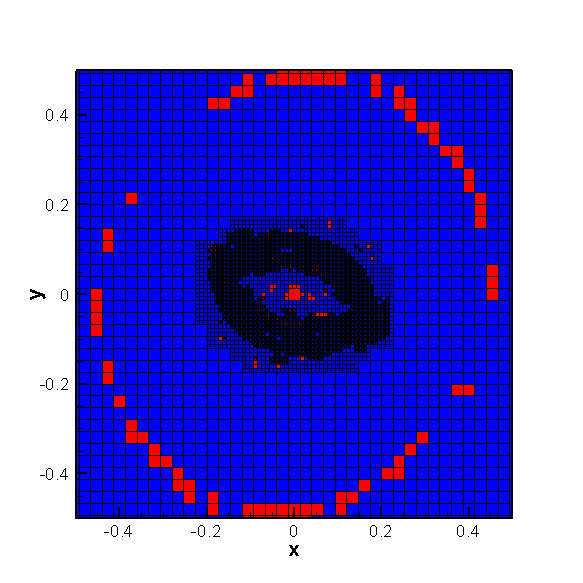}\\
	\includegraphics[width=0.33\linewidth]{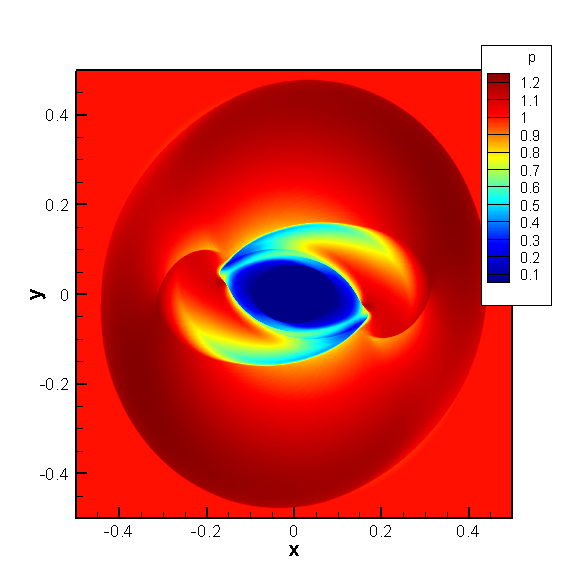}%
	\includegraphics[width=0.33\linewidth]{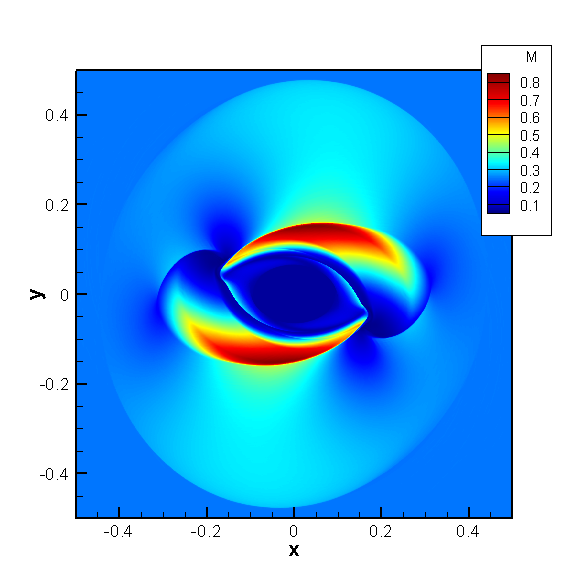}%
	\includegraphics[width=0.33\linewidth]{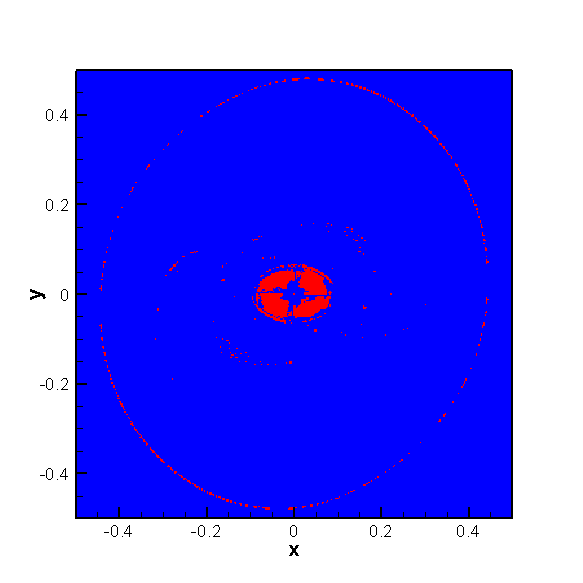}%
	\caption{MHD Rotor problem at the final time $t_f =0.25$ solved with our $P4P5$ sixth order scheme. 
		In the left column we plot the pressure contours, in the central column the magnetic density profile $M = \frac{(B_x^2 + B_y^2+B_z^2)}{(8\pi)}$ 
		and in the right column we depict in red the troubled cells and in blue unlimited cells in blue. 
		The results on the first row are obtained with a coarse mesh of $45\times 45$ cells and $\ell_{\max}=2$ levels of refinement with $\rRef=3$. 
		The results on the second row are obtained by using a fine uniform grid of $405\times405$ elements corresponding to the finest AMR grid level. 
		The computation on the finer grids takes twice the time of the computation on a coarse mesh with AMR.	}
	\label{fig:MHDRotor_P4P5_45x45_ref2_p}  
\end{figure}

The MHD rotor problem is a classical benchmark for MHD that was first proposed by Balsara and Spicer in~\cite{BalsaraSpicer1999}. 
It consists of a rapidly rotating fluid of high density embedded in a fluid at rest with low density. 
Both fluids are subject to an initially constant magnetic field. 

The rotor produces torsional Alfv\'en waves that are launched into the outer fluid at rest, resulting in a decrease
of angular momentum of the spinning rotor. 
We consider as computational domain the square $\Omega=[-0.5, 0.5]\times[-0.5, 0.5]$
and as initial condition we take the density inside a circle of radius $r \leq 0.1$ equal to $\rho=10$, while the density of the ambient fluid at rest is
set to $\rho=1$. 
The rotor has an angular velocity of $\omega=10$. 
The pressure is $p=1$ and the magnetic field vector is set to $\B = (2.5, 0, 0)^T$ in the entire domain.
As proposed by Balsara and Spicer we apply a linear taper to the velocity and to the density in the range from 
$0.1 \leq r \leq 0.105$ so that density and velocity match those of the ambient fluid at rest at a radius of $r=0.105$.  
The speed for the hyperbolic divergence cleaning is set to $c_h=8$ and $\gamma=1.4$ is used. Wall boundary conditions are applied everywhere. 

We run this problem on a coarse mesh made of $45\times45$ elements activating the AMR procedure with $\ell_{\max}=2$ levels of refinement and $\rRef=3$,
and for comparison we also employ a finer uniform mesh of $405\times405$ elements corresponding to the finest AMR grid level.
In particular, we have employed the $\PN{2}{4}$ fifth order scheme and the $\PN{4}{5}$ sixth order scheme with the Rusanov numerical flux and our \textit{a posteriori} subcell WENO FV limiter. 
In all the cases, we can observe a good agreement between the obtained numerical results and those available in the literature, see Figures~\ref{fig:MHDRotor_P2P4_45x45_ref2_p}-\ref{fig:MHDRotor_P4P5_45x45_ref2_p}. 

%
%
%
%
%
%

\subsection{MHD Orszag-Tang vortex}

We consider now the the vortex system of Orszag and Tang~\cite{OrszagTang,DahlburgPicone,PiconeDahlburg} for the ideal MHD equations. 
We choose as computational domain the square $\Omega=[0,2\pi]\times[0,2\pi]$ with periodic boundary conditions set everywhere; 
we cover it with a uniform grid of $128\times 128$ elements.

The initial condition written in terms of primitive variables are the following
\begin{equation}
\begin{aligned} 
& \left(\rho,u,v,w,p,B_x,B_y,B_z\right)
= \\
& \left( \gamma^2, -\sin(y), \sin(x), 0, \gamma, -\sqrt{4\pi} \sin(y), \sqrt{4\pi} \sin(2x), 0\right)
\label{eqn.Orszag-IC}                                                   
\end{aligned} 
\end{equation}
with $\gamma=5/3$. 
The divergence cleaning speed is set to $c_h=2$ and the final time of the simulation is taken to be $t_f=3$ as in~\cite{Dumbser2008}.

\begin{figure}[h] 
	\centering
	\includegraphics[width=0.33\linewidth]{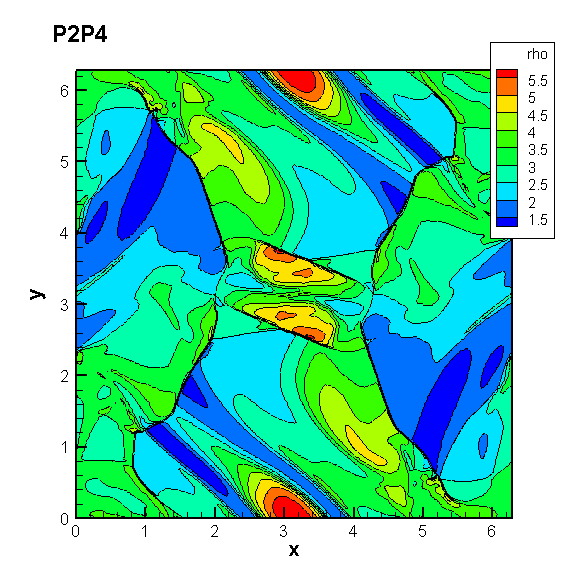}%
	\includegraphics[width=0.33\linewidth]{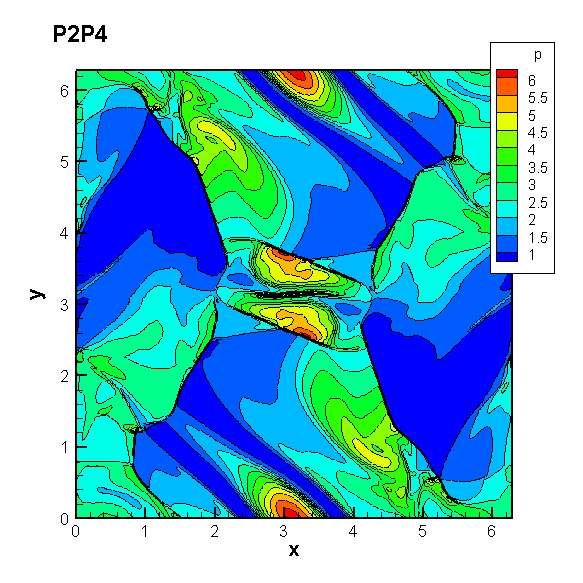}%
	\includegraphics[width=0.33\linewidth]{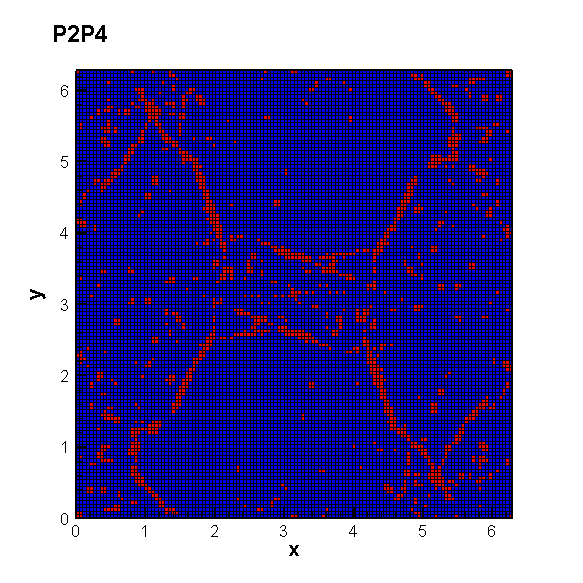}\\
	\includegraphics[width=0.33\linewidth]{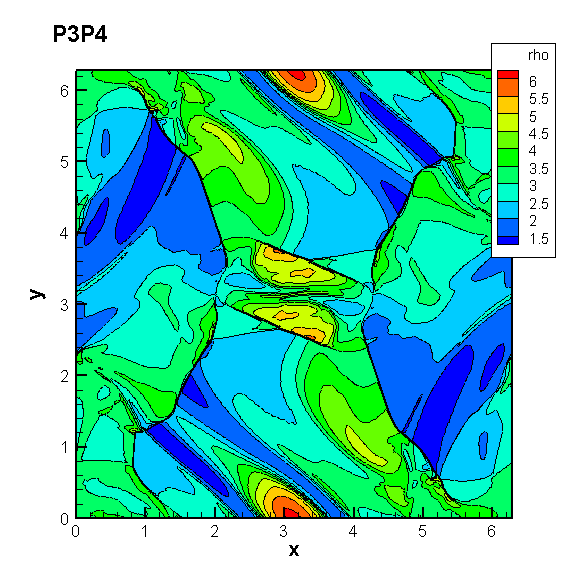}%
	\includegraphics[width=0.33\linewidth]{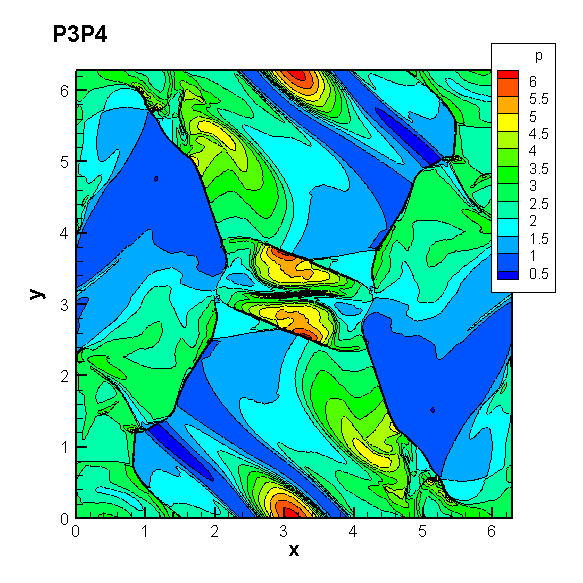}%
	\includegraphics[width=0.33\linewidth]{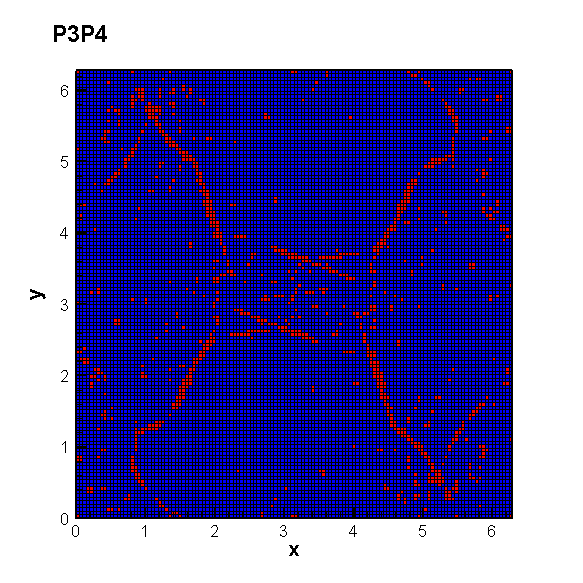}\\
	\includegraphics[width=0.33\linewidth]{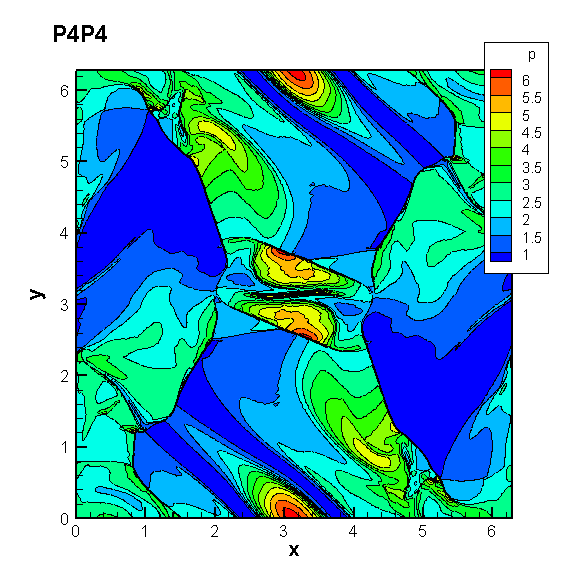}%
	\includegraphics[width=0.33\linewidth]{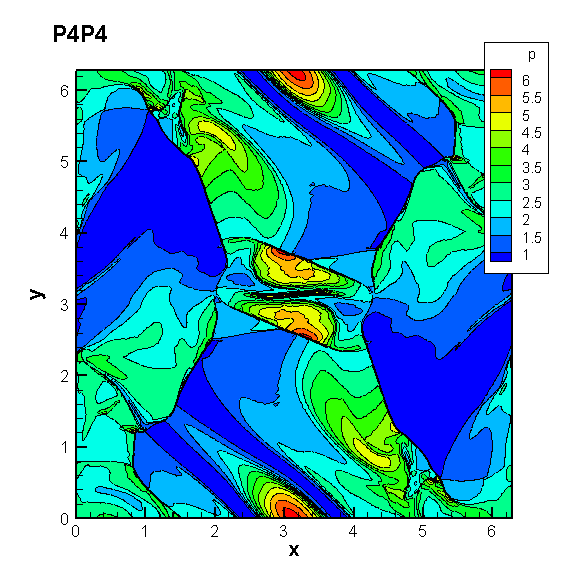}%
	\includegraphics[width=0.33\linewidth]{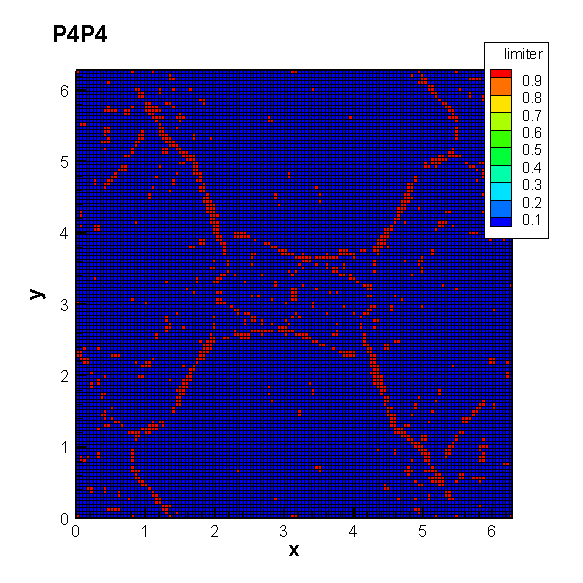}%
	\caption{MHD Orszag-Tang at the final time of $t_f = 3.0$. We compare the results obtained with three fifth order schemes, 
		namely the $P_2P_4$, $P_3P_4$, $P_4P_4$ schemes. In the first column we depict the density contours, 
		in the second column we depict the pressure contours and in the third column we depict the troubled cells in red and the unlimited cells in blue. 
		One can notice that the three methods lead to similar results but with respect to $P_4P_4$ the $P_2P_4$ is $3.39$ times faster and the $P_3P_4$ is $1.60$ time faster.} 
	\label{fig:MHDOrszagTang_P2P4_128x128_rho}
\end{figure}

We solve this test by employing three different fifth order schemes, namely the hybrid $\PN{2}{4}$ and $\PN{3}{4}$ schemes and the pure DG $\PN{4}{4}$ scheme,
with the Rusanov numerical flux and equipped with our \textit{a posteriori } subcell TVD finite volume limiter.
The obtained numerical results and the cells on which the limiter is activated are presented in Figure~\ref{fig:MHDOrszagTang_P2P4_128x128_rho}. 
One can notice that the three methods produce similar results with a good qualitative agreement compared to the solutions provided in~\cite{Dumbser2008,MHDdivFree2015,DGLimiter2}; 
moreover, the hybrid schemes are computationally more efficient than the pure DG scheme.

%

\subsection{Special relativistic MHD equations}

The system of equations of special relativistic magnetohydrodynamics (RMHD) is supposed to provide a sufficiently accurate description of the dynamics 
of those astrophysical plasma that move close to the speed of light and which are
subject to electromagnetic forces that dominate over the gravitational forces.
For example this is the case of high energy astrophysical phenomena like
extragalactic jets~\cite{Begelman1984}, gamma-ray bursts~\cite{Kouveliotou1993} or magnetospheres of neutron stars~\cite{Michel1991}. 

For a more detailed description of this model and a review of the numerical methods used in its approximation we refer to~\cite{Zanotti2015d} and the reference therein. 
Here we briefly recall only the main terms appearing in the equations, which indeed can be written under the general hyperbolic form~\eqref{eq.generalform} by choosing
\be
{\bf Q}=\left[\begin{array}{c}
	D \\ S_j \\ U \\ B^j 
\end{array}\right]\!,\quad
{\bf f}^i=\left[\begin{array}{c}
	v^i D \\
	W^i_j \\
	S^i \\
	\epsilon^{jik}E^k 
\end{array}\right]\!,
\label{eq:fluxes}
\ee
where we have employed the classical tensor index notation based on the Einstein summation convention, which implies summation over two equal indices.

The conserved variables $(D,S_j,U,B^j)$ 
are related to the rest-mass density $\rho$, to the thermal
pressure $p$, to the fluid velocity $v_i$ and to the magnetic field $B^i$ 
by
\be
\label{eq:cons1}
&D   = \rho W,\\
&S_i = \rho h W^2 v_i + \epsilon_{ijk}E_j B_k, \\
&U   = \rho h W^2 - p + \frac{1}{2}(E^2 + B^2),
\ee
where $\epsilon_{ijk}$ is  the spatial Levi--Civita tensor and $\delta_{ij}$ is the Kronecker symbol. As usual in ideal MHD, the electric field is given by $\mathbf{E} = -\mathbf{v} \times \mathbf{B}$.        
The spatial tensor $W^i_j$ in (\ref{eq:fluxes}), representing the momentum  flux density, is 
\be
W_{ij} \equiv \rho h W^2 v_i v_j - E_i E_j - B_i B_j + \left[p +\frac{1}{2}(E^2+B^2)\right]\delta_{ij}, \\
\label{eq:W} 
\ee
where $\delta_{ij}$ is the Kronecker delta.

The above equations include the divergence free condition $\vec \nabla\cdot\vec B=0$ for the magnetic field, 
which, although is guaranteed by the Maxwell equations at a continuous level, 
is not automatically satisfied from a numerical point of view.
Different strategies can be adopted in order to solve this problem (see~\cite{Toth2000} for a review). 
Here, as for the MHD case of Section~\ref{sec.MHD}, we have  adopted the so called \textit{divergence-cleaning approach} 
presented in~\cite{MunzCleaning,Dedneretal}, which considers an augmented system with an additional equation for a scalar 
field $\Phi$, in order to propagate away the deviations from $\vec \nabla\cdot\vec B=0$
\be
\label{eq:divB}
\partial_t \Phi + \partial_i B^i = -\kappa \Phi\,,
\ee
while the fluxes for the evolution of the magnetic field are also modified, namely ${\bf f}^i(B^j)\rightarrow \epsilon^{jik}E^k + \Phi \delta^{ij}$.

\subsubsection{Alfv\'en wave}

As for the previous set of equations we first of all check the convergence of our new numerical scheme.
In the case of RMHD we can consider the propagation of a circularly polarized Alfven wave, 
for which an analytic solution can be computed, see~\cite{Komissarov1997,DelZanna2007}. 

As initial condition we impose the following profile for the magnetic field and the velocity field 
\be
B_x &= B_0 \\
B_y &= \eta B_0\cos[k(x-v_A t)]\\
B_z &= \eta B_0\sin[k(x-v_A t)]\\
v_x &= 0\\
v_y &=-v_A B_y/B_0\\
v_z &=-v_A B_z/B_0.
\ee
where $B_0$ is the uniform magnetic field along $x$, $\rho=p=B_0=\eta=1$, 
$k$ is the wave number, while $v_A$ is the Alfven speed at which the wave propagates
\[
v_A^2 = \frac{B_0^2}{\rho h + B_0^2\left( 1 + \eta^2 \right)} \left[ \frac{1}{2} \left( 1 + \sqrt{1-\left( \frac{2 \eta B_0^2}{\rho h + B_0^2\left(1+\eta^2\right)}   \right )^2        }        \,       \right ) \right ]^{-1}
\]
and $\gamma=5/3$.

For the computational domain, we consider the 2D square $\Omega=[0, 2\pi]\times[0, 2\pi]$
with periodic boundary conditions set everywhere, and we run our simulation up to the final time  $t_f=L/v_A=2\pi/v_A$ corresponding to one period.

In Table~\ref{tab.orderOfconvergenceFV_Alfven}, we report the $L_2$ norm of the errors between our numerical results and the analytical solution for the variable $\rho$.
The convergence rates from second up to sixth order of accuracy are confirmed 
both for the hybrid schemes $\PN{N}{M}$ with $M>N$ and for the pure DG schemes with $N=M$.

\begin{landscape} 
	\begin{table*}
		\caption{Numerical convergence table for general $P_NP_M$ schemes for the relativistic Alfv\'en wave. The error norms refer to the variable $\rho$ at time $t_f~=~2\pi/v_A$ in $L_2$ norm.} 
		\label{tab.orderOfconvergenceFV_Alfven}
		\begin{center} 	
			\begin{tabular}{|c|ccc|ccc|ccc|ccc|ccc|ccc|} 
				\hline 
				h                        & \CPU      & $L_2$             & \OLdue        &   \CPU      & $L_2$             & \OLdue       &     \CPU      & $L_2$             & \OLdue       &      \CPU      & $L_2$             & \OLdue       &    \CPU      & $L_2$             & \OLdue       &        \CPU     & $L_2$            & \OLdue          \\[1pt]
				\hline                                                                                                                                                                                                                                                                                                                                   
				$\mathbf{\mathcal{O}2}$      &               \mcolt{0}{1}                    &                \mcolt{1}{1}                    &               &                  &               &                &                  &               &              &                   &              &                 &                  &                   \\
				\smm 2.0e-01       \smmm &\sm 3  \sm & \sm 1.9e-01   \sm & \smm     \smm & \sm 9  \sm & \sm 1.2e-03   \sm & \smm     \smm &   \sm     \sm & \sm          \sm & \smm     \smm &    \sm     \sm & \sm          \sm & \smm     \smm &  \sm     \sm & \sm          \sm & \smm     \smm &      \sm    \sm & \sm          \sm & \smm     \smm    \\
				\smm 1.2e-01       \smmm &\sm 13 \sm & \sm 9.0e-02   \sm & \smm 1.6 \smm & \sm 38 \sm & \sm 4.8e-04   \sm & \smm 2.0 \smm &   \sm     \sm & \sm          \sm & \smm     \smm &    \sm     \sm & \sm          \sm & \smm     \smm &  \sm     \sm & \sm          \sm & \smm     \smm &      \sm    \sm & \sm          \sm & \smm     \smm    \\
				\smm 1.0e-01       \smmm &\sm 31 \sm & \sm 6.2e-02   \sm & \smm 1.7 \smm & \sm 84 \sm & \sm 3.1e-04   \sm & \smm 2.0 \smm &   \sm     \sm & \sm          \sm & \smm     \smm &    \sm     \sm & \sm          \sm & \smm     \smm &  \sm     \sm & \sm          \sm & \smm     \smm &      \sm    \sm & \sm          \sm & \smm     \smm    \\
				\smm 8.3e-02       \smmm &\sm 49 \sm & \sm 4.5e-02   \sm & \smm 1.7 \smm & \sm 135\sm & \sm 2.1e-04   \sm & \smm 2.0 \smm &   \sm     \sm & \sm          \sm & \smm     \smm &    \sm     \sm & \sm          \sm & \smm     \smm &  \sm     \sm & \sm          \sm & \smm     \smm &      \sm    \sm & \sm          \sm & \smm     \smm    \\[1pt]
				\hline                                                                                                                                                                                                                                                                                                                               
				$\mathbf{\mathcal{O}3}$      &               \mcolt{0}{2}                    &                \mcolt{1}{2}                    &                  \mcolt{2}{2}                    &               &                   &               &             &                   &               &                 &                   &                    \\
				\smm 5.0e-01       \smmm &\sm 0.7\sm & \sm 4.4e-02   \sm & \smm     \smm & \sm 2.0\sm & \sm 6.0e-04   \sm & \smm     \smm &   \sm 3.8\sm & \sm 9.6e-04   \sm & \smm     \smm &    \sm    \sm & \sm           \sm & \smm     \smm &  \sm    \sm & \sm           \sm & \smm     \smm &      \sm    \sm & \sm           \sm & \smm     \smm    \\
				\smm 3.3e-01       \smmm &\sm 2.4\sm & \sm 1.3e-02   \sm & \smm 2.9 \smm & \sm 6.4\sm & \sm 1.8e-04   \sm & \smm 2.9 \smm &   \sm 12 \sm & \sm 2.9e-04   \sm & \smm 2.8 \smm &    \sm    \sm & \sm           \sm & \smm     \smm &  \sm    \sm & \sm           \sm & \smm     \smm &      \sm    \sm & \sm           \sm & \smm     \smm    \\
				\smm 2.5e-01       \smmm &\sm 5.3\sm & \sm 5.6e-03   \sm & \smm 2.9 \smm & \sm 15 \sm & \sm 7.5e-05   \sm & \smm 3.0 \smm &   \sm 27 \sm & \sm 1.2e-04   \sm & \smm 2.9 \smm &    \sm    \sm & \sm           \sm & \smm     \smm &  \sm    \sm & \sm           \sm & \smm     \smm &      \sm    \sm & \sm           \sm & \smm     \smm    \\
				\smm 2.0e-01       \smmm &\sm 9.3\sm & \sm 2.8e-03   \sm & \smm 3.0 \smm & \sm 26 \sm & \sm 3.9e-05   \sm & \smm 2.9 \smm &   \sm 52 \sm & \sm 6.6e-05   \sm & \smm 2.9 \smm &    \sm    \sm & \sm           \sm & \smm     \smm &  \sm    \sm & \sm           \sm & \smm     \smm &      \sm    \sm & \sm           \sm & \smm     \smm    \\[1pt]
				\hline                                                                                                                                                                                                                                                                                                                               
				$\mathbf{\mathcal{O}4}$      &               \mcolt{0}{3}                    &                \mcolt{1}{3}                    &                  \mcolt{2}{3}                    &                   \mcolt{3}{3}                    &             &                   &               &                 &                   &                    \\
				\smm 5.0e-01       \smmm &\sm 1.7\sm & \sm 1.8e-03   \sm & \smm     \smm & \sm 4  \sm & \sm 2.1e-04   \sm & \smm     \smm &   \sm 8  \sm & \sm 1.1e-05   \sm & \smm     \smm &    \sm 12 \sm & \sm 5.2e-06   \sm & \smm     \smm &  \sm    \sm & \sm           \sm & \smm     \smm &      \sm    \sm & \sm           \sm & \smm     \smm    \\
				\smm 3.3e-01       \smmm &\sm 5.5\sm & \sm 1.9e-04   \sm & \smm 5.4 \smm & \sm 13 \sm & \sm 4.3e-05   \sm & \smm 3.9 \smm &   \sm 27 \sm & \sm 1.6e-06   \sm & \smm 4.8 \smm &    \sm 41 \sm & \sm 7.6e-07   \sm & \smm 4.7 \smm &  \sm    \sm & \sm           \sm & \smm     \smm &      \sm    \sm & \sm           \sm & \smm     \smm    \\
				\smm 2.5e-01       \smmm &\sm 13 \sm & \sm 4.3e-05   \sm & \smm 5.3 \smm & \sm 31 \sm & \sm 1.3e-05   \sm & \smm 4.0 \smm &   \sm 65 \sm & \sm 4.3e-07   \sm & \smm 4.6 \smm &    \sm 113\sm & \sm 2.5e-07   \sm & \smm 3.8 \smm &  \sm    \sm & \sm           \sm & \smm     \smm &      \sm    \sm & \sm           \sm & \smm     \smm    \\
				\smm 2.0e-01       \smmm &\sm 27 \sm & \sm 1.3e-05   \sm & \smm 5.1 \smm & \sm 76 \sm & \sm 5.6e-06   \sm & \smm 3.9 \smm &   \sm 127\sm & \sm 1.5e-07   \sm & \smm 4.6 \smm &    \sm 205\sm & \sm 1.1e-07   \sm & \smm 3.4 \smm &  \sm    \sm & \sm           \sm & \smm     \smm &      \sm    \sm & \sm           \sm & \smm     \smm    \\[1pt]
				\hline                                                                                                                                                                                                                                                                                                                               
				$\mathbf{\mathcal{O}5}$      &               \mcolt{0}{4}                    &                \mcolt{1}{4}                    &                  \mcolt{2}{4}                    &                   \mcolt{3}{4}                    &                 \mcolt{4}{4}                    &                 &                  &                     \\
				\smm 5.0e-01       \smmm &\sm 3  \sm & \sm 8.8e-04   \sm & \smm     \smm & \sm 8  \sm & \sm 5.9e-06   \sm & \smm     \smm &   \sm 16 \sm & \sm 8.1e-06   \sm & \smm     \smm &    \sm 23  \sm & \sm 1.6e-07   \sm & \smm    \smm &  \sm 33  \sm & \sm 2.6e-07   \sm & \smm    \smm &      \sm    \sm & \sm          \sm & \smm    \smm      \\
				\smm 4.0e-01       \smmm &\sm 7  \sm & \sm 2.9e-04   \sm & \smm 4.9 \smm & \sm 16 \sm & \sm 1.9e-06   \sm & \smm 4.9 \smm &   \sm 33 \sm & \sm 2.7e-06   \sm & \smm 4.8 \smm &    \sm 49  \sm & \sm 5.9e-08   \sm & \smm 4.4\smm &  \sm 70  \sm & \sm 8.8e-08   \sm & \smm 4.8\smm &      \sm    \sm & \sm          \sm & \smm    \smm      \\
				\smm 3.3e-01       \smmm &\sm 10 \sm & \sm 1.1e-04   \sm & \smm 4.9 \smm & \sm 25 \sm & \sm 7.9e-07   \sm & \smm 4.9 \smm &   \sm 69 \sm & \sm 1.1e-06   \sm & \smm 4.9 \smm &    \sm 97  \sm & \sm 2.2e-08   \sm & \smm 5.4\smm &  \sm 132 \sm & \sm 3.7e-08   \sm & \smm 4.7\smm &      \sm    \sm & \sm          \sm & \smm    \smm      \\
				\smm 2.8e-01       \smmm &\sm 17 \sm & \sm 5.4e-05   \sm & \smm 4.9 \smm & \sm 42 \sm & \sm 3.6e-07   \sm & \smm 5.0 \smm &   \sm 90 \sm & \sm 5.2e-07   \sm & \smm 4.8 \smm &    \sm 129 \sm & \sm 9.7e-09   \sm & \smm 5.2\smm &  \sm 203 \sm & \sm 1.7e-08   \sm & \smm 5.0\smm &      \sm    \sm & \sm          \sm & \smm    \smm      \\[1pt]
				\hline                                                                                                                                                                                                                                                                                                                               
				$\mathbf{\mathcal{O}6}$      &               \mcolt{0}{5}                    &                \mcolt{1}{5}                    &                  \mcolt{2}{5}                    &                   \mcolt{3}{5}                    &                 \mcolt{4}{5}                    &                     \mcolt{5}{5}                       \\
				\smm 1.2e+00       \smmm &\sm 0.8\sm & \sm 1.8e-02   \sm & \smm     \smm & \sm 1.0\sm & \sm 1.5e-04   \sm & \smm     \smm &   \sm 1.8\sm & \sm 2.5e-05   \sm & \smm     \smm &    \sm 3.0\sm & \sm 1.3e-05   \sm & \smm     \smm &  \sm 4.6\sm & \sm 1.0e-06   \sm & \smm     \smm &      \sm 7.6\sm & \sm 2.6e-07   \sm & \smm     \smm    \\			
				\smm 1.0e+00       \smmm &\sm 1.0\sm & \sm 3.1e-03   \sm & \smm 8.0 \smm & \sm 2.0\sm & \sm 3.9e-05   \sm & \smm 6.1 \smm &   \sm 3.6\sm & \sm 5.8e-06   \sm & \smm 6.5 \smm &    \sm 5.8\sm & \sm 3.5e-06   \sm & \smm 5.9 \smm &  \sm 9  \sm & \sm 2.2e-07   \sm & \smm 6.7 \smm &      \sm 11 \sm & \sm 6.3e-08   \sm & \smm 6.3 \smm    \\
				\smm 8.3e-01       \smmm &\sm 2.3\sm & \sm 7.9e-04   \sm & \smm 7.4 \smm & \sm 3.8\sm & \sm 1.2e-05   \sm & \smm 6.2 \smm &   \sm 6.5\sm & \sm 1.8e-06   \sm & \smm 6.2 \smm &    \sm 10 \sm & \sm 1.2e-06   \sm & \smm 6.0 \smm &  \sm 21 \sm & \sm 6.9e-08   \sm & \smm 6.3 \smm &      \sm 19 \sm & \sm 2.1e-08   \sm & \smm 6.0 \smm    \\
				\smm 7.1e-01       \smmm &\sm 3.4\sm & \sm 2.6e-04   \sm & \smm 7.1 \smm & \sm 6.6\sm & \sm 5.0e-06   \sm & \smm 6.0 \smm &   \sm 12 \sm & \sm 6.9e-07   \sm & \smm 6.3 \smm &    \sm 16 \sm & \sm 4.4e-07   \sm & \smm 6.3 \smm &  \sm 30 \sm & \sm 2.6e-08   \sm & \smm 6.2 \smm &      \sm 34 \sm & \sm 8.7e-09   \sm & \smm 5.7 \smm    \\
				\hline 
			\end{tabular}		
		\end{center}
	\end{table*}
\end{landscape}

\subsubsection{Riemann problems}
\begin{table}[t] 
	\caption{
		\label{tab:RP1D}
		Initial conditions for the one--dimensional Riemann problems.} 
	\begin{center} 
		\begin{tabular}{|c|c||ccccc|c|c|} 
			\hline
			Problem    && $\rho$ &$(v_x$&$v_y$&$v_z)$ & $p$ & $t_f$ & $\gamma$ \\
			\hline
			\multirow{2}{*}{\rotatebox{0}{\textbf{RP1}}} 
			&$x \leq 0$     & 1   &  0.9   & 0   &  0   & 1  &    \multirow{2}{*}{0.4}&\multirow{2}{*}{$\frac{5}{3}$}\\ 
			&$x  > 0$  & 1   &  0     & 0   &  0   & 10 &                        &\\ 
			\hline
			\multirow{2}{*}{\rotatebox{0}{\textbf{RP2}}}
			&$x \leq 0$     & 1   &  -0.6   & 0   &  0   & 10  &    \multirow{2}{*}{0.4}&\multirow{2}{*}{$\frac{4}{3}$}\\ 
			&$x > 0$  & 10  &   0.5   & 0   &  0   & 20 &                        &\\ 
			\hline
		\end{tabular} 		
	\end{center}
\end{table} 
\begin{figure}[b]
	\centering
	\includegraphics[width=0.49\linewidth]{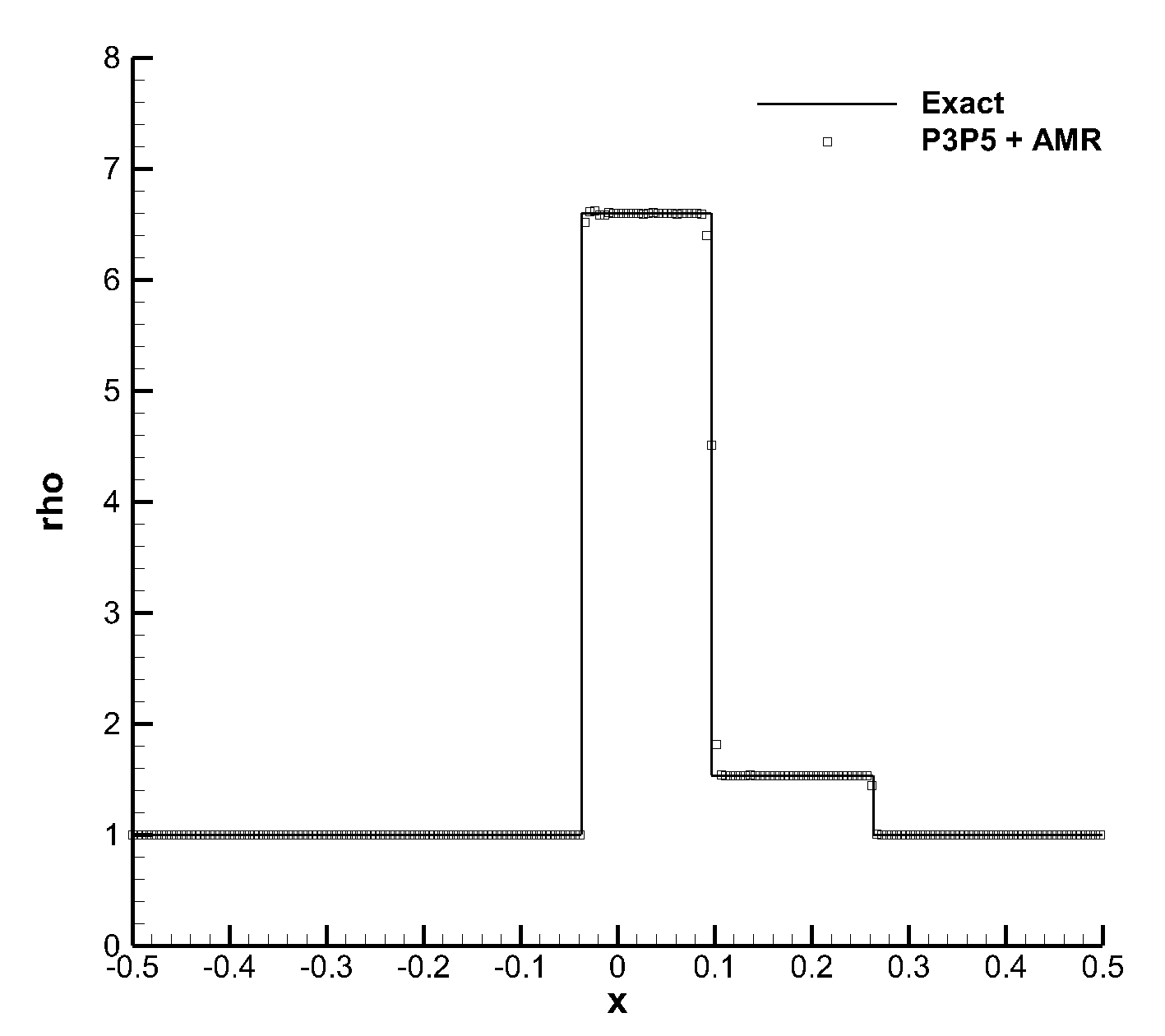}
	\includegraphics[width=0.49\linewidth]{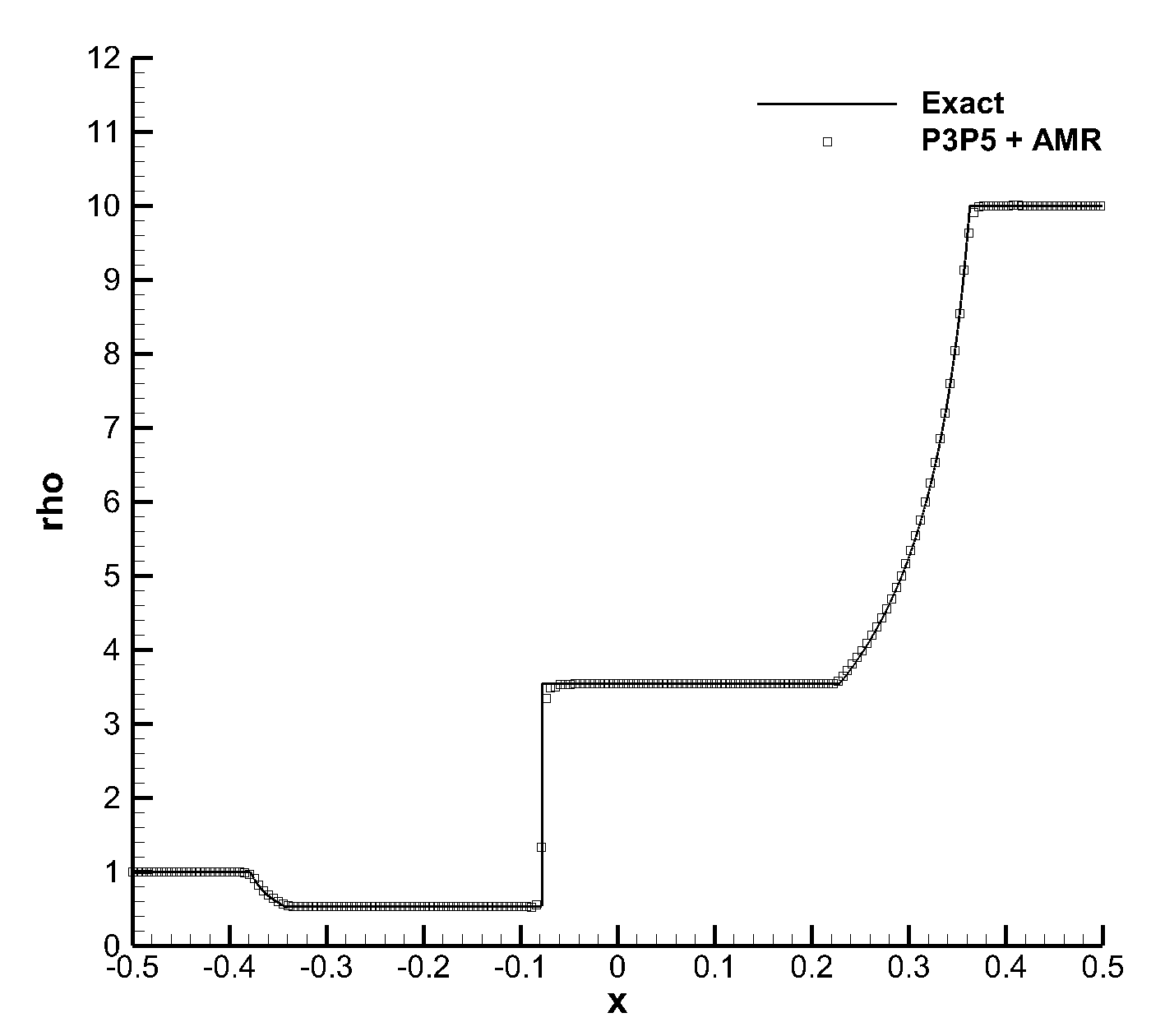}
	\caption{RHD equations and Riemann problems \textbf{RP1} (left) and \textbf{RP2} (right), see Table~\ref{tab:RP1D} for the initial conditions. 
		In the Figure, we compare our numerical results for the density $\rho$ (squares) with the exact solution (continuous line).
	}
	\label{fig:rhdrp1vsexact}
\end{figure}

Now, in order to check the robustness and accuracy of our \textit{a posteriori} subcell FV limiter for the general class of the $\PN{N}{M}$ schemes, 
we solve two classical Riemann problems of RHD (i.e. RMHD with $\B=\0$) for which also an exact solution is available.

We consider the computational domain $\Omega=[-0.5,0.5]\times[0,1]$ and as initial condition we impose the discontinuous values given in Table~\ref{tab:RP1D}.
We solve these two test cases with a fifth order $\PN{3}{5}$ scheme and the HLLEM numerical flux,
over a mesh of $20\times10$ elements with $\ell_{\max}=2$ levels of refinement and $\rRef=3$. 

The obtained numerical results, see Figure~\ref{fig:rhdrp1vsexact}, show once again that our limiter procedure preserves the resolution of the underlying $\PN{N}{M}$ scheme even on a coarse mesh.

%
%

\subsubsection{Cylindrical blast wave}
\begin{figure}[h]
	\centering
	\begin{tabular}{cc} 
		\includegraphics[width=0.45\linewidth]{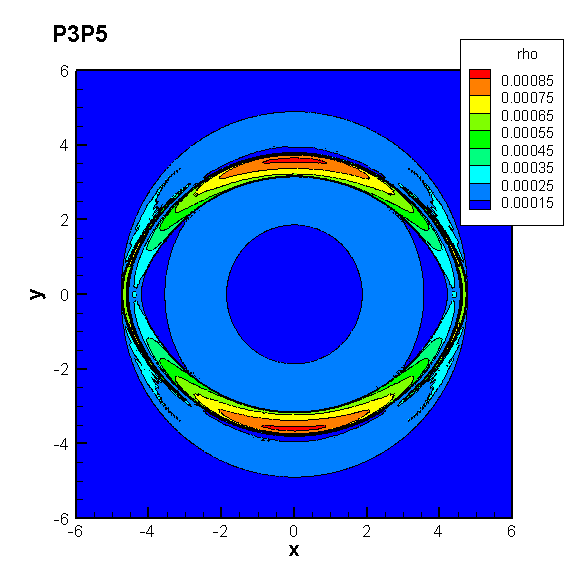} & 
		\includegraphics[width=0.45\linewidth]{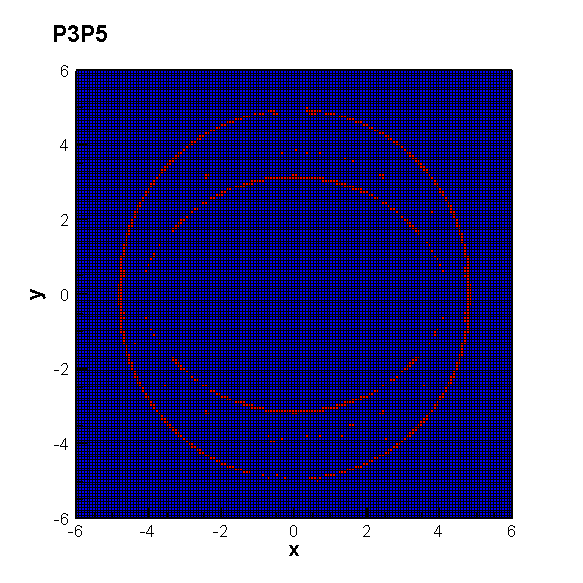} \\ 
		\includegraphics[width=0.45\linewidth]{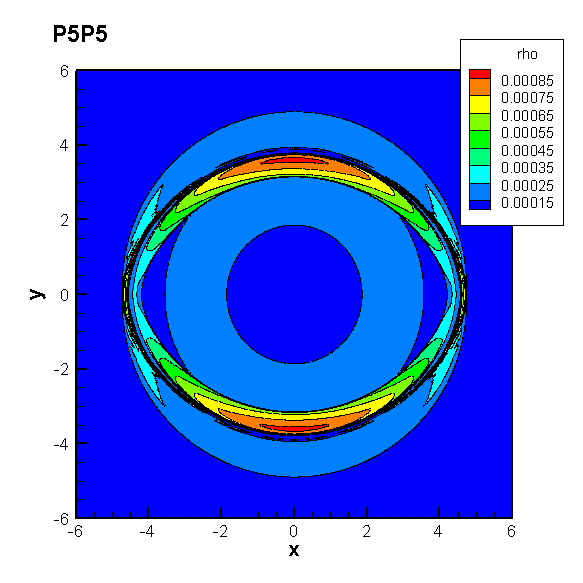}  & 
		\includegraphics[width=0.45\linewidth]{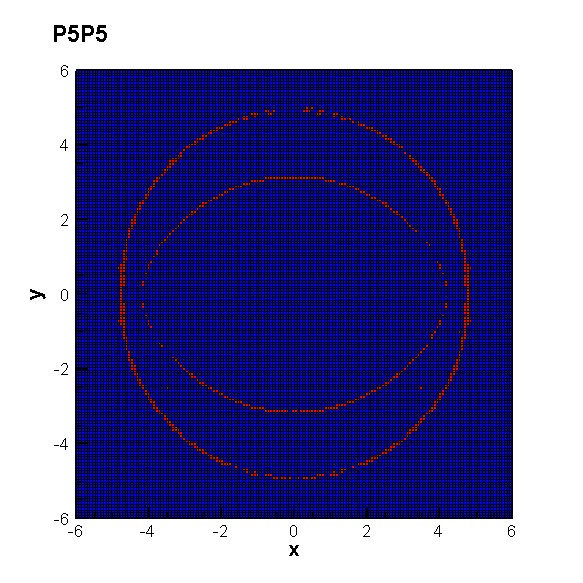}
	\end{tabular} 
	\caption{RMHD blast wave at time $t_f=4.0$. We show the results obtained with two sixth order schemes, namely the $P_3P_5$ and $P_5P_5$ schemes. 
		In the left column we plot the density contours and in the right column we depict the troubled cells in red and the unlimited cells in blue. 	
		The resolution and the number of limited cells are quite similar with the two approaches but the hybrid $\PN{3}{5}$ scheme is $2.79$ times faster than the pure DG $\PN{5}{5}$ scheme.}
	\label{fig:RMHD_Blast_P3P5_160x160_Ref0_rho}
\end{figure}

We now take into account a truly two dimensional test in RMHD, i.e. the cylindrical expansion of a blast wave in a plasma with an initially uniform magnetic field.  
This is a severe test proposed in~\cite{Komissarov1999}, and subsequently also solved in  \cite{Leismann2005,DelZanna2007,DumbserZanotti,Zanotti2015d}. 

For the initial condition we set 
the rest-mass density and the pressure equal to $\rho=0.01$ and $p=1$
within a cylinder of radius $r=1.0$, and $\rho=10^{-4}$ and $p=5\times10^{-4}$ outside.
Like in~\cite{Komissarov1999} and in~\cite{DelZanna2007}, the inner and outer values are joined
through a smooth ramp function between $r=0.8$ and $r=1$, to avoid a sharp discontinuity in the initial conditions.
The plasma is initially at rest and subject to a constant magnetic field along the $x$-direction, i.e. $B_x=0.1, B_y=0, B_z=0$.

We have solved this problem on the computational domain 
$\Omega = [-6,6]\times[-6,6]$, with a uniform mesh of $160\times 160$ elements.
We have used the Rusanov numerical flux and two sixth order schemes, namely the $\PN{3}{5}$ and the $\PN{5}{5}$ schemes, equipped with 
the robust \textit{a posteriori} subcell second-order TVD FV limiter. 
The obtained numerical results, which agree with those available in the literature, are reported in Figure~\ref{fig:RMHD_Blast_P3P5_160x160_Ref0_rho}.



\subsubsection{RMHD Orszag-Tang vortex}
\begin{figure}[b]
	\centering
	\includegraphics[width=0.45\linewidth]{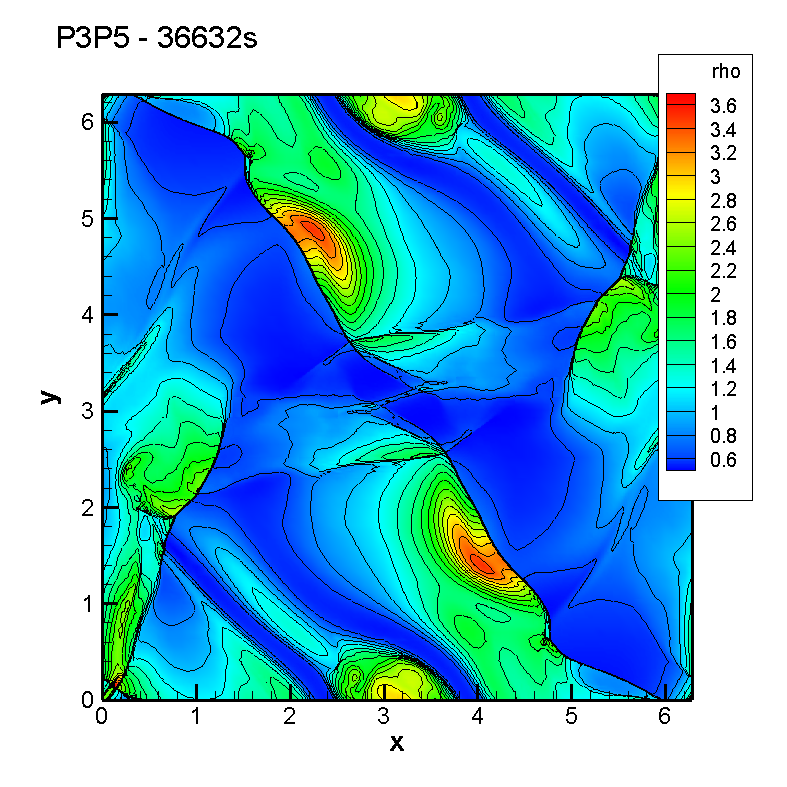}
	\includegraphics[width=0.45\linewidth]{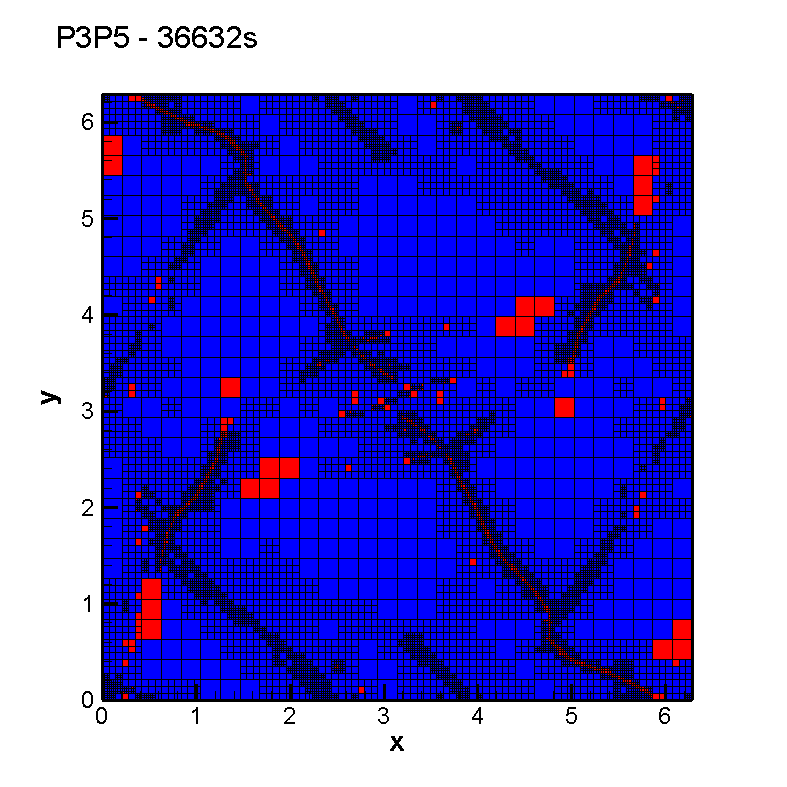}\\
	\includegraphics[width=0.45\linewidth]{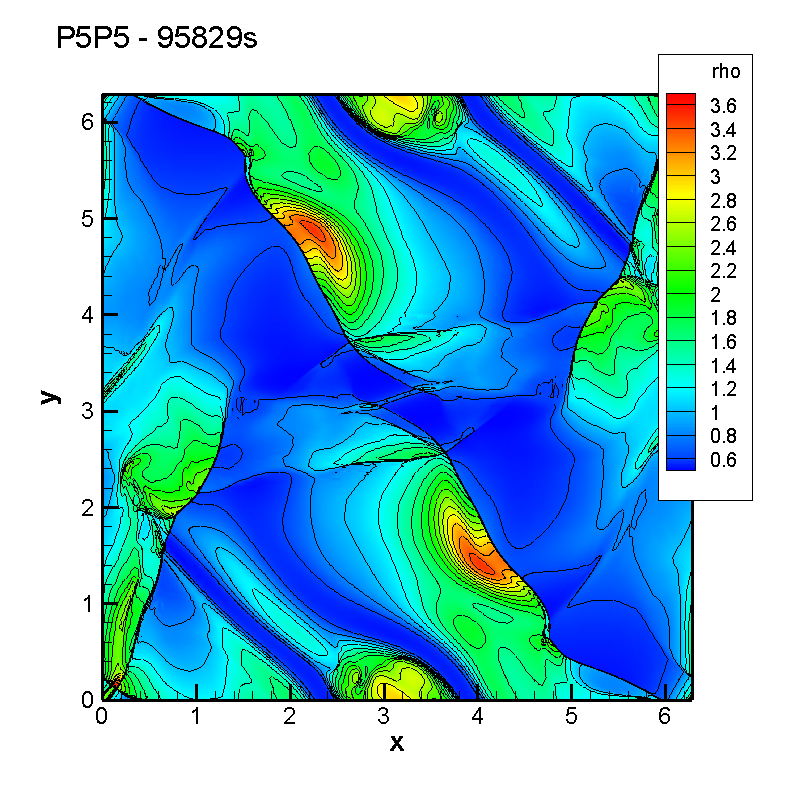}
	\includegraphics[width=0.45\linewidth]{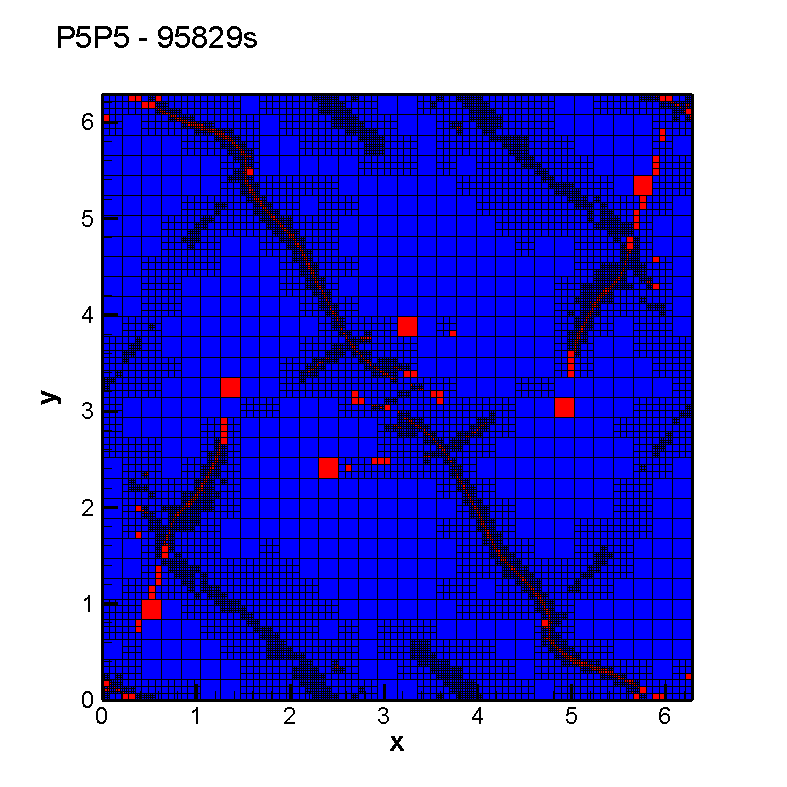}
	\caption{RMHD Orszag Tang at time $t_f=5$. We present the density contours (left column) and the limited cells (in red in the right column) obtained with two sixth order schemes, 
		namely the $\PN{3}{5}$ and the $\PN{5}{5}$ schemes, on an adaptive mesh with $45\times45$ control volume on the coarsest level, $\ell_{\max}=2$ and $\rRef=3$.}
	\label{fig:rmhdotp3p5}
\end{figure}

Finally, we have chosen the relativistic version of the well known Orszag-Tang vortex problem, proposed by~\cite{OrszagTang}, 
and adapted to the relativistic case in~\cite{DumbserZanotti,Zanotti2015d}.
The computational domain is $\Omega = [0,2\pi]\times[0,2\pi]$ and the initial conditions are given by
\be
& \left( \rho, u, v, w, p, B_x ,B_y,B_z \right) = \\
& \left(  1 , - \frac{3}{4\sqrt{2}}\sin\left(y\right), \frac{3}{4\sqrt{2}}\sin \left(x \right), 0, 1, - \sin\left(y\right), \sin \left(2x \right), 0 \right),
\label{eq:OrszagTang_ic}
\ee
with $\gamma=4/3$. 

To solve this system we employ the sixth order hybrid scheme $\PN{3}{5}$ over a level zero mesh of $45\times45$ elements, activating the AMR feature with $\ell_{\max}=2$ levels of refinement and $\rRef=3$. 
The obtained numerical results are reported in Figure~\ref{fig:rmhdotp3p5}: once again we can notice 
that the two schemes have a similar resolution but the hybrid scheme is $2.62$ times faster than a pure DG scheme, of the same order. Furthermore, we can observe that the proposed \textit{a posteriori} subcell limiter procedure is robust and maintains the high resolution of the underlying $\PN{N}{M}$ scheme even on coarse meshes. 
%
%

\section{Conclusion}
\label{sec.conclusion}

In this paper we have proposed a new simple, robust, accurate and computationally efficient limiting strategy for the general family of ADER $\PN{N}{M}$ schemes,
allowing, for the first time in literature, the use of hybrid reconstructed methods ($N>0, M>N$) in the modeling of discontinuous phenomena.
The key ideas behind our limiter are: 
i) its local activation only where the linear schemes introduces oscillations through an \textit{a posteriori} detector, 
ii) its robustness due to the use of strong stability preserving TVD or WENO FV schemes as limiter, 
iii) its resolution due to the use of the limiter on $2N+1$ subcells.
Thus, we have been able to apply this new approach to many different systems of hyperbolic conservation laws, providing highly accurate numerical results in all cases.
Moreover, we had the possibility to compare the performance of the class of intermediate $P_NP_M$ schemes with $M>N>0$ with pure DG schemes ($M=N$). We have observed that in most cases the intermediate $P_NP_M$ schemes offer a similar resolution compared to pure DG methods, but at a reduced computational cost. 

Future work will consider the extension of $\PN{N}{M}$ scheme with $N~>~0$, $M~>~N$ to unstructured moving meshes~\cite{Lagrange2D,GaburroNonConforming}, 
in particular for regenerating Voronoi tessellations~\cite{GaburroAREPO,gaburroReview}, and to the three-dimensional case.
Finally, due to their low memory consumption and their gain in computational efficiency compared to DG schemes, 
they will be also considered for astrophysical applications~\cite{gaburro2018well,ADERCCZ4,dumbser2020glm}
and the unified first order hyperbolic model of continuum mechanics proposed in  \cite{PeshRom2014,GPRmodel,GPRmodelMHD,FrontierADERGPR}, where a large number of conserved variables has to be discretized. Due to their accuracy and compact stencil in the future we also plan to use $\PN{N}{M}$ schemes with \textit{a posteriori} subcell finite volume limiter in the context of hyperbolic reformulations of nonlinear dispersive systems and wave propagation problems, see e.g.  \cite{Dhaouadi2018,favrie:2017,EM20,DispersiveSWE,DispersiveCoupling}.

\section*{Acknowledgments}

The research presented in this paper has been partially financed by the European Research Council (ERC) under the 
European Union's Seventh Framework Programme (FP7/2007-2013) with the research project \textit{STiMulUs},  
ERC Grant agreement no. 278267. 

E.~G. has been supported by a national mobility grant for young researchers in Italy, funded by GNCS-INdAM. E.~G. has also received funding from the University of Trento via the Strategic 
Initiative \textit{Starting Grant Giovani Ricercatori 2019}.

M.~D. has been funded by the European Union's Horizon 2020 Research and Innovation  Programme under the project \textit{ExaHyPE}, grant no. 671698 (call FETHPC-1-2014). M.~D. also acknowledges the financial support received from 
the Italian Ministry of Education, University and Research (MIUR) in the frame of the Departments of Excellence  Initiative 2018--2022 attributed to DICAM of the University of Trento (grant L. 232/2016) and in the frame of the 
PRIN 2017 project \textit{Innovative numerical methods for evolutionary partial differential equations and  applications}. Furthermore, M.~D. has also received funding from the University of Trento via the Strategic 
Initiative \textit{Modeling and Simulation}.

%
%

\bibliographystyle{spmpsci}      
\bibliography{references}  

\begin{thebibliography}{100}
\providecommand{\url}[1]{{#1}}
\providecommand{\urlprefix}{URL }
\expandafter\ifx\csname urlstyle\endcsname\relax
  \providecommand{\doi}[1]{DOI~\discretionary{}{}{}#1}\else
  \providecommand{\doi}{DOI~\discretionary{}{}{}\begingroup
  \urlstyle{rm}\Url}\fi

\bibitem{Alic:2009}
{Alic}, D., {Bona}, C., {Bona-Casas}, C.: {Towards a gauge-polyvalent numerical
  relativity code}.
\newblock Phys. Rev. D \textbf{79}(4), 044026 (2009)

\bibitem{Alic:2012}
{Alic}, D., {Bona-Casas}, C., {Bona}, C., {Rezzolla}, L., {Palenzuela}, C.:
  {Conformal and covariant formulation of the Z4 system with
  constraint-violation damping}.
\newblock Phys. Rev. D \textbf{85}(6) (2012)

\bibitem{Aloy1999c}
{Aloy}, M.A., {Ib{\'a}{\~n}ez}, J.M., {Mart{\'{\i}}}, J.M., {M{\"u}ller}, E.:
  {GENESIS: A High-Resolution Code for Three-dimensional Relativistic
  Hydrodynamics}.
\newblock Astrohys. J. Suppl. \textbf{122}, 151--166 (1999)

\bibitem{BaerNunziato1986}
Baer, M., Nunziato, J.: A two-phase mixture theory for the
  deflagration-to-detonation transition {(DDT)} in reactive granular materials.
\newblock J. Multiphase Flow \textbf{12}, 861--889 (1986)

\bibitem{BalsaraRMHD}
Balsara, D.: Total variation diminishing scheme for relativistic
  magnetohydrodynamics.
\newblock The Astrophysical Journal Supplement Series \textbf{132}, 83--101
  (2001)

\bibitem{Balsara2004}
Balsara, D.: Second-order accurate schemes for magnetohydrodynamics with
  divergence-free reconstruction.
\newblock The Astrophysical Journal Supplement Series \textbf{151}, 149--184
  (2004)

\bibitem{MHDdivFree2015}
Balsara, D., Dumbser, M.: {Divergence-free MHD on unstructured meshes using
  high order finite volume schemes based on multidimensional Riemann solvers}.
\newblock Journal of Computational Physics \textbf{299}, 687 -- 715 (2015)

\bibitem{balsarashu}
Balsara, D., Shu, C.: Monotonicity preserving weighted essentially
  non-oscillatory schemes with increasingly high order of accuracy.
\newblock Journal of Computational Physics \textbf{160}, 405--452 (2000)

\bibitem{BalsaraSpicer1999}
Balsara, D., Spicer, D.: A staggered mesh algorithm using high order godunov
  fluxes to ensure solenoidal magnetic fields in magnetohydrodynamic
  simulations.
\newblock Journal of Computational Physics \textbf{149}, 270--292 (1999)

\bibitem{Banyuls97}
Banyuls, F., Font, J.A., Ib{\'a}{\~n}ez, J.M., Mart{\'\i}, J.M., Miralles,
  J.A.: Numerical {3+1} general-relativistic hydrodynamics: A local
  characteristic approach.
\newblock Astrophys. J. \textbf{476}, 221 (1997)

\bibitem{DispersiveSWE}
Bassi, C., Bonaventura, L., Busto, S., Dumbser, M.: A hyperbolic reformulation
  of the {Serre-Green-Naghdi} model for general bottom topographies.
\newblock Computers \& Fluids \textbf{212}, 104716 (2020)

\bibitem{DispersiveCoupling}
Bassi, C., Busto, S., Dumbser, M.: High order {ADER-DG} schemes for the
  simulation of linear seismic waves induced by nonlinear dispersive
  free-surface water waves.
\newblock Applied Numerical Mathematics \textbf{158}, 236--263 (2020)

\bibitem{Begelman1984}
{Begelman}, M.C., {Blandford}, R.D., {Rees}, M.J.: {Theory of extragalactic
  radio sources}.
\newblock Reviews of Modern Physics \textbf{56}, 255--351 (1984)

\bibitem{Berger-Colella1989}
{Berger}, M.J., {Colella}, P.: {Local adaptive mesh refinement for shock
  hydrodynamics}.
\newblock Journal of Computational Physics \textbf{82}, 64--84 (1989)

\bibitem{Berger-Oliger1984}
{Berger}, M.J., {Oliger}, J.: {Adaptive Mesh Refinement for Hyperbolic Partial
  Differential Equations}.
\newblock Journal of Computational Physics \textbf{53}, 484 (1984)

\bibitem{boscheri2017efficient}
Boscheri, W.: An efficient high order direct ale ader finite volume scheme with
  a posteriori limiting for hydrodynamics and magnetohydrodynamics.
\newblock International Journal for Numerical Methods in Fluids \textbf{84}(2),
  76--106 (2017)

\bibitem{WAO-ALE}
Boscheri, W., Balsara, D.: High order direct arbitrary-lagrangian-eulerian
  (ale) {PNPM} schemes with weno adaptive-order reconstruction on unstructured
  meshes.
\newblock Journal of Computational Physics \textbf{398}, 108899 (2019)

\bibitem{Lagrange2D}
Boscheri, W., Dumbser, M.: {Arbitrary--Lagrangian--Eulerian One--Step WENO
  Finite Volume Schemes on Unstructured Triangular Meshes}.
\newblock Communications in Computational Physics \textbf{14}, 1174--1206
  (2013)

\bibitem{ALEDG}
Boscheri, W., Dumbser, M.: {Arbitrary-Lagrangian-Eulerian discontinuous
  Galerkin schemes with a posteriori subcell finite volume limiting on moving
  unstructured meshes}.
\newblock Journal of Computational Physics \textbf{346}, 449 -- 479 (2017)

\bibitem{boscheri2014high}
Boscheri, W., Dumbser, M., Balsara, D.: High-order ader-weno ale schemes on
  unstructured triangular meshes—application of several node solvers to
  hydrodynamics and magnetohydrodynamics.
\newblock International Journal for Numerical Methods in Fluids
  \textbf{76}(10), 737--778 (2014)

\bibitem{ALEMOOD2}
Boscheri, W., Loub\`ere, R.: {High order accurate direct
  Arbitrary-Lagrangian-Eulerian ADER-MOOD finite volume schemes for
  non-conservative hyperbolic systems with stiff source terms}.
\newblock Communications in Computational Physics \textbf{21}, 271--312 (2017)

\bibitem{ALEMOOD1}
Boscheri, W., Loub\`ere, R., Dumbser, M.: {Direct Arbitrary-Lagrangian-Eulerian
  ADER-MOOD finite volume schemes for multidimensional hyperbolic conservation
  laws}.
\newblock Journal of Computational Physics \textbf{292}, 56--87 (2015)

\bibitem{DGCWENO}
Boscheri, W., Semplice, M., Dumbser, M.: {Central WENO Subcell Finite Volume
  Limiters for ADER Discontinuous Galerkin Schemes on Fixed and Moving
  Unstructured Meshes}.
\newblock Communications in Computational Physics \textbf{25}, 311--346 (2019)

\bibitem{Peano1}
Bungartz, H., Mehl, M., Neckel, T., Weinzierl, T.: {The PDE framework Peano
  applied to fluid dynamics: An efficient implementation of a parallel
  multiscale fluid dynamics solver on octree-like adaptive Cartesian grids}.
\newblock Computational Mechanics \textbf{46}, 103--114 (2010)

\bibitem{FrontierADERGPR}
Busto, S., Chiocchetti, S., Dumbser, M., Gaburro, E., Peshkov, I.: High order
  {ADER} schemes for continuum mechanics.
\newblock Frontiers in Physics \textbf{8}, 32 (2020)

\bibitem{BTVC2016}
Busto, S., Toro, E., V{\'a}zquez-Cend{\'o}n, E.: Design and analysis of
  {ADER}-type schemes for model advection--diffusion--reaction equations.
\newblock Journal of Computational Physics \textbf{327}, 553--575 (2016)

\bibitem{Casulli1990}
Casulli, V.: Semi-implicit finite difference methods for the two-dimensional
  shallow water equations.
\newblock Journal of Computational Physics \textbf{86}, 56--74 (1990)

\bibitem{chiravalle20193d}
Chiravalle, V., Morgan, N.: A {3D Lagrangian} cell-centered hydrodynamic method
  with higher-order reconstructions for gas and solid dynamics.
\newblock Computers \& Mathematics with Applications \textbf{78}(2), 298--317
  (2019)

\bibitem{CDL1}
Clain, S., Diot, S., Loub{\`e}re, R.: A high-order finite volume method for
  systems of conservation laws—multi-dimensional optimal order detection
  ({MOOD}).
\newblock Journal of Computational Physics \textbf{230}(10), 4028 -- 4050
  (2011)

\bibitem{DahlburgPicone}
Dahlburg, R.B., Picone, J.M.: Evolution of the {Orszag--Tang} vortex system in
  a compressible medium. {I.} initial average subsonic flow.
\newblock Phys. Fluids B \textbf{1}, 2153--2171 (1989)

\bibitem{Dedneretal}
Dedner, A., Kemm, F., Kr\"oner, D., Munz, C.D., Schnitzer, T., Wesenberg, M.:
  Hyperbolic divergence cleaning for the {MHD} equations.
\newblock Journal of Computational Physics \textbf{175}, 645--673 (2002)

\bibitem{DelZanna2007}
{Del Zanna}, L., {Zanotti}, O., {Bucciantini}, N., {Londrillo}, P.: {ECHO: a
  Eulerian conservative high-order scheme for general relativistic
  magnetohydrodynamics and magnetodynamics}.
\newblock Astron. Astrophys. \textbf{473}, 11--30 (2007)

\bibitem{Dhaouadi2018}
Dhaouadi, F., Favrie, N., Gavrilyuk, S.: {Extended Lagrangian approach for the
  defocusing nonlinear Schr{\"{o}}dinger equation}.
\newblock Studies in Applied Mathematics pp. 1--20 (2018)

\bibitem{CDL2}
Diot, S., Clain, S., Loub{\`e}re, R.: Improved detection criteria for the
  multi-dimensional optimal order detection ({MOOD}) on unstructured meshes
  with very high-order polynomials.
\newblock Computers and Fluids \textbf{64}, 43 -- 63 (2012)

\bibitem{CDL3}
Diot, S., Loub{\`e}re, R., Clain, S.: The {MOOD} method in the
  three-dimensional case: Very-high-order finite volume method for hyperbolic
  systems.
\newblock International Journal of Numerical Methods in Fluids \textbf{73},
  362--392 (2013)

\bibitem{ADERNSE}
Dumbser, M.: Arbitrary high order {PNPM} schemes on unstructured meshes for the
  compressible {Navier--Stokes} equations.
\newblock Computers \& Fluids \textbf{39}, 60--76 (2010)

\bibitem{dumbser2013diffuse}
Dumbser, M.: A diffuse interface method for complex three-dimensional free
  surface flows.
\newblock Computer Methods in Applied Mechanics and Engineering \textbf{257},
  47--64 (2013)

\bibitem{Dumbser2015}
Dumbser, M., Balsara, D.: A new, efficient formulation of the {HLLEM Riemann}
  solver for general conservative and non-conservative hyperbolic systems.
\newblock Journal of Computational Physics \textbf{304}, 275--319 (2016)

\bibitem{Dumbser2008}
Dumbser, M., Balsara, D., Toro, E., Munz, C.: A unified framework for the
  construction of one--step finite--volume and discontinuous {Galerkin}
  schemes.
\newblock Journal of Computational Physics \textbf{227}, 8209--8253 (2008)

\bibitem{dumbser2020glm}
Dumbser, M., Fambri, F., Gaburro, E., Reinarz, A.: On glm curl cleaning for a
  first order reduction of the ccz4 formulation of the einstein field
  equations.
\newblock Journal of Computational Physics \textbf{404}, 109088 (2020)

\bibitem{DFTBW2018}
Dumbser, M., Fambri, F., Tavelli, M., Bader, M., Weinzierl, T.: Efficient
  implementation of {ADER} discontinuous {G}alerkin schemes for a scalable
  hyperbolic pde engine.
\newblock Axioms \textbf{7}(3), 63 (2018)

\bibitem{ADERCCZ4}
Dumbser, M., Guercilena, F., K\"oppel, S., Rezzolla, L., Zanotti, O.:
  {Conformal and covariant Z4 formulation of the Einstein equations: strongly
  hyperbolic first--order reduction and solution with discontinuous Galerkin
  schemes}.
\newblock Physical Review D \textbf{97}, 084053 (2018)

\bibitem{AMR3DNC}
Dumbser, M., Hidalgo, A., Zanotti, O.: {High order space--time adaptive
  ADER--WENO finite volume schemes for non--conservative hyperbolic systems}.
\newblock Computer Methods in Applied Mechanics and Engineering \textbf{268},
  359--387 (2014)

\bibitem{DumbserKaeser06b}
Dumbser, M., K\"aser, M.: Arbitrary high order non-oscillatory finite volume
  schemes on unstructured meshes for linear hyperbolic systems.
\newblock Journal of Computational Physics \textbf{221}, 693--723 (2007)

\bibitem{DGLimiter3}
Dumbser, M., Loub{\`{e}}re, R.: {A simple robust and accurate a posteriori
  sub-cell finite volume limiter for the discontinuous Galerkin method on
  unstructured meshes}.
\newblock Journal of Computational Physics \textbf{319}, 163--199 (2016)

\bibitem{GPRmodel}
Dumbser, M., Peshkov, I., Romenski, E., Zanotti, O.: {High order ADER schemes
  for a unified first order hyperbolic formulation of continuum mechanics:
  Viscous heat-conducting fluids and elastic solids}.
\newblock Journal of Computational Physics \textbf{314}, 824--862 (2016)

\bibitem{GPRmodelMHD}
Dumbser, M., Peshkov, I., Romenski, E., Zanotti, O.: {H}igh order {ADER}
  schemes for a unified first order hyperbolic formulation of {N}ewtonian
  continuum mechanics coupled with electro-dynamics.
\newblock Journal of Computational Physics \textbf{348}, 298--342 (2017)

\bibitem{DumbserZanotti}
Dumbser, M., Zanotti, O.: Very high order {PNPM} schemes on unstructured meshes
  for the resistive relativistic {MHD} equations.
\newblock Journal of Computational Physics \textbf{228}, 6991--7006 (2009)

\bibitem{AMR3DCL}
Dumbser, M., Zanotti, O., Hidalgo, A., Balsara, D.: {ADER-WENO Finite Volume
  Schemes with Space-Time Adaptive Mesh Refinement}.
\newblock Journal of Computational Physics \textbf{248}, 257--286 (2013)

\bibitem{DGLimiter1}
Dumbser, M., Zanotti, O., Loub{\`{e}}re, R., Diot, S.: {A posteriori subcell
  limiting of the discontinuous Galerkin finite element method for hyperbolic
  conservation laws}.
\newblock Journal of Computational Physics \textbf{278}, 47--75 (2014)

\bibitem{HLLEM}
Einfeldt, B., Munz, C.D., Roe, P.L., Sj\"{o}green, B.: On godunov-type methods
  near low densities.
\newblock J. Comput. Phys. \textbf{92}(2), 273--295 (1991)

\bibitem{EM20}
Escalante, C., Morales, T.: A general non--hydrostatic hyperbolic formulation
  for {B}oussinesq dispersive shallow flows and its numerical approximation.
\newblock Journal of Scientific Computing \textbf{83}, 62 (2020)

\bibitem{ADERGRMHD}
Fambri, F., Dumbser, M., K\"oppel, S., Rezzolla, L., Zanotti, O.: {ADER
  discontinuous Galerkin schemes for general-relativistic ideal
  magnetohydrodynamics}.
\newblock Monthly Notices of the Royal Astronomical Society (MNRAS)
  \textbf{477}, 4543--4564 (2018)

\bibitem{ADERDGVisc}
Fambri, F., Dumbser, M., Zanotti, O.: Space-time adaptive {ADER}-{DG} schemes
  for dissipative flows: {C}ompressible {N}avier-{S}tokes and resistive {MHD}
  equations.
\newblock Computer Physics Communications \textbf{220}, 297--318 (2017)

\bibitem{favrie:2017}
Favrie, N., Gavrilyuk, S.: {A rapid numerical method for solving
  Serre-Green-Naghdi equations describing long free surface gravity waves.}
\newblock Nonlinearity \textbf{30}, 2718--2736 (2017)

\bibitem{gaburroReview}
Gaburro, E.: A unified framework for the solution of hyperbolic pde systems
  using high order direct arbitrary-lagrangian-eulerian schemes on moving
  unstructured meshes with topology change.
\newblock Archives of Computational Methods in Engineering  (2020).
\newblock \doi{10.1007/s11831-020-09411-7}

\bibitem{GaburroAREPO}
Gaburro, E., Boscheri, W., Chiocchetti, S., Klingenberg, C., Springel, V.,
  Dumbser, M.: High order direct {Arbitrary-Lagrangian-Eulerian} schemes on
  moving {V}oronoi meshes with topology changes.
\newblock Journal of Computational Physics \textbf{407}, 109167 (2020)

\bibitem{gaburro2018well}
Gaburro, E., Castro, M.J., Dumbser, M.: Well-balanced
  arbitrary-lagrangian-eulerian finite volume schemes on moving nonconforming
  meshes for the euler equations of gas dynamics with gravity.
\newblock Monthly Notices of the Royal Astronomical Society \textbf{477}(2),
  2251--2275 (2018)

\bibitem{gaburro2018diffuse}
Gaburro, E., Castro, M.J., Dumbser, M.: A well balanced diffuse interface
  method for complex nonhydrostatic free surface flows.
\newblock Computers \& Fluids \textbf{175}, 180--198 (2018)

\bibitem{GaburroNonConforming}
Gaburro, E., Dumbser, M., Castro, M.: {Direct Arbitrary-Lagrangian-Eulerian
  finite volume schemes on moving nonconforming unstructured meshes}.
\newblock Computers and Fluids \textbf{159}, 254--275 (2017)

\bibitem{godunov}
Godunov, S.: Finite difference methods for the computation of discontinuous
  solutions of the equations of fluid dynamics.
\newblock Mathematics of the USSR: Sbornik \textbf{47}, 271--306 (1959)

\bibitem{GodRom2003}
Godunov, S., Romenski, E.: Elements of continuum mechanics and conservation
  laws.
\newblock Kluwer Academic/Plenum Publishers (2003)

\bibitem{GodRom1972}
Godunov, S.K., Romenskii, E.I.: {Nonstationary equations of nonlinear
  elasticity theory in Eulerian coordinates}.
\newblock Journal of Applied Mechanics and Technical Physics \textbf{13}(6),
  868--884 (1972)

\bibitem{guermond2018second}
Guermond, J.L., Nazarov, M., Popov, B., Tomas, I.: Second-order invariant
  domain preserving approximation of the euler equations using convex limiting.
\newblock SIAM Journal on Scientific Computing \textbf{40}(5), A3211--A3239
  (2018)

\bibitem{luo3}
Halashi, B., Luo, H.: {A reconstructed discontinuous Galerkin method for
  magnetohydrodynamics on arbitrary grids}.
\newblock Journal of Computational Physics \textbf{326}, 258--277 (2016)

\bibitem{eno}
Harten, A., Engquist, B., Osher, S., Chakravarthy, S.: Uniformly high order
  essentially non-oscillatory schemes, {III}.
\newblock Journal of Computational Physics \textbf{71}, 231--303 (1987)

\bibitem{HartenENO}
Harten, A., Osher, S.: Uniformly high-order accurate nonoscillatory schemes
  {I}.
\newblock SIAM J. Num. Anal. \textbf{24}, 279--309 (1987)

\bibitem{HuShuVortex1999}
Hu, C., Shu, C.: A high-order weno finite difference scheme for the equations
  of ideal magnetohydrodynamics.
\newblock Journal of Computational Physics \textbf{150}, 561 -- 594 (1999)

\bibitem{ryan2}
Ji, L., Xu, Y., Ryan, J.: Accuracy enhancement of the linear
  convection--diffusion equation in multiple dimensions.
\newblock Mathematics of Computation \textbf{81}, 1929--1950 (2012)

\bibitem{JiangShu1996}
Jiang, G., Shu, C.: Efficient implementation of weighted {ENO} schemes.
\newblock Journal of Computational Physics \textbf{126}(1), 202--228 (1996)

\bibitem{SedovExact}
Kamm, J., Timmes, F.: On efficient generation of numerically robust sedov
  solutions.
\newblock Technical Report LA-UR-07-2849,LANL  (2007)

\bibitem{SolidBodies2020}
Kemm, F., Gaburro, E., Thein, F., Dumbser, M.: A simple diffuse interface
  approach for compressible flows around moving solids of arbitrary shape based
  on a reduced baer-nunziato model.
\newblock Computers \& Fluids \textbf{204}, 104536 (2020)

\bibitem{Khokhlov1998}
Khokhlov, A.: Fully threaded tree algorithms for adaptive refinement fluid
  dynamics simulations.
\newblock Journal of Computational Physics \textbf{143}(2), 519 -- 543 (1998)

\bibitem{ryan4}
King, J., Mirzaee, H., Ryan, J., Kirby, R.: Smoothness--increasing
  accuracy--conserving {SIAC)} filtering for discontinuous {Galerkin}
  solutions: improved errors versus higher-order accuracy.
\newblock Journal of Scientific Computing \textbf{53}, 129--149 (2012)

\bibitem{Komissarov1997}
{Komissarov}, S.S.: {On the properties of Alfv{\'e}n waves in relativistic
  magnetohydrodynamics}.
\newblock Physics Letters A \textbf{232}, 435--442 (1997)

\bibitem{Komissarov1999}
{Komissarov}, S.S.: {A Godunov-type scheme for relativistic
  magnetohydrodynamics}.
\newblock Mon. Not. R. Astron. Soc. \textbf{303}, 343--366 (1999)

\bibitem{Kouveliotou1993}
{Kouveliotou}, C., {Meegan}, C.A., {Fishman}, G.J., {Bhat}, N.P., {Briggs},
  M.S., {Koshut}, T.M., {Paciesas}, W.S., {Pendleton}, G.N.: {Identification of
  two classes of gamma-ray bursts}.
\newblock Astrophys. J. \textbf{413}, L101--L104 (1993)

\bibitem{lax2}
Lax, P.: Weak solutions of nonlinear hyperbolic equations and their numerical
  approximation.
\newblock Comm. Pure Appl. Math. \textbf{7}, 159--193 (1954)

\bibitem{vanLeerDGdiffusion}
van Leer, B., Nomura, S.: Discontinuous {Galerkin} for diffusion.
\newblock In: Proceedings of 17th AIAA Computational Fluid Dynamics Conference
  (June 6--9 2005), AIAA-2005-5108 (2005)

\bibitem{Leismann2005}
{Leismann}, T., {Ant{\'o}n}, L., {Aloy}, M.A., {M{\"u}ller}, E.,
  {Mart{\'{\i}}}, J.M., {Miralles}, J.A., {Ib{\'a}{\~n}ez}, J.M.: {Relativistic
  MHD simulations of extragalactic jets}.
\newblock Astronomy and Astrophyiscs \textbf{436}, 503--526 (2005)

\bibitem{Lohner1987}
L\"ohner, R.: {An adaptive finite element scheme for transient problems in
  CFD}.
\newblock Computer Methods in Applied Mechanics and Engineering \textbf{61},
  323--338 (1987)

\bibitem{ADERMOOD}
Loubere, R., Dumbser, M., Diot, S.: A new family of high order unstructured
  mood and ader finite volume schemes for multidimensional systems of
  hyperbolic conservation laws.
\newblock Communications in Computational Physics \textbf{16}(3), 718--763
  (2014)

\bibitem{LoubereSedov3D}
Loub\`ere, R., Maire, P., V\'achal, P.: {3D staggered Lagrangian hydrodynamics
  scheme with cell-centered Riemann solver-based artificial viscosity.}
\newblock International Journal for Numerical Methods in Fluids \textbf{72}, 22
  -- 42 (2013)

\bibitem{luo1}
Luo, H., Luo, L., Nourgaliev, R., Mousseau, V., Dinh, N.: {A reconstructed
  discontinuous Galerkin method for the compressible Navier--Stokes equations
  on arbitrary grids}.
\newblock Journal of Computational Physics \textbf{229}, 6961--6978 (2010)

\bibitem{luo2}
Luo, H., Xia, Y., Spiegel, S., Nourgaliev, R., Jiang, Z.: {A reconstructed
  discontinuous Galerkin method based on a Hierarchical WENO reconstruction for
  compressible flows on tetrahedral grids }.
\newblock Journal of Computational Physics \textbf{236}, 477--492 (2013)

\bibitem{Michel1991}
{Michel}, F.C.: {Theory of neutron star magnetospheres} (1991)

\bibitem{ryan5}
Mirzaee, H., King, J., Ryan, J., Kirby, R.: Smoothness--increasing
  accuracy--conserving {(SIAC)} filters for discontinuous {Galerkin} solutions
  over unstructured triangular meshes.
\newblock SIAM Journal on Scientific Computing \textbf{35}, A212--A230 (2013)

\bibitem{MunzCleaning}
Munz, C., Omnes, P., Schneider, R., Sonnendr\"ucker, E., Voss, U.: {Divergence
  Correction Techniques for Maxwell Solvers Based on a Hyperbolic Model}.
\newblock Journal of Computational Physics \textbf{161}, 484--511 (2000)

\bibitem{OrszagTang}
Orszag, S.A., Tang, C.M.: Small--scale structure of two--dimensional
  magnetohydrodynamic turbulence.
\newblock Journal of Fluid Mechanics \textbf{90}, 129 (1979)

\bibitem{PeshRom2014}
Peshkov, I., Romenski, E.: A hyperbolic model for viscous {{N}ewtonian} flows.
\newblock Continuum Mechanics and Thermodynamics \textbf{28}, 85--104 (2016)

\bibitem{PiconeDahlburg}
Picone, J.M., Dahlburg, R.B.: Evolution of the {Orszag--Tang} vortex system in
  a compressible medium. {II.} supersonic flow.
\newblock Phys. Fluids B \textbf{3}, 29--44 (1991)

\bibitem{rannabauer2018ader}
Rannabauer, L., Dumbser, M., Bader, M.: Ader-dg with a-posteriori finite-volume
  limiting to simulate tsunamis in a parallel adaptive mesh refinement
  framework.
\newblock Computers \& Fluids \textbf{173}, 299--306 (2018)

\bibitem{Exahype}
Reinarz, A., et~al.: Exahype: an engine for parallel dynamically adaptive
  simulations of wave problems.
\newblock Computer Physics Communications p. 107251 (2020)

\bibitem{DeLaRosaMunzDGMHD}
de~la Rosa, J.N., Munz, C.D.: Hybrid dg/fv schemes for magnetohydrodynamics and
  relativistic hydrodynamics.
\newblock Computer Physics Communications \textbf{222}, 113 -- 135 (2018)

\bibitem{ryan3}
Ryan, J., Cockburn, B.: Local derivative post-processing for the discontinuous
  {Galerkin} method.
\newblock Journal of Computational Physics \textbf{228}, 8642--8664 (2009)

\bibitem{ryan1}
Ryan, J., Shu, C., Atkins, H.: Extension of a post-processing technique for the
  discontinuous {Galerkin} method for hyperbolic equations with applications to
  an aeroacoustic problem.
\newblock SIAM Journal on Scientific Computing \textbf{26}, 821--843 (2005)

\bibitem{Sedov}
Sedov, L.: Similarity and Dimensional Methods in Mechanics.
\newblock Academic Press, New York (1959)

\bibitem{shi2003resolution}
Shi, J., Zhang, Y.T., Shu, C.W.: Resolution of high order weno schemes for
  complicated flow structures.
\newblock Journal of Computational Physics \textbf{186}(2), 690--696 (2003)

\bibitem{shu1}
Shu, C.: Essentially non-oscillatory and weighted essentially non-oscillatory
  schemes for hyperbolic {Conservation} {Laws}.
\newblock NASA/CR-97-206253 ICASE Report No.97-65  (1997)

\bibitem{shuosher1}
Shu, C., Osher, S.: Efficient implementation of essentially non-oscillatory
  shock capturing schemes.
\newblock Journal of Computational Physics \textbf{77}, 439--471 (1988)

\bibitem{shu2016high}
Shu, C.W.: High order {WENO} and {DG} methods for time-dependent
  convection-dominated {PDEs}: A brief survey of several recent developments.
\newblock Journal of Computational Physics \textbf{316}, 598--613 (2016)

\bibitem{SonntagDG}
Sonntag, M., Munz, C.: Shock capturing for discontinuous galerkin methods using
  finite volume subcells.
\newblock In: J.~Fuhrmann, M.~Ohlberger, C.~Rohde (eds.) Finite Volumes for
  Complex Applications VII, pp. 945--953. Springer (2014)

\bibitem{stroud}
Stroud, A.: Approximate Calculation of Multiple Integrals.
\newblock Prentice-Hall Inc., Englewood Cliffs, New Jersey (1971)

\bibitem{tavelli2014high}
Tavelli, M., Dumbser, M.: A high order semi-implicit discontinuous galerkin
  method for the two dimensional shallow water equations on staggered
  unstructured meshes.
\newblock Applied Mathematics and Computation \textbf{234}, 623--644 (2014)

\bibitem{toro3}
Titarev, V., Toro, E.: {ADER}: Arbitrary high order {Godunov} approach.
\newblock Journal of Scientific Computing \textbf{17}(1-4), 609--618 (2002)

\bibitem{titarevtoro}
Titarev, V., Toro, E.: {ADER} schemes for three-dimensional nonlinear
  hyperbolic systems.
\newblock Journal of Computational Physics \textbf{204}, 715--736 (2005)

\bibitem{toro-book}
Toro, E.: {Riemann} Solvers and Numerical Methods for Fluid Dynamics, second
  edn.
\newblock Springer (1999)

\bibitem{ToroBook}
Toro, E.: Riemann Solvers and Numerical Methods for Fluid Dynamics: a Practical
  Introduction.
\newblock Springer (2009)

\bibitem{toro4}
Toro, E., Titarev, V.: Solution of the generalized {Riemann} problem for
  advection-reaction equations.
\newblock Proc. Roy. Soc. London pp. 271--281 (2002)

\bibitem{Toro:2006a}
Toro, E.F., Titarev, V.A.: {Derivative Riemann solvers for systems of
  conservation laws and ADER methods}.
\newblock Journal of Computational Physics \textbf{212}(1), 150--165 (2006)

\bibitem{Toth2000}
Toth, G.: The div b=0 constraint in shock-capturing magnetohydrodynamics codes.
\newblock Journal of Computational Physics \textbf{161}(2), 605 -- 652 (2000)

\bibitem{luo6}
Wang, C., Cheng, J., Berndt, M., Carlson, N., Luo, H.: {Application of
  nonlinear Krylov acceleration to a reconstructed discontinuous Galerkin
  method for compressible flows}.
\newblock Computers and Fluids \textbf{163}, 32--49 (2018)

\bibitem{luo4}
Wang, C., Luo, H., Shashkov, M.: {A reconstructed discontinuous Galerkin method
  for compressible flows in Lagrangian formulation}.
\newblock Computers and Fluids \textbf{202}, 104522 (2020)

\bibitem{wang2020reconstructed}
Wang, C., Luo, H., Shashkov, M.: A reconstructed discontinuous galerkin method
  for compressible flows in lagrangian formulation.
\newblock Computers \& Fluids p. 104522 (2020)

\bibitem{luo5}
Wang, X., Cheng, C., Luo, H., Zhao, Q.: {A reconstructed direct discontinuous
  Galerkin method for simulating the compressible laminar and turbulent flows
  on hybrid grids}.
\newblock Computers and Fluids \textbf{168}, 216--231 (2018)

\bibitem{Peano2}
Weinzierl, T., Mehl, M.: {Peano-A traversal and storage scheme for octree-like
  adaptive Cartesian multiscale grids}.
\newblock SIAM Journal on Scientific Computing \textbf{33}, 2732--2760 (2011)

\bibitem{woodwardcol84}
Woodward, P., Colella, P.: The numerical simulation of two-dimensional fluid
  flow with strong shocks.
\newblock Journal of Computational Physics \textbf{54}, 115--173 (1984)

\bibitem{RMHD}
Zanna, L.D., Bucciantini, N., Londrillo, P.: An efficient shock-capturing
  central-type scheme for multidimensional relativistic flows {II}.
  magnetohydrodynamics.
\newblock Astronomy and Astrophysics \textbf{400}, 397--413 (2003)

\bibitem{Zanotti2015d}
Zanotti, O., Fambri, F., Dumbser, M.: Solving the relativistic
  magnetohydrodynamics equations with {ADER} discontinuous {G}alerkin methods,
  a posteriori subcell limiting and adaptive mesh refinement.
\newblock Mon. Not. R. Astron. Soc. \textbf{452}, 3010--3029 ({2015})

\bibitem{DGLimiter2}
Zanotti, O., Fambri, F., Dumbser, M., Hidalgo, A.: {Space--time adaptive ADER
  discontinuous Galerkin finite element schemes with a posteriori sub--cell
  finite volume limiting}.
\newblock Computers and Fluids \textbf{118}, 204--224 (2015)

\bibitem{ZhangShu3D}
Zhang, Y., Shu, C.: Third order {WENO} scheme on three dimensional tetrahedral
  meshes.
\newblock Communications in Computational Physics \textbf{5}, 836--848 (2009)

\end{thebibliography}

\end{document}